\RequirePackage{ifthen}
\newboolean{MOR}
\setboolean{MOR}{false}

\ifthenelse{\boolean{MOR}}
{
\documentclass[ijds,nonblindrev]{informs4}
} 
{
\documentclass[10pt]{article}
\usepackage[margin=1in]{geometry}
}

\usepackage{amsmath}
\usepackage{graphicx}
\usepackage{psfrag}
\usepackage{amssymb}
\usepackage{amsmath}
\usepackage{verbatim,booktabs}
\usepackage{subfigure,epsfig,url}
\usepackage{float}
\usepackage{mathrsfs}

\usepackage{enumitem} 

\usepackage[usenames,dvipsnames]{pstricks}
\usepackage{epsfig}
\usepackage{pst-grad} 
\usepackage{pst-plot} 

\ifthenelse{\boolean{MOR}}
{
\usepackage{amsfonts}

\usepackage{natbib}
 \NatBibNumeric
 \def\newblock{\ }%
 \bibpunct[, ]{[}{]}{,}{n}{}{,}%

\usepackage{rotating}
\usepackage{fancyvrb}

\TheoremsNumberedThrough     
\ECRepeatTheorems
\JOURNAL{Mathematics of Operations Research}

\EquationsNumberedThrough    

\MANUSCRIPTNO{MOR-2022-201}

\newenvironment{prf}[1][]
{\proof{Proof.}}
{\hfill\ensuremath{\Halmos} \endproof}

\newenvironment{prfc}[1][]
{\proof{Proof #1.}}
{\hfill\ensuremath{\Halmos} \endproof}
}
{
\usepackage{hyperref}
\usepackage{amsthm}

\newtheorem{theorem}{Theorem}
\newtheorem{lemma}{Lemma}
\newtheorem{proposition}{Proposition}
\newtheorem{corollary}{Corollary}
\newtheorem{remark}{Remark}
\newtheorem{claim}{Claim}

\newenvironment{prf}[1][]
{\begin{proof}}
{\end{proof}}

\newenvironment{prfc}[1][]
{\begin{proof}[Proof #1]}
{\end{proof}}
}

\usepackage[capitalize,noabbrev]{cleveref}

\newenvironment{cpf}
{\begin{trivlist} \item[] {\em Proof of claim. }}
{$\hfill\diamond$ \end{trivlist}}

\newenvironment{spf}
{\begin{trivlist} \item[] {\em Proof of step. }}
{$\hfill\triangleleft$ \end{trivlist}}

\newtheorem{step}{Step}

\crefname{claim}{Claim}{Claims}

\usepackage{xspace}

\newcommand{\card}[1]{\lvert#1\rvert}

\newcommand{\ceil}[1]{\lceil#1\rceil}
\newcommand{\pare}[1]{\left(#1\right)}

\def\st{{\rm s.t.}}

\DeclareMathOperator{\prob}{\mathbb P}
\DeclareMathOperator{\avg}{\mathbb E}
\DeclareMathOperator{\obj}{obj}

\newcommand{\pr}{p}
\newcommand{\R}{\mathbb R}

\allowdisplaybreaks

\def\ie{{i.e.,} }

\def\G{{\mathcal G}}
\def\S{{\mathcal S}}

\def\N{{\mathcal N}}

\def\T{{\mathcal T}}

\def\M{{\mathcal M}}
\def\P{{\mathcal P}}
\def\W{{\mathcal W}}
\def\F{{\mathcal F}}
\def\01{\ensuremath{0\mathord{-}1}}

\newcommand{\ground}{\bar {\mathcal W}}
\newcommand{\noisy}{\mathcal G}
\newcommand{\adversary}{\mathcal G'}

\def\MP{\text{MP}}

\newcommand{\ttl}{Rank-one Boolean tensor factorization and the multilinear polytope}

\newcommand{\bstrct}{We consider the NP-hard problem of
finding the closest rank-one binary tensor to a given binary tensor, which we refer to as the rank-one Boolean tensor factorization (BTF) problem.
This optimization problem can be used to recover a planted rank-one tensor from noisy observations.
We formulate rank-one BTF as the problem of minimizing a linear function over a highly structured multilinear set. Leveraging on our prior results regarding the facial structure of multilinear polytopes, we propose novel linear programming relaxations for rank-one BTF. 
We then establish deterministic sufficient conditions under which our proposed linear programs recover a planted rank-one tensor.
To analyze the effectiveness of these deterministic conditions, we consider a semi-random model for the noisy tensor, and obtain high probability recovery guarantees for the linear programs.
Our theoretical results as well as numerical simulations indicate that certain facets of the multilinear polytope significantly improve the recovery properties of linear programming relaxations for rank-one BTF.}

\newcommand{\kywrds}{Rank-one Boolean tensor factorization, multilinear polytope, linear programming relaxation, recovery guarantee, semi-random models}

\newcommand{\MSCcodes}{Primary: 90C10; secondary: 90C26, 90C57}

\ifthenelse{\boolean{MOR}}
{
\TITLE{\ttl}
\RUNTITLE{\ttl}

\ARTICLEAUTHORS{
\AUTHOR{Alberto Del Pia}
\AFF{Department of Industrial and Systems Engineering \& Wisconsin Institute for Discovery, University of Wisconsin-Madison, \EMAIL{delpia@wisc.edu}}
\AUTHOR{Aida Khajavirad}
\AFF{Department of Industrial \& Systems Engineering, Lehigh University, \EMAIL{aida@lehigh.edu}}}
\RUNAUTHOR{Del Pia and Khajavirad}

\ABSTRACT{\bstrct}
\MSCCLASS{\MSCcodes}
}
{
\title{\ttl}

\author{Alberto Del Pia
\thanks{Department of Industrial and Systems Engineering \& Wisconsin Institute for Discovery,
University of Wisconsin-Madison.
E-mail: {\tt delpia@wisc.edu}}
\and
Aida Khajavirad
\thanks{Department of Industrial and Systems Engineering,
Lehigh University.
E-mail: {\tt aida@lehigh.edu}}}
}

\date{May 2, 2024}

\begin{document}

\maketitle

\ifthenelse{\boolean{MOR}}
{}
{
\begin{abstract}
\bstrct
\end{abstract}



\bigskip

\noindent
\hangindent=2cm
\emph{Key words:}
\kywrds

\bigskip

\noindent
\emph{Mathematics Subject Classification:} 
\MSCcodes
}

\section{Introduction}
A \emph{tensor} of order $N$ is an $N$-dimensional array.
Factorizations of high-order tensors, i.e., $N \geq 3$, as products of low-rank matrices, have applications in signal processing, numerical linear algebra, computer vision, data mining, neuroscience, and elsewhere~\cite{KolBad09,Hore16,Zhou13,Bi13}.
We consider the problem of factorizing a high-order tensor with binary entries, henceforth, referred to as a \emph{binary tensor}. Such problems arise in applications such as neuro-imaging, recommendation systems, topic modeling, and sensor network localization~\cite{Wang19,Maz14,Rai15,Hong20}. 
In \emph{Boolean tensor factorization} (BTF), the binary tensor is approximated by products of low rank \emph{binary} matrices using Boolean algebra~\cite{Miet11}. 
BTF is a very useful tool for analyzing binary tensors to discover latent factors from them~\cite{ErdMie13,Park17,RukHolYau18,Wan20}.
Furthermore, BTF produces more interpretable and sparser results than normal factorization methods~\cite{Miet11}.
BTF is NP-hard in general~\cite{GilVav18};  all existing methods to tackle this problem rely on heuristics and hence do not provide any guarantee on the quality of the solution~\cite{Miet11,ErdMie13,beloh13,Park17,RukHolYau18}.

In order to formally define BTF, we first introduce some notation. For an integer $n$, we denote by $[n]:=\{1,2,\dots,n\}$.
All the tensors that we consider in this work have order three, meaning that each element of the tensor has three indices.
Given a tensor $\W$, we denote its $(i,j,k)$th element by $w_{ijk}$.
We denote by $\otimes$ the vector outer product.
That is, if $x \in \R^n$, $y \in \R^m$, $z \in \R^l$, then $\W = x \otimes y \otimes z$ is a $n \times m \times l$ tensor defined by $w_{ijk} = x_i y_j z_k$, for $i \in [n]$, $j \in [m]$, $k \in [l]$.
The \emph{Frobenius norm} of a tensor $\W$, is defined as $\Vert\W\Vert:=\sqrt{\sum_{i,j,k}{w^2_{ijk}}}$.
The \emph{rank} (or \emph{Boolean rank}) of a binary tensor $\W$ is the smallest integer $r$ such that there exist $3r$ binary vectors $x^t,y^t,z^t$, for $t \in [r]$, with
$$
\W = \bigvee^r_{t=1} (x^t \otimes y^t \otimes z^t),
$$
where $\vee$ denotes the component-wise ``or'' operation.
In particular, a binary tensor has rank one if it is the outer product of three binary vectors.
Computing the Boolean rank of a binary tensor is NP-complete~\cite{Hastad90,Miet11}.
Interestingly, unlike matrices, there exist $n \times m \times l$ binary tensors whose Boolean rank is larger than $\max\{n,m,l\}$.
Indeed, a tight upper bound on the Boolean rank of a binary tensor is given by $\min\{nm,nl,ml\}$~\cite{Miet11}.

\subsection{Problem statement}
The \emph{rank-$r$ BTF} is the problem of finding the closest rank-$r$ binary tensor to a given binary tensor.
Precisely, we are given a $n \times m \times l$ binary tensor $\G=(g_{ijk})$ and an integer $r$, and we seek $3r$ binary vectors $x^t \in \{0,1\}^n$, $y^t \in \{0,1\}^m$, $z^t \in \{0,1\}^l$, for all $t \in [r]$, that minimize
\begin{align*}
\Big\lVert\G - \bigvee^r_{t=1} x^t \otimes y^t \otimes z^t\Big\rVert^2.
\end{align*}
In this paper, we focus on the simplest case of BTF; namely, the case with $r=1$, referred to as the~\emph{rank-one BTF}.
This problem can be formulated as the following optimization problem:
\begin{align}
\label{pr comb}
\tag{P}
\min \quad & \Big\Vert\G - x \otimes y \otimes z\Big\Vert^2 \\
\st \quad & x \in \{0,1\}^n, \ y \in \{0,1\}^m, \ z \in \{0,1\}^l. \nonumber
\end{align}
Rank-one BTF is NP-hard in general~\cite{Miet11} and to this date no algorithm with theoretical guarantees is known for this problem. In this paper, we introduce novel linear programming (LP) relaxations with theoretical performance guarantees for rank-one BTF. To this end, in the following, we present an equivalent integer programming reformulation of Problem~\eqref{pr comb} in an extended space.
Define
\begin{equation}\label{s1s2}
S_0: = \{(i,j,k) \in [n]\times[m]\times[l]: g_{ijk}=0\}, \quad
S_1: = \{(i,j,k) \in [n]\times[m]\times[l]: g_{ijk}=1\}.
\end{equation}
Since $\G, x, y, z$ are all binary valued, the objective function of Problem~\eqref{pr comb} can be written as
\begin{align*}
\Big\Vert\G - x \otimes y \otimes z\Big\Vert^2
=
 \sum_{(i,j,k) \in S_0}{x_i y_j z_k}+
\sum_{(i,j,k) \in S_1}{(1-x_i y_j z_k)}.
\end{align*}
Subsequently, we introduce auxiliary variables $w_{ijk} := x_i y_j z_k$, for $i\in [n]$, $j\in [m]$, $k\in [l]$.
It then follows that rank-one BTF can be equivalently written, in an extended space, as the problem of minimizing a linear function over a highly structured~\emph{multilinear set}~\cite{dPKha16}:
\begin{align}
\label{tfr1}
\tag{extP}
\min \quad & \sum_{(i,j,k) \in S_0}{w_{ijk}}+\sum_{(i,j,k) \in S_1}{(1-w_{ijk})} \\
\st  \quad & w_{ijk} = x_i y_j z_k, \qquad \forall i \in [n], j \in [m], k\in [l] \nonumber\\
& x \in \{0,1\}^n, \ y \in \{0,1\}^m, \ z \in \{0,1\}^l. \nonumber
\end{align}
The above reformulation enables us to leverage on our previous results regarding the facial structure of the convex hull of multilinear sets~\cite{dPKha16,dPIda16,dPIda18,dPIda19,dPIdaSah20,dPKha23MPA}, and develop strong LP relaxations for rank-one BTF.

\subsection{Recovery of a planted model} It is widely accepted that worst-case guarantees for algorithms are often too pessimistic, as the input data in most real-world applications is highly structured.
Motivated by this observation, a recent stream of research in mathematical data science is focused on obtaining theoretical guarantees for various existing and new algorithms under stochastic models that better reflect \emph{typical} problem instances (see for example~\cite{AbbBanHal16,MixVilWar17,LinStr19,dPIdaTim20,AntoAida20,AidaTonio21}). More specifically, in our context, we assume that the input tensor $\noisy$ is obtained by corrupting a binary rank-one tensor, referred to as the \emph{ground truth}. Our goal is to obtain sufficient conditions under which the solution returned by the algorithm corresponds to the ground truth. Such conditions are often referred to as~\emph{recovery guarantees}~\cite{AbbBanHal16,MixVilWar17,LinStr19,dPIdaTim20,AntoAida20,AidaTonio21}.

In this paper, we show that given a corrupted tensor $\noisy$, Problem~\eqref{pr comb} coincides with the maximum likelihood estimator for recovering the ground truth (see~\cref{MLEMAP}).
Consider an optimization problem for rank-one BTF.
We say that this optimization problem \emph{recovers the ground truth}, if the optimization problem has a unique optimal solution and it corresponds to the ground truth rank-one tensor.
We are interested in addressing the following question: under a suitable generative model for the input, what is the maximum level of corruption, under which our optimization problems recover the ground truth with high probability?
Throughout this paper, when we write that an event happens \emph{with high probability}, we mean that the event happens with probability that goes to one, as $n,m,l \to \infty$.
Our analysis is based on an important assumption that the optimization problem only receives the tensor $\G$ as input and does not have the knowledge of how it was generated, hence it is, in particular, a ``parameter-free'' algorithm.


\subsection{Fully random versus semi-random models}
We now define two stochastic models for rank-one BTF.
First, we introduce the \emph{fully-random corruption model}.
Consider binary vectors $\bar x \in \{0,1\}^n$, $\bar y \in \{0,1\}^m$, $\bar z \in \{0,1\}^l$ and define the \emph{ground truth} rank-one tensor $\bar \W = (\bar w_{ijk}) := \bar x \otimes \bar y \otimes \bar z  \in \{0,1\}^{n \times m \times l}$.
Given $\pr \in [0,1]$, the noisy tensor $\G$ is constructed as follows: for each $(i,j,k) \in [n] \times [m] \times [l]$, $g_{ijk}$ is corrupted with probability $\pr$, \ie $g_{ijk} := 1- \bar x_i \bar y_j \bar z_k$, and $g_{ijk}$ is not corrupted with probability $1-\pr$, \ie $g_{ijk} := \bar x_i \bar y_j \bar z_k$.
This model has been used in~\cite{Miet11,Park17,RukHolYau18} for conducting empirical analysis of some heuristics for BTF.

It is well-understood that many popular algorithms, including spectral methods, crucially rely on the fully-random model and on the knowledge of the error probability parameter $p$ to succeed
(see for example~\cite{feige01}). 
To address this shortcoming, \emph{semi-random models} have been introduced~\cite{feige01}; these  
models are generated by combining random and adversarial components. 
Semi-random models are often used to measure the ``robustness'' of algorithms as they better reflect real-world instances than fully-random models. 
It has been shown that, unlike some popular heuristics such as spectral methods that fail under semi-random models, optimization-based techniques such as semi-definite programming continue to work well under these models~\cite{feige01,moiperwei16,ricJavMon16}. In this paper, we prove that for rank-one BTF, all ``reasonable'' optimization algorithms share this robustness property (see~\cref{robustLPs}).

We now define the \emph{semi-random corruption model} for rank-one BTF.
In this model, first the noisy tensor $\noisy$ is constructed from the ground truth rank-one tensor $\ground$ according to the fully-random model.
Then, an adversary modifies $\noisy$ by applying an arbitrary sequence of~\emph{monotone transformations} defined as follows:
\begin{itemize}
\item
If an entry is 1 in $\ground$ and is 0 in $\noisy$, the adversary may revert to 1 the entry in $\noisy$;
\item
If an entry is 0 in $\ground$ and is 1 in $\noisy$, the adversary may revert to 0 the entry in $\noisy$.
\end{itemize}
It is important to note that while it seems the adversary makes helpful changes to the noisy tensor which should make the problem easier to solve, such changes break the special properties associated with the fully-random model and hence heuristics that highly rely on such unrealistic properties do not succeed in semi-random models. 


Denote by $r_{\bar x}$ (resp. $r_{\bar y}$, $r_{\bar z}$)
the ratio of ones in $\bar x$ (resp. $\bar y$, $\bar z$) to the number of elements in $\bar x$ (resp. $\bar y$, $\bar z$).
In this paper, we mainly focus on the case where $\pr$ is a constant, and we often consider the ``dense regime,'' in which $r_{\bar x}$, $r_{\bar y}$, and $r_{\bar z}$ are positive constants.

\subsection{Our contribution} In this paper, we introduce novel LP relaxations for rank-one BTF.
We investigate the recovery properties of the proposed LPs under the semi-random corruption model.
Clearly, any convex relaxation can recover the ground truth only if the original NP-hard problem succeeds in doing so.
We start by establishing the recovery threshold for rank-one BTF under our two corruption models.
Namely, we obtain necessary and sufficient conditions under which Problem~\eqref{pr comb} recovers the ground truth with high probability. In particular, our results imply that, under mild assumptions on the growth rate of $n,m,l$, Problem~\eqref{pr comb} recovers the ground truth with high probability if and only if $\pr < 1/2$ (see Theorems~\ref{th it lb} and~\ref{th it ub}).

We then study the recovery properties of our proposed LPs. To this end, we first establish deterministic sufficient conditions under which the ground truth tensor is the unique optimal solution of each LP relaxation. 
Subsequently, we obtain high probability recovery guarantees under the semi-random corruption model. 
We start by considering a simple LP relaxation for rank-one BTF, referred to as the \emph{standard LP}, and obtain a recovery guarantee for it. Our result in particular implies that if $r_{\bar x} = r_{\bar y} = r_{\bar z}$,  under mild assumptions on the growth rate of $n,m,l$, the standard LP recovers the ground truth with high probability if $\pr < \frac{r_{\bar w}}{2(1+r_{\bar w})}$, where we define  $r_{\bar w} := r_{\bar x} r_{\bar y} r_{\bar z}$ (see Theorem~\ref{th LP1}).
We refer to our strongest proposed LP relaxation as the \emph{complete LP}.
Since the theoretical analysis of the complete LP is rather complex, we consider a relaxation of this LP, which we refer to as the \emph{flower LP}.
Roughly speaking, this intermediate LP is obtained by adding flower inequalities~\cite{dPIda18} to the standard LP.
We prove that under mild assumptions on the growth rate of $n,m,l$, flower LP recovers the ground truth with high probability if $\pr < \frac{r_{\bar w}}{1+2 r_{\bar w}}$ (see Theorem~\ref{th LP2}). That is, utilizing a stronger LP relaxation results in an improvement of up to $33\%$ in the recovery threshold. Numerical experiments suggest that our recovery guarantees are fairly tight, and that the complete LP significantly outperforms the flower LP.
We believe that obtaining recovery guarantees for the complete LP is an interesting open question.
We remark that all proposed LP relaxations can be solved efficiently both in theory (\ie in polynomial time) and in practice.

\paragraph{Outline.} In \cref{sec basics}, we establish some basic properties of our optimization framework. In \cref{sec LP4BTF}, we introduce our LP relaxations for rank-one BTF. In \cref{sec main results}, we present the statements of our recovery guarantees. Preliminary numerical results are provided in Section~\ref{sec numerics}.
In \cref{sec IP recovery}, we prove our necessary and sufficient recovery conditions for rank-one BTF.
In \cref{sec SL recovery,sec FR recovery}, we prove our recovery guarantees for the standard LP and the flower LP, respectively.
Finally, the proof of a technical result omitted from \cref{sec LP4BTF} is given in \cref{sec facets}.

\section{Optimization and recovery}
\label{sec basics}
In this section, we establish two fundamental properties of optimization algorithms that convey the effectiveness of our approach for recovering a planted tensor from noisy observations.

\subsection{Maximum likelihood and maximum a posteriori estimators}


In the following, we present the connections between Problem~\eqref{pr comb} and the maximum likelihood (ML) estimator  and the maximum a posteriori (MAP) estimator for the recovery problem of a planted tensor.
By definition, the MAP estimator maximizes the
probability of recovering the planted rank-one tensor. Hence, if the MAP estimator fails in recovering the ground truth, no algorithm, efficient or not, will succeed in doing so.

\begin{proposition}\label{MLEMAP}
Consider the fully-random corruption model with $p < \frac{1}{2}$.
Then solving Problem~\eqref{pr comb} is equivalent to finding the ML estimator.
Furthermore, if we assume that all ground truth rank-one binary tensors are equally likely a priori, solving Problem~\eqref{pr comb} is equivalent to finding the MAP estimator.
\end{proposition}

\begin{prf}
Let $\W$ be the ground truth binary tensor and let $\G$ be the observed noisy binary tensor.
First, we show that solving Problem~\eqref{pr comb} is equivalent to finding the ML estimator.
We use the notation $\prob[\G \mid \W]$ to denote the probability that the noisy binary tensor $\G$ was observed given that $\W$ was the ground truth rank-one binary tensor.
Note that we have
\begin{align*}
& \prob\left[g_{ijk}=1 \mid w_{ijk}=0\right]=\prob\left[g_{ijk}=0 \mid w_{ijk}=1\right] = p \\
& \prob\left[g_{ijk}=1 \mid w_{ijk}=1\right]=\prob\left[g_{ijk}=0 \mid w_{ijk}=0\right] = 1-p,
\end{align*}
thus
\begin{align*}
\prob\left[g_{ijk} \mid w_{ijk}\right]=
\begin{cases}
p^{w_{ijk}} (1-p)^{1-w_{ijk}} & \text{if $g_{ijk}=0$} \\
p^{1-w_{ijk}} (1-p)^{w_{ijk}} & \text{if $g_{ijk}=1$}.
\end{cases}
\end{align*}
Given the noisy binary tensor $\G$, we seek a rank-one binary tensor $\W$ that maximizes the likelihood function
$$
\prob[\G \mid \W].
$$
Since the probability of corruption of the entries of $\W$ are independent, we have
\begin{align*}
\prob[\G \mid \W]
& =\prod_{(i,j,k) \in [n]\times [m]\times [l]} \prob\left[g_{ijk} \mid w_{ijk}\right] \\
& =\prod_{(i,j,k) \in S_0} \pare{p^{w_{ijk}} (1-p)^{1-w_{ijk}}} \prod_{(i,j,k) \in S_1} \pare{p^{1-w_{ijk}} (1-p)^{w_{ijk}}}.
\end{align*}
Maximizing the likelihood function is equivalent to maximizing the log-likelihood function:
\begin{align*}
\log \prob[\G \mid \W]
& = \sum_{(i,j,k) \in S_0} \pare{w_{ijk} \log p + (1-w_{ijk}) \log (1-p)} 
+ \sum_{(i,j,k) \in S_1} \pare{(1-w_{ijk}) \log p + w_{ijk} \log (1-p)} \\
& = \sum_{(i,j,k) \in S_0} \pare{w_{ijk} \log \frac{p}{1-p} + \log (1-p)} 
- \sum_{(i,j,k) \in S_1} \pare{w_{ijk} \log \frac{p}{1-p} - \log p} \\
& = \pare{\sum_{(i,j,k) \in S_0} w_{ijk} - \sum_{(i,j,k) \in S_1} w_{ijk}} \log \frac{p}{1-p} + |S_0|\log (1-p)  + |S_1|\log p.
\end{align*}
We then observe that the objective function of Problem~\eqref{pr comb} is obtained from the log-likelihood function via translation and scaling by $\log \frac{p}{1-p}$.
This scaling factor is negative since $p < 1/2$, which is desired since in Problem~\eqref{pr comb} we minimize the objective function.
This concludes the proof that solving Problem~\eqref{pr comb} is equivalent to finding ML estimator.

If we assume that all ground truth rank-one binary tensors are equally likely a priori, then by Bayes' rule, 
$$
\arg\max \prob[\W \mid \G] = \arg\max \prob[\G \mid \W],
$$
thus solving Problem~\eqref{pr comb} is equivalent to finding the MAP estimaor.
\end{prf}

\subsection{Robustness of optimization algorithms}
We say that an optimization problem is \emph{robust}, if whenever it recovers the ground truth for an input tensor $\G$, then it also recovers the ground truth if an adversary modifies $\G$ by applying an arbitrary sequence of monotone
transformations. 
As a consequence, if a robust optimization problem recovers the ground truth with high probability under the fully-random corruption model, then it also recovers the ground truth with high probability under the associated semi-random model.
In the following, we prove that for rank-one BTF, any reasonable optimization algorithm is robust. To the best of our knowledge, there exists no robust spectral method for this problem.

We start by introducing some notation. As before, given binary vectors
$\bar x \in \{0,1\}^n$, $\bar y \in \{0,1\}^m$, $\bar z \in \{0,1\}^l$, define the ground truth rank-one tensor $\bar \W = (\bar w_{ijk}) := \bar x \otimes \bar y \otimes \bar z  \in \{0,1\}^{n \times m \times l}$.
Define:
\begin{align*}
& \P := \{(i,j,k) \in [n] \times [m] \times [l]: \bar x_i = \bar y_j = \bar z_k = 1\}, \\
& \N := \{(i,j,k)\in [n] \times [m] \times [l] : \bar x_i = 0 \; \lor \; \bar y_j = 0 \; \lor \;  \bar z_k = 0\}, \\
& \T := \{(i,j,k) \in [n] \times [m] \times [l]: g_{ijk} = \bar x_i \bar y_j \bar z_k \}, \\
& \F: = \{(i,j,k) \in [n] \times [m] \times [l]: g_{ijk} = 1-\bar x_i \bar y_j \bar z_k\}.
\end{align*}
It then follows that
$$S_0 = (\T \cap \N) \cup (\F \cap \P), \qquad
S_1 = (\T \cap \P) \cup (\F \cap \N),$$
where $S_0, S_1$ are defined in~\eqref{s1s2}. 

\begin{proposition}
\label{robustLPs}
Denote by Problem~(Rec) a minimization problem whose objective function is identical to that of Problem~\eqref{tfr1} and whose feasible region satisfies $0 \leq w_{ijk} \leq 1$ for all $i\in [n]$, $j \in [m]$, $k \in [l]$. Then Problem~(Rec) is robust.
\end{proposition}

\begin{prf}
For notational simplicity, we assume that the optimization variables are $(x,y,z,\W)$. However, the proof follows from the same line of arguments, if Problem~(Rec) contains additional variables. 
Consider $\noisy, \adversary \in \{0,1\}^{n \times m \times l}$ such that $\adversary$ can be obtained from $\noisy$ via a sequence of monotone transformations. 
Denote by OP1, Problem~(Rec) with input tensor $\G$ and denote by OP2, Problem~(Rec) with input tensor $\G'$.
We assume that $(\bar x, \bar y, \bar z, \bar \W)$ is the unique optimal solution of OP1. 
We show that $(\bar x, \bar y, \bar z, \bar \W)$ is the unique optimal solution of OP2 as well.

Define $\T' := \{(i,j,k) \in [n] \times [m] \times [l]: g'_{ijk} = \bar x_i \bar y_j \bar z_k \}$ and define $Q := \T' \setminus \T$. Notice that by definition $\T' \supset \T$.
Since $(\bar x, \bar y, \bar z, \bar \W)$ is the unique optimal solution of OP1, the optimal value of this optimization problem is given by $|\F|$. Moreover, $(\bar x, \bar y, \bar z, \bar \W)$ is a feasible solution of OP2 with the objective value equal to $|\F|-|Q|$. Suppose that $(\bar x, \bar y, \bar z, \bar \W)$ is not the unique optimal solution of OP2; then there exists a solution $(\tilde x, \tilde y, \tilde z, \tilde \W)$ different from $(\bar x, \bar y, \bar z, \bar \W)$ whose objective value is less that or equal to $|\F|-|Q|$. 
Define $S'_0 := \{(i,j,k) \in [n]\times[m]\times[l]: g'_{ijk}=0\}=(\T' \cap \N) \cup (\F' \cap \P)$ and
$S'_1 := \{(i,j,k) \in [n]\times[m]\times[l]: g'_{ijk}=1\}=(\T' \cap \P) \cup (\F' \cap \N)$, where $\F' = \F \cup Q$.
Let us examine the objective value $\tilde f$ of OP2 at $(\tilde x, \tilde y, \tilde z, \tilde \W)$:    
\begin{align}\label{fa1}
    \tilde f = & \sum_{(i,j,k) \S'_0}{\tilde w_{ijk}}+
    \sum_{(i,j,k) \S'_1}{(1-\tilde w_{ijk})}\\
    =& \sum_{(i,j,k) \in \S_0}{\tilde w_{ijk}}+\sum_{(i,j,k) \in Q \cap \N}{\tilde w_{ijk}}-\sum_{(i,j,k) \in Q \cap \P}{\tilde w_{ijk}}+
    \sum_{(i,j,k) \in \S_1}{(1-\tilde w_{ijk})}+\sum_{(i,j,k) \in Q \cap \P}{(1-\tilde w_{ijk})}\nonumber\\
    &-\sum_{(i,j,k) \in Q \cap \N}{(1-\tilde w_{ijk})}\nonumber\\
    =& \sum_{(i,j,k) \in \S_0}{\tilde w_{ijk}} +
    \sum_{(i,j,k) \in \S_1}{(1-\tilde w_{ijk})} -
    \Big(\sum_{(i,j,k) \in Q \cap \P}{(2\tilde w_{ijk}-1)}-\sum_{(i,j,k) \in Q \cap \N}{(2\tilde w_{ijk}-1)}\Big)\nonumber\\
    \leq & |F|-|Q|.
\end{align}
Moreover, since by assumption $0 \leq \tilde w_{ijk} \leq 1$ for all $(i,j,k) \in [n] \times [m] \times [l]$, we have:
\begin{equation}\label{fa2}
 \sum_{(i,j,k) \in Q \cap \P}{(2\tilde w_{ijk}-1)}-\sum_{(i,j,k) \in Q \cap \N}{(2\tilde w_{ijk}-1)} \leq |Q \cap \P|+|Q \cap \N| =|Q|.  \end{equation}
From~\eqref{fa1} and~\eqref{fa2} it follows that
\begin{equation}
\sum_{(i,j,k) \in \S_0}{\tilde w_{ijk}} +
    \sum_{(i,j,k) \in \S_1}{(1-\tilde w_{ijk})} \leq |F|,    
\end{equation}
which contradicts with the assumption that $(\bar x, \bar y, \bar z, \bar \W)$ is the unique optimal solution of OP1. Therefore,
$(\bar x, \bar y, \bar z, \bar \W)$ is the unique optimal solution of OP2 as well and this completes the proof.
\end{prf}

In particular, \cref{robustLPs} 
enables us to establish recovery under the semi-random corruption model by proving recovery under the simpler fully-random corruption model.

\section{LP relaxations for rank-one BTF}
\label{sec LP4BTF}
A simple LP relaxation of Problem~\eqref{tfr1} can be obtained by replacing each multilinear term
$w_{ijk} = x_i y_j z_k$, $x_i, y_j, z_k \in \{0,1\}$, by its convex hull~\cite{Crama93}.
Using the sign of objective function coefficients, we remove a subset of constraints that are never active at an optimal solution to obtain:
\begin{align}
\label{LP1}
\tag{sLP}
\min \quad & \sum_{(i,j,k) \in S_0}{w_{ijk}}+\sum_{(i,j,k) \in S_1}{(1-w_{ijk})}\nonumber \\
\st  \quad & w_{ijk} \leq x_i, \; w_{ijk} \leq y_j, \; w_{ijk} \leq  z_k, \qquad \forall (i,j,k) \in S_1\label{cp}\\
& w_{ijk}\geq 0, \; w_{ijk} \geq x_i + y_j + z_k -2, \qquad \forall (i,j,k) \in S_0 \label{cn}\\
& (x,y,z) \in [0,1]^{n+m+l} \label{ub}.
\end{align}
Throughout this paper, we refer to Problem~\eqref{LP1} as the \emph{standard LP}.
In Sections~\ref{sec main results} and~\ref{sec numerics}, we investigate the effectiveness of this LP theoretically and numerically, respectively.  Next, leveraging on our previous results regarding the facial structure of the convex hull of multilinear
sets~\cite{dPKha16,dPIda16,dPIda18,dPIda19,dPIdaSah20,dPKha23MPA}, we propose stronger LP relaxations for rank-one BTF.

\subsection{The multilinear polytope and new LP relaxations}
We start by providing a brief overview of the multilinear polytope.
Subsequently, we propose new LP relaxations for rank-one BTF.
Consider a hypergraph $G=(V,E)$, where $V$ is the set of nodes, and $E$ is the set of edges, where each edge is a subset of $V$ of cardinality at least two.
The \emph{multilinear set} $\S_G$ is defined as the set of binary points $(u, w) \in \{0,1\}^{V} \times \{0,1\}^{E}$ satisfying the collection of equations $w_e = \prod_{v \in e}{u_v}$, for all $e \in E$.
The \emph{multilinear polytope} $\MP_G$ is defined as the convex hull of the multilinear set $\S_G$.
It then follows that the feasible region of Problem~\eqref{tfr1} is a highly structured multilinear set and hence understanding the facial structure of its convex hull is key to constructing strong LP relaxations for rank-one BTF.

In~\cite{dPIda18,dPIda19,dPKha23MPA} we obtain sufficient conditions, in terms of acyclicity degree of the hypergraph, under which the multilinear polytope has a compact extended formulation. As a byproduct, in these papers we introduce new classes of valid inequalities for the multilinear polytope; namely, flower inequalities~\cite{dPIda18}, running intersection inequalities~\cite{dPIda19} and their extensions~\cite{dPKha23MPA}.
For more theoretical results on the facial structure of the multilinear polytope, we refer the reader to
\cite{SheAda90,sherali97,TawRich13,dPKha16,dPIda16,CraRod17,BucCraRod18,dPDiG21,dPWal22}.

Henceforth, we refer to the convex hull of the feasible region of Problem~\eqref{tfr1} as the \emph{multilinear polytope of rank-one BTF}, and we denote by $G^{\rm BT}= (V, E^{\rm BT})$ its hypergraph.
In particular, $G^{\rm BT}$ is a \emph{tripartite} hypergraph: it has $n+m+l$ nodes, and $nml$ edges, each edge contains three nodes: one from the first $n$, one from the second $m$, and one from the last $l$. Now define the hypergraph $\bar G= (V, E^{\rm BT} \cup \bar E)$, where $\bar E$ contains all subsets of cardinality two of each $e \in E^{\rm BT}$. 
It then follows that the projection of $\MP_{\bar G}$ onto the space defined by $G^{\rm BT}$ coincides with $\MP_{G^{\rm BT}}$. Hence, to obtain a polyhedral relaxation for $\S_{G^{\rm BT}}$, it suffices to obtain a polyhedral relaxation for $\S_{\bar G}$. To this end, for each $(i,j,k) \in [n]\times[m]\times[l]$, we replace the multilinear set defined by the three equations $w_{ijk}= x_i y_j z_k$, $w^1_{ij}= x_i y_j$, $w^2_{ik}= x_i z_k$, $w^3_{jk}=y_j z_k$, $x_i,y_j,z_k \in \{0,1\}$ by its convex hull~\cite{SheAda90} to obtain the following LP relaxation of Problem~\eqref{tfr1}:
\begin{align}
\label{LPadd}
\tag{cLP}
\min \quad & \sum_{(i,j,k) \in S_0}{w_{ijk}}+\sum_{(i,j,k) \in S_1}{(1-w_{ijk})}\nonumber \\
\st  \quad & w_{ijk} \leq w^1_{ij}, \; w_{ijk} \leq w^2_{ik}, \; w_{ijk} \leq  w^{3}_{jk}, \quad \forall (i,j,k) \in S_1\label{addeq1}\\
&x_i+y_j+z_k-w^1_{ij}-w^2_{ik}-w^3_{jk}+w_{ijk} \leq 1, \quad \forall (i,j,k) \in S_1\nonumber\\
& w^1_{ij}+w^2_{ik}-w_{ijk} \leq x_i, \quad \forall (i,j,k) \in S_0\label{addeq2}, \\
& w^1_{ij}+w^3_{jk}-w_{ijk} \leq y_j, \quad \forall (i,j,k) \in S_0\label{addeq3}, \\
& w^2_{ik}+w^3_{jk}-w_{ijk} \leq z_k, \quad \forall (i,j,k) \in S_0\label{addeq4}, \\
& w_{ijk}\geq 0, \quad \forall (i,j,k) \in S_0 \nonumber\\
& (x,y,z) \in [0,1]^{n+m+l}, \; (w^1, w^2, w^3) \in [0,1]^{nm+nl+ml} \nonumber,
\end{align}
where as before we removed some constraints that are never active at an optimal solution. Henceforth, we refer to Problem~\eqref{LPadd} as the \emph{complete LP}.
As we will show in \cref{sec numerics}, the complete LP is significantly stronger than the standard LP. Notice that the standard LP consists of $nml+n+m+l$ variables, and at most $3nml$ constraints, while the complete LP consists of $nml+nm+nl+ml+n+m+l$ variables, and at most $4nml$ constraints. This moderate increase in size has a rather insignificant impact on the computational cost of solving Problem~\eqref{LPadd}. Indeed, the complete LP can be solved very fast using off-the-shelf LP solvers such as {\tt CPLEX} and {\tt Gurobi}.

Let us detail on the connections between the complete LP and alternative LP relaxations for Problem~\eqref{tfr1}. 
First, consider the Reformulation Relaxation Technique (RLT) hierarchy~\cite{SheAda90} for Problem~\eqref{tfr1}. It can be checked that
the constraints of Problem~\eqref{LPadd} are present in level-2 RLT of Problem~\eqref{tfr1}. However, level-2 RLT of Problem~\eqref{tfr1} contains $\Theta((n+m+l)^3)$ variables and constraints, and hence is significantly more expensive to solve.

As we mentioned before, flower inequalities~\cite{dPIda18} and running intersection inequalities~\cite{dPIda19} are among the most popular valid inequalities for the multilinear polytope and both theoretical and computational benefits of these inequalities for binary polynomial optimization have been investigated in the literature~\cite{dPIdaSah20,Aida22}. However, for the multilinear polytope of rank-one BTF $\MP_{G^{\rm BT}}$, it can be shown that all flower inequalities and running intersection inequalities are implied by the feasible region of Problem~\eqref{LPadd} (see proof of Theorem~3 in~\cite{dPIda19}).  In fact, inequalities~\eqref{addeq1} are flower inequalities of $\MP_{\bar G}$, and inequalities~\eqref{addeq2}--\eqref{addeq4} are 
running intersection inequalities of $\MP_{\bar G}$.

Obtaining recovery guarantees for the complete LP is rather involved.
To investigate the impact of the inequalities defining the complete LP in improving recovery properties of the standard LP, 
we consider a specific relaxation of Problem~\eqref{LPadd}, which we will refer as the \emph{flower LP}:
\begin{align}\label{LP3}
\tag{fLP}
\min \quad & \sum_{(i,j,k) \in S_0}{w_{ijk}}+\sum_{(i,j,k) \in S_1}{(1-w_{ijk})}\nonumber \\
\st  \quad & w_{ijk} \leq  x_i, \ w_{ijk} \leq y_j, \ w_{ijk} \leq  z_k, \qquad \forall (i,j,k) \in S_1  \label{f1}\\
& w_{ijk} \geq 0, \quad  w_{ijk} \geq x_i + y_j + z_k -2, \qquad \forall (i,j,k) \in S_0 \label{f2}\\
& w_{i'jk} - w_{ijk} \leq 1 - x_i, \qquad \forall (i,j,k) \in S_0, \; (i',j,k) \in S_1\label{f3} \\
& w_{ij'k} - w_{ijk} \leq 1 - y_j, \qquad \forall (i,j,k) \in S_0, \; (i,j',k) \in S_1 \label{f4} \\
& w_{ijk'} - w_{ijk} \leq 1 - z_k, \qquad \forall (i,j,k) \in S_0, \; (i,j,k') \in S_1 \label{f5} \\
& (x,y,z) \in [0,1]^{n+m+l}.\nonumber
\end{align}
Inequalities~\eqref{f3}-~\eqref{f5} are flower inequalities of $\MP_{G^{\rm BT}}$ and as we mentioned before are implied by the feasible region of the complete LP.
As the feasible region of Problem~\eqref{LPadd} is contained in that of Problem~\eqref{LP3}, we conclude that if the flower LP succeeds in recovering the ground truth, so does the complete LP. 
Due to its simple formulation, Problem~\eqref{LP3} is simpler to analyze than Problem~\eqref{LPadd}, yet, as we detail in Section~\ref{sec main results}, it significantly outperforms Problem~\eqref{LP1} in recovering the ground truth.
It is important to note that  while in general, flower inequalities are not facet-defining for the multilinear polytope, as we show next, they define facets of the multilinear polytope of rank-one BTF.

\begin{proposition}
\label{prop facets}
All inequalities defining the feasible region of flower LP, \ie inequalities~\eqref{f1}-~\eqref{f5}, define facets of the multilinear polytope of rank-one BTF.
\end{proposition}


The proof of \cref{prop facets} relies on standard techniques and is given in \cref{sec facets}.

\section{Statements of recovery guarantees}
\label{sec main results}
In this section, we summarize the main results of this paper. The proofs are deferred to \cref{sec IP recovery,sec SL recovery,sec FR recovery}.

\subsection{Recovery conditions for rank-one BTF}

We start by characterizing the corruption range, in terms of $p$, for which 
rank-one BTF 
can recover the ground truth with high probability.
Such conditions 
serve as a reference point for assessing the effectiveness of our LP relaxations.
Our 
results for rank-one BTF essentially indicate that Problem~\eqref{pr comb} recovers the ground truth with high probability if and only if $p < 1/2$.
These results can be seen as tight, since in the fully-random model for $p=1/2$, each entry of the noisy tensor $\noisy$ is zero or one with equal probability, no matter what the original ground truth rank-one tensor $\ground$ is.

First, we present 
necessary conditions under which Problem~\eqref{pr comb} recovers the ground truth under fully-random corruption model.
Note that  a similar condition for the semi-random model cannot be given as the adversary may choose to undo all corruptions.

\begin{theorem}
[Necessary conditions for recovery]
\label{th it lb}
Consider the fully-random corruption model.
If $\pr \ge 1/2$, then the probability that Problem~\eqref{pr comb} recovers the ground truth is at most $1/2$.
Furthermore, if $r_{\bar x}, r_{\bar y}, r_{\bar z}, p$ are positive constants and $p > 1/2$, then with high probability Problem~\eqref{pr comb} does not recover the ground truth.
\end{theorem}

Next, we give 
sufficient conditions under which Problem~\eqref{pr comb} recovers the ground truth with high probability.
This result holds for the more general semi-random corruption model.

\begin{theorem}
[Sufficient conditions for recovery]
\label{th it ub}
Consider the semi-random corruption model.
Assume that $r_{\bar x}, r_{\bar y}, r_{\bar z}$ are positive constants and
\begin{equation}\label{lass}
\lim_{n, m, l \to \infty} \frac{n+m+l}{\min\{nm, nl, ml\}} = 0.
\end{equation}
If $\pr$ is a constant satisfying $p < 1/2$, then Problem~\eqref{pr comb} recovers the ground truth with high probability.
\end{theorem}

Proofs of the above theorems are given in \cref{sec IP recovery}.
In fact, in \cref{sec IP recovery}, we also present more general 
conditions for recovery
in which we do not assume that $r_{\bar x}, r_{\bar y}, r_{\bar z}, p$ are constants (see \cref{prop it lb,prop it ub}).

\begin{remark}\label{rem1}
The limit assumption~\eqref{lass} in \cref{th it ub} is not too restrictive.
Consider $m$ and $l$ as functions of $n$, i.e., $m = m(n)$ and $l = l(n)$.
Furthermore, assume that $m(n)$ grows faster than $n$ and that $l(n)$ grows faster than $m(n)$, i.e. $n \in O(m(n))$, $m(n) \in O(l(n))$.
Then assumption~\eqref{lass} is satisfied if $l(n) \in o(nm(n))$.
Intuitively, sufficient conditions of~\cref{th it ub} require that the functions $m(n)$ and $l(n)$ grow similarly as $n$ increases.
Two simple examples of functions that satisfy these assumptions are:
(\ref{th it ub}.i). $m(n), l(n) \in \Theta(n^k)$, for any positive integer $k$;
(\ref{th it ub}.ii). $m(n), l(n) \in \Theta(\exp(n))$.
An example of functions that do not satisfy the assumptions of~\cref{th it ub} is: $m(n) \in \Theta(n^k)$ and $l(n) \in \Theta(n^{k+1})$, for any positive integer $k$.
\end{remark}

\subsection{Recovery conditions for the standard LP}

Next, we present recovery guarantees for the standard LP; namely, Problem~\eqref{LP1}.
In particular, we obtain a sufficient condition in terms of $p, r_{\bar x}, r_{\bar y}, r_{\bar z}$ under which the standard LP recovers the ground truth with high probability.

\begin{theorem}
\label{th LP1}
Consider the semi-random corruption model. 
Assume that $r_{\bar x}, r_{\bar y}, r_{\bar z}$ are positive constants, and without loss of generality, assume $r_{\bar x} \geq r_{\bar y} \geq  r_{\bar z}$.
Assume that, as $n, m, l \to \infty$, we have $n \exp (-m)$, $n \exp (-l)$, $m \exp (-n)$, $m \exp (-l)$, $l \exp (-n)$, $l \exp (-m) \to 0$.
If $\pr$ is a constant satisfying
\begin{equation}\label{cond2}
p < \frac{r_{\bar x} r_{\bar y} r_{\bar z}}{2 (1+ r_{\bar x} r_{\bar y} r_{\bar z}) + \delta},
\end{equation}
where
$$\delta := \frac{1}{2} r_{\bar x} r_{\bar y} + \frac{1}{2} r_{\bar x} r_{\bar z} - r_{\bar y} r_{\bar z} + r_{\bar x} - \frac{1}{2} r_{\bar y} - \frac{1}{2} r_{\bar z} \geq 0,$$
then Problem~\eqref{LP1} recovers the ground truth with high probability.
\end{theorem}

Note that $\delta \geq 0$ since the assumption $r_{\bar x} \geq r_{\bar y} \geq r_{\bar z}$ implies that
$r_{\bar x} r_{\bar y} + r_{\bar x} r_{\bar z} \geq  2 r_{\bar y} r_{\bar z}$ and $2 r_{\bar x} \geq r_{\bar y} + r_{\bar z}$, and these inequalities are satisfied tightly if and only if $r_{\bar x} = r_{\bar y} = r_{\bar z}$. Now suppose that the ground truth satisfies $r_{\bar x} = r_{\bar y} = r_{\bar z}$, and denote by $r_{\bar w}$ the \emph{tensor density}; \ie $r_{\bar w} := r_{\bar x} r_{\bar y} r_{\bar z}$.
We then obtain the following corollary of \cref{th LP1}:
\begin{corollary}
\label{corLP1}
Consider the semi-random corruption model.
Suppose that $r_{\bar x}, r_{\bar y}, r_{\bar z}$ are positive constants satisfying $r_{\bar x}= r_{\bar y}= r_{\bar z}$ and let $r_{\bar w} := r_{\bar x} r_{\bar y} r_{\bar z}$.
Assume that, as $n, m, l \to \infty$, we have $n \exp (-m)$, $n \exp (-l)$, $m \exp (-n)$, $m \exp (-l)$, $l \exp (-n)$, $l \exp (-m) \to 0$.
If $\pr$ is a constant satisfying
$$
p < \frac{r_{\bar w}}{2 (1+r_{\bar w})},
$$
then Problem~\eqref{LP1} recovers the ground truth with high probability.
\end{corollary}

The proof of \cref{th LP1} is given in \cref{sec SL recovery}.
To prove this result, we first obtain a deterministic sufficient condition for recovery (see \cref{deter1,unique1}).
Then, using the deterministic condition together with~\cref{robustLPs}, we derive a recovery guarantee under the semi-random corruption model.

\begin{remark}\label{rem2}
The limit assumptions in \cref{th LP1} are not too restrictive.
Consider $m$ and $l$ as functions of $n$, i.e., $m = m(n)$ and $l = l(n)$.
Two simple examples of functions that satisfy the assumptions of~\cref{th LP1} are:
(\ref{th LP1}.i). $m(n) \in \Theta(n^h)$ and $l(n) \in \Theta(n^k)$, for positive integers $h,k$;
(\ref{th LP1}.ii). $m(n), l(n) \in \Theta(\exp(n/2))$.
It is important to note that the limit assumptions in \cref{th LP1} are not comparable to those in \cref{th it ub}.
In particular, (\ref{th LP1}.i) contains as special case the example given in Remark~\ref{rem1} (\ie $m(n) \in \Theta(n^k)$, $l(n) \in \Theta(n^{k+1})$) that does not satisfy assumption~\eqref{lass}.
On the other hand, example~(\ref{th it ub}.ii) does not satisfy the limit assumptions in \cref{th LP1}.
Clearly, any LP relaxation of rank-one BTF recovers the ground truth, only if the original integer program succeeds in doing so. The fact that the assumptions of Theorem~\ref{th LP1} are not implied by those of~\cref{th it ub} is merely an artifact of our analysis. 
\end{remark}

\subsection{Recovery conditions for the flower LP}

Next, we present recovery a guarantee for the flower LP; namely, Problem~\eqref{LP3}:

\begin{theorem}
\label{th LP2}
Consider the semi-random corruption model.
Assume that $r_{\bar x}, r_{\bar y}, r_{\bar z}$ are positive constants and let $r_{\bar w} := r_{\bar x} r_{\bar y} r_{\bar z}$.
Assume that, as $n, m, l \to \infty$, we have $nml \exp (-n)$, $nml \exp (-m)$, $nml \exp (-l) \to 0$.
If $\pr$ is a constant satisfying
\begin{equation}\label{condaux simple}
\pr <
\frac{r_{\bar w}}{1+2 r_{\bar w}},
\end{equation}
then Problem~\eqref{LP3} recovers the ground truth with high probability.
\end{theorem}

The proof of \cref{th LP2} is given in \cref{sec FR recovery}.
Our proof scheme is similar to that of Theorem~\ref{th LP1}:
we first obtain a deterministic sufficient condition for recovery (see~\cref{deter2,unique2}); next we consider the semi-random corruption model.

\begin{remark}
The limit assumptions in \cref{th LP2} are not too restrictive, even though they are stronger than the limit assumptions in \cref{th LP1}.
Consider $m$ and $l$ as functions of $n$, i.e., $m = m(n)$ and $l = l(n)$.
A simple example of functions that satisfies the assumptions of \cref{th LP2} is example~(\ref{th LP1}.i) of Remark~\ref{rem2}.
On the other hand, example~(\ref{th LP1}.ii) of Remark~\ref{rem2} does not satisfy the assumptions of~\cref{th LP1}.
An example of functions that satisfies the assumptions in all \cref{th it ub,th LP1,th LP2} is example~(\ref{th it ub}.i) of Remark~\ref{rem1}. Again notice that the standard LP recovers the ground truth only if the flower LP succeeds in doing so, and that fact that the assumptions of Theorem~\ref{th LP2} are not implied by those of Theorem~\ref{th LP1} is merely an artifact of our analysis.
\end{remark}

\begin{figure}[htbp]
 \centering
 \psfrag{p}{$p$}
 \psfrag{r}{$r_{\bar w}$}
\epsfig{figure=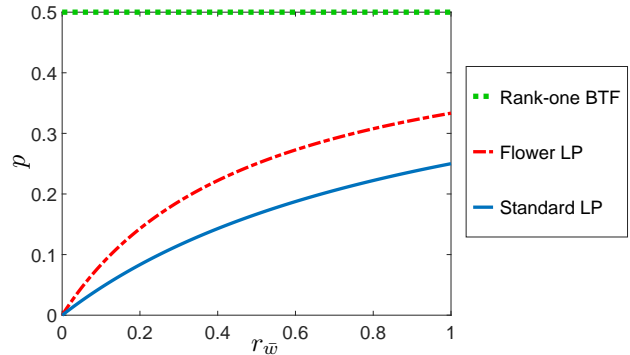, scale=0.25, trim=10mm 0mm 10mm 0mm, clip}
 \caption{Comparison of the recovery threshold for two LP relaxations of Rank-one BTF under the semi-random corruption model as given by Corollary~\ref{corLP1} and Theorem~\ref{th LP2}.}
\label{figure1}
\end{figure}

The recovery thresholds of the two LP relaxations for rank-one BTF together with rank-one BTF are depicted in \cref{figure1}: as the recovery threshold of Corollary~1 serves as an upper bound for the recovery threshold of standard LP, we conclude that the addition of flower inequalities significantly improves the recovery properties of the LP relaxation.

We conclude this section by observing that the parameter $p$ is not an input to our LP relaxation schemes.
Hence, the knowledge of such parameter is not necessary for our recovery guarantees, which is an important and desirable property. Moreover, to obtain recovery guarantees we do not make any assumption on how the ground truth tensor $\bar \W$ was generated.

\section{Numerical experiments}
\label{sec numerics}
In this section, we conduct a preliminary numerical study to compare the recovery properties of the proposed LP relaxations for rank-one BTF. A comprehensive computational study that includes real data sets from the literature is a topic of future research.

\begin{figure}[htbp]
 \centering
 \psfrag{p}{$p$}
 \psfrag{r}{$r_{\bar w}$}
\subfigure[\ Standard LP~\label{fig2a}]{\epsfig{figure=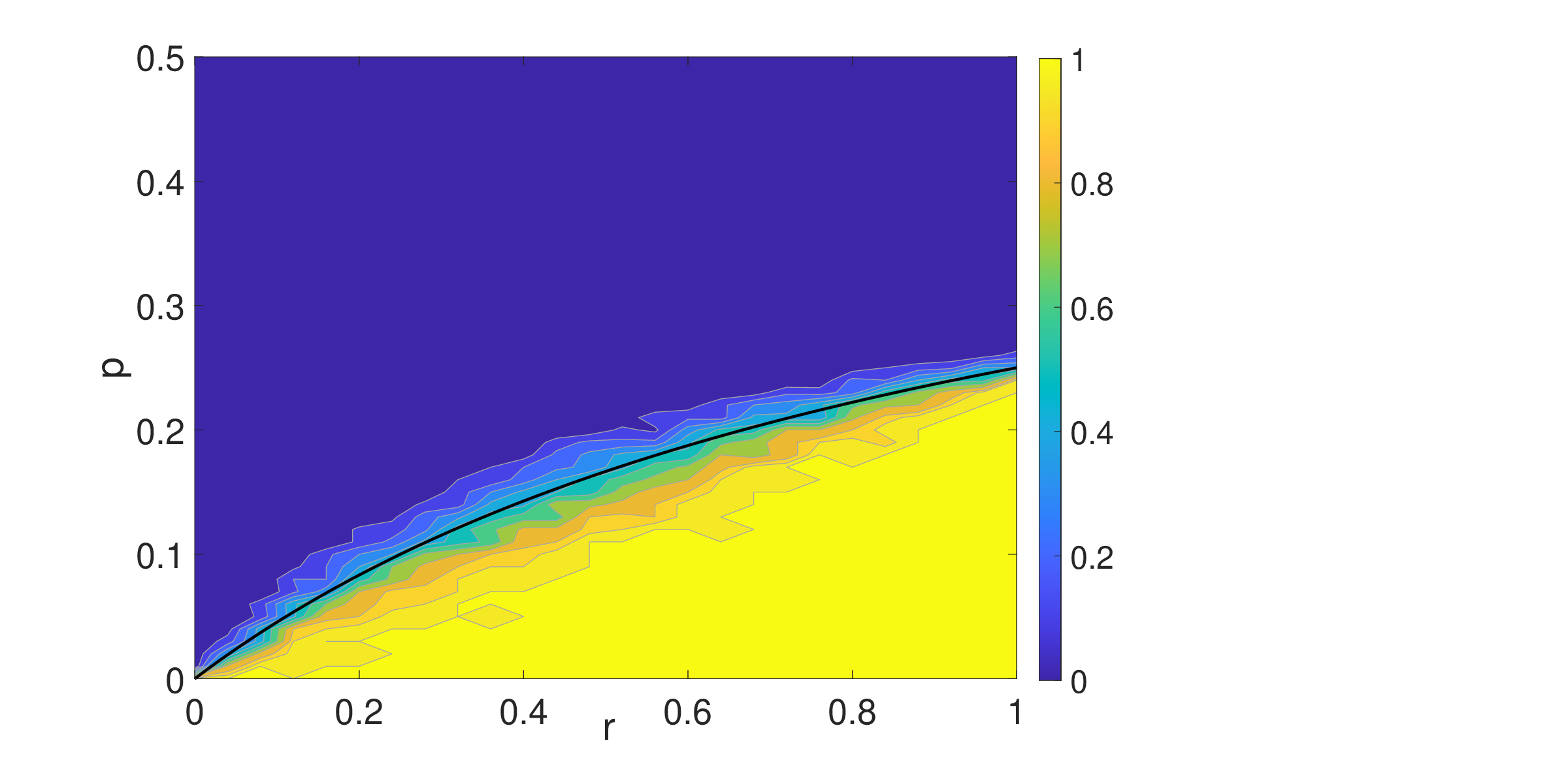, scale=0.25, trim=20mm 0mm 110mm 0mm, clip}}
\subfigure[\ Flower LP~\label{fig2b}]{\epsfig{figure=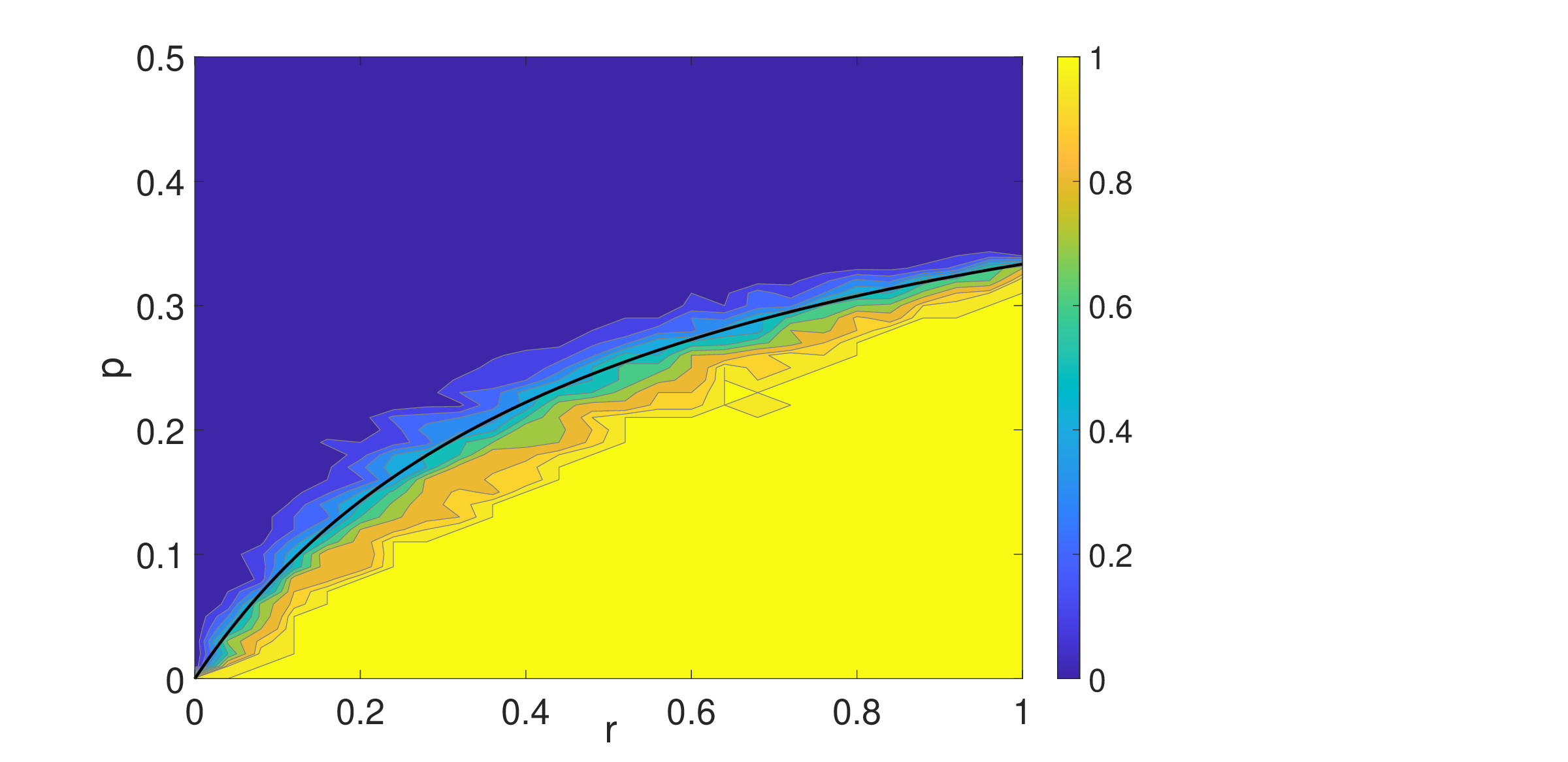, scale=0.25, trim=20mm 0mm 110mm 0mm, clip}}
\subfigure[\ Complete LP~\label{fig2c}]{\epsfig{figure=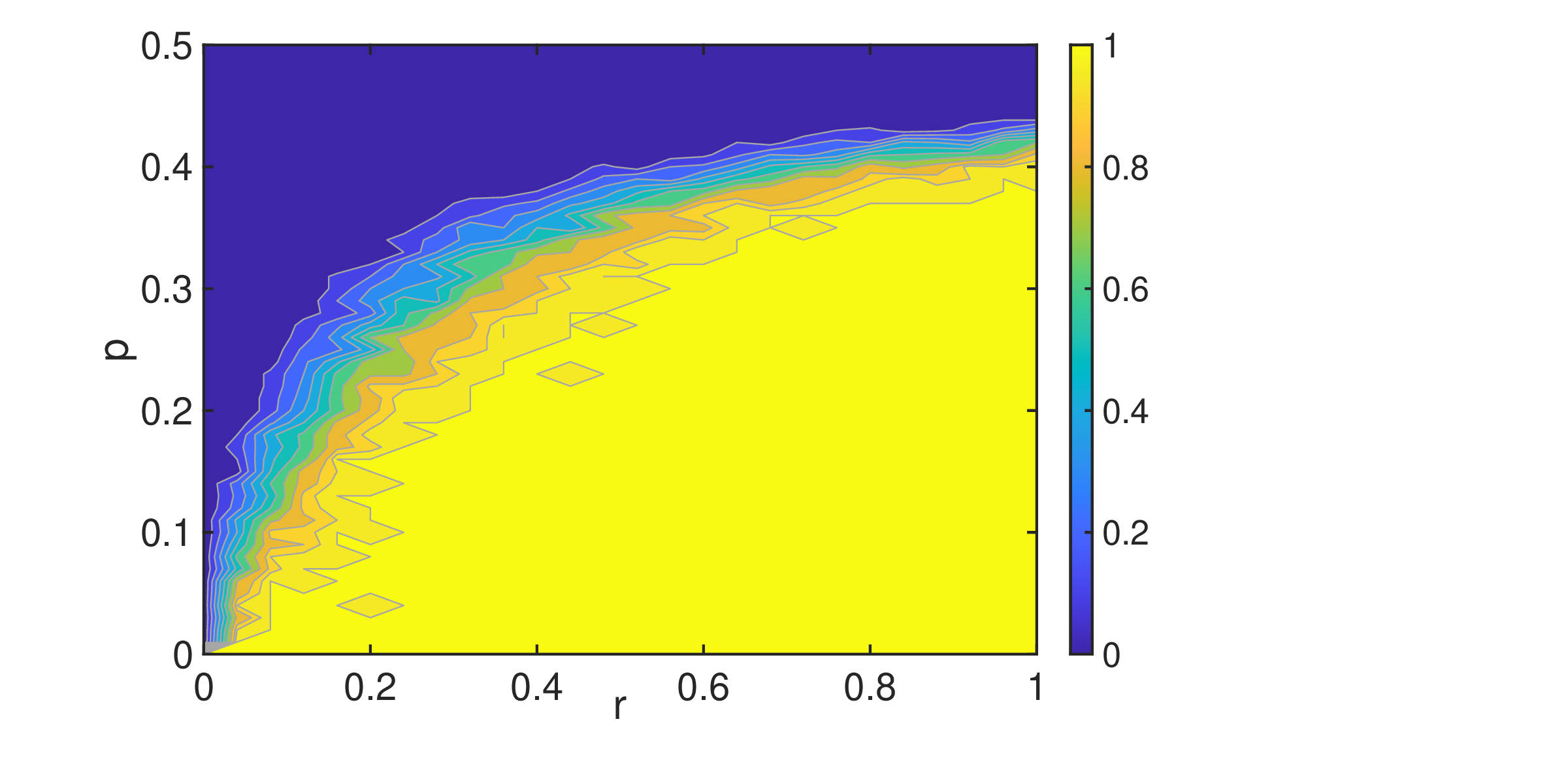, scale=0.25, trim=20mm 0mm 110mm 0mm, clip}}
 \caption{The empirical rate of recovery for LP relaxations of Rank-one BTF under the fully-random corruption model. The black curves in Figures~\ref{fig2a} and~\ref{fig2b} correspond to the recovery guarantees given by Corollary~\ref{corLP1} and Theorem~\ref{th LP2}, respectively.
}
\label{figure2}
\end{figure}

We consider three LP relaxations for rank-one BTF: (i) the standard LP, defined by Problem~\eqref{LP1}, (ii) the flower LP, defined by Problem~\eqref{LP3}, and (iii) the complete LP, defined by Problem~\eqref{LPadd}. We generate the input tensor $\G$ according to the fully-random corruption model defined before. For our numerical experiments we set $n= m = l= 15$, $r_{\bar w} \in [0:0.04:1.0]$, and $p \in [0:0.01:0.5]$. 
Given $r_{\bar w}$, we construct the ground truth tensor $\bar \W = \bar x \otimes \bar y \otimes \bar z$ as follows: we set $q := \sqrt[3]{r_{\bar w}}$ and let $\bar x$, $\bar y$, $\bar z$ be vectors of Bernoulli random variables with parameter $q$. 
It then follows that $\avg[\sum_{ijk}{\bar w_{ijk}}] = q^3 = r_{\bar w}$.
For each fixed pair $(r_{\bar w}, p)$, we run 40 random trials. 
We then count the number of times each LP relaxation recovers the ground truth. Dividing by the number of trials, we
obtain the empirical rate of recovery. All experiments are performed on the {\tt NEOS} server~\cite{neos98} and all LPs are solved with {\tt GAMS/CPLEX}~\cite{cplex}. Results are shown in Figure~\ref{figure2}; as can be seen from Figures~\ref{fig2a} and~\ref{fig2b}, our recovery guarantees given by~\cref{corLP1} and~\cref{th LP2} are fairly tight, and we conjecture that these conditions are necessary for recovery as well. Moreover, Figure~\ref{fig2c} indicates that the complete LP significantly outperforms the flower LP and hence understanding its recovery properties is a topic of future research.

Finally, we would like to remark that while the recovery threshold of the original NP-hard problem is independent of the density of the input tensor, our theoretical and numerical results indicate that for all considered LP relaxations of rank-one BTF, the recovery threshold improves as the tensor density increases. 


\section{Recovery proofs for rank-one BTF}
\label{sec IP recovery}

The main goal of this section is to prove \cref{th it lb,th it ub}.
We start by introducing some notation that will be used in this section.
Let $(\bar x, \bar y, \bar z) \in \{0,1\}^{n+m+l}$.
Like in the corruption models, we denote by $r_{\bar x}$ (resp. $r_{\bar y}$, $r_{\bar z}$) the ratio of ones in $\bar x$ (resp. $\bar y$, $\bar z$) to the number of elements in $\bar x$ (resp. $\bar y$, $\bar z$).
For every $(x,y,z) \in \{0,1\}^{n+m+l}$, we denote by
\begin{align*}
\Delta_{(x,y,z)} := \{(i,j,k) \in [n] \times [m] \times [l] \mid \bar x_i \bar y_j \bar z_k \neq x_i y_j z_k \}, \qquad \delta_{(x,y,z)} := \card{\Delta_{(x,y,z)}}.
\end{align*}
For every $(x,y,z) \in \{0,1\}^{n+m+l}$ and every $s,t \in \{0,1\}$, we define the following sets:
\begin{align*}
X_{st} := \{i \in [n] \mid \bar x_i = s, \, x_i = t\}, \
Y_{st} := \{j \in [m] \mid \bar y_j = s, \, y_j = t\}, \
Z_{st} := \{k \in [l] \mid \bar z_k = s, \, z_k = t\}.
\end{align*}
For every $(x,y,z) \in \{0,1\}^{n+m+l}$, we denote by $\obj(x,y,z)$ the objective function of \eqref{pr comb}.
Such function is explicitly given by
\begin{align*}
\obj(x,y,z) := 
\Big\Vert\G - x \otimes y \otimes z \Big\Vert^2
=
\sum_{\substack{i \in [n], \ j \in [m], \ k \in [l]: \\ g_{ijk} = 0}}{x_i y_j z_k}+
\sum_{\substack{i \in [n], \ j \in [m], \ k \in [l]: \\ g_{ijk} = 1}}{(1-x_i y_j z_k)},
\end{align*}
and it gives the number of triples $(i,j,k) \in [n] \times [m] \times [l]$ for which $x_iy_jz_k \neq g_{ijk}$.
%
In the remainder of the paper we also denote by $\prob[A]$ the probability of an event $A$, and by $\avg[X]$ the expected value of a random variable $X$.

\subsection{Some useful lemmas}
\label{sec lemmas}

Next, we present two lemmas that will be useful in our analysis.
In \cref{lem card}, we study the quantity $\delta_{(x,y,z)}$ and provide a lower bound.

\begin{lemma}
\label{lem card}
Let $(\bar x, \bar y, \bar z) \in \{0,1\}^{n+m+l}$.
\begin{enumerate}[label=(\roman*)]
\item
\label{lem card edge sets}
For every $(x,y,z) \in \{0,1\}^{n+m+l}$, we have
\begin{align*}
\delta_{(x,y,z)}
& =
\card{X_{11}} \card{Y_{11}} \card{Z_{10}} +
\card{X_{11}} \card{Y_{10}} \card{Z_{11}} +
\card{X_{10}} \card{Y_{11}} \card{Z_{11}} +
\card{X_{11}} \card{Y_{10}} \card{Z_{10}} +
\card{X_{10}} \card{Y_{11}} \card{Z_{10}} \\
& \, +
\card{X_{10}} \card{Y_{10}} \card{Z_{11}} +
\card{X_{10}} \card{Y_{10}} \card{Z_{10}} +
\card{X_{11}} \card{Y_{11}} \card{Z_{01}} +
\card{X_{11}} \card{Y_{01}} \card{Z_{11}} +
\card{X_{01}} \card{Y_{11}} \card{Z_{11}} \\
& \, +
\card{X_{11}} \card{Y_{01}} \card{Z_{01}} +
\card{X_{01}} \card{Y_{11}} \card{Z_{01}} +
\card{X_{01}} \card{Y_{01}} \card{Z_{11}} +
\card{X_{01}} \card{Y_{01}} \card{Z_{01}}.
\end{align*}
\item
\label{lem card all}
For every $(x,y,z) \in \{0,1\}^{n+m+l}$ such that $(x,y,z) \neq (\bar x, \bar y, \bar z)$, we have
$$
\delta_{(x,y,z)} \ge \min\{r_{\bar x} n r_{\bar y} m, r_{\bar x} n r_{\bar z} l, r_{\bar y} m r_{\bar z} l\}.
$$
\end{enumerate}
\end{lemma}

\begin{prf}
\ref{lem card edge sets}.
Let $(x,y,z) \in \{0,1\}^{n+m+l}$.
We have
\begin{align*}
\Delta_{(x,y,z)} & =
\{(i,j,k) \in [n] \times [m] \times [l] \mid \bar x_i \bar y_j \bar z_k =1, x_i y_j z_k =0 \} \\
& \quad \cup
\{(i,j,k) \in [n] \times [m] \times [l] \mid \bar x_i \bar y_j \bar z_k =0, x_i y_j z_k =1 \},
\end{align*}
where the union is disjoint.
Thus,
\begin{align*}
\delta_{(x,y,z)}
& = (
\card{X_{11}} \card{Y_{11}} \card{Z_{10}} +
\card{X_{11}} \card{Y_{10}} \card{Z_{11}} +
\card{X_{10}} \card{Y_{11}} \card{Z_{11}} +
\card{X_{11}} \card{Y_{10}} \card{Z_{10}} +
\card{X_{10}} \card{Y_{11}} \card{Z_{10}} \\
& \, +
\card{X_{10}} \card{Y_{10}} \card{Z_{11}} +
\card{X_{10}} \card{Y_{10}} \card{Z_{10}}
)
+
(
\card{X_{11}} \card{Y_{11}} \card{Z_{01}} +
\card{X_{11}} \card{Y_{01}} \card{Z_{11}} +
\card{X_{01}} \card{Y_{11}} \card{Z_{11}} \\
& \, +
\card{X_{11}} \card{Y_{01}} \card{Z_{01}} +
\card{X_{01}} \card{Y_{11}} \card{Z_{01}} +
\card{X_{01}} \card{Y_{01}} \card{Z_{11}} +
\card{X_{01}} \card{Y_{01}} \card{Z_{01}}
).
\end{align*}

\ref{lem card all}.
Let $(x,y,z) \in \{0,1\}^{n+m+l}$ such that $(x,y,z) \neq (\bar x, \bar y, \bar z)$.
From \ref{lem card edge sets}, we have
\begin{align*}
\delta_{(x,y,z)}
& \ge
\card{X_{10}} \card{Y_{11}} \card{Z_{11}} +
\card{X_{10}} \card{Y_{11}} \card{Z_{10}} +
\card{X_{10}} \card{Y_{10}} \card{Z_{11}} +
\card{X_{10}} \card{Y_{10}} \card{Z_{10}} \\
& =
\card{X_{10}} (\card{Y_{11}} + \card{Y_{10}}) (\card{Z_{11}} + \card{Z_{10}})
=
\card{X_{10}} r_{\bar y} m r_{\bar z} l.
\end{align*}
Thus, if we assume $\card{X_{10}} \ge 1$, we obtain $\delta_{(x,y,z)} \ge r_{\bar y} m r_{\bar z} l$ and we are done.
Symmetrically, if we assume $\card{Y_{10}} \ge 1$, we obtain $\delta_{(x,y,z)} \ge r_{\bar x} n r_{\bar z} l$ and we are done.
Symmetrically, if we assume $\card{Z_{10}} \ge 1$, we obtain $\delta_{(x,y,z)} \ge r_{\bar x} n r_{\bar y} m$ and we are done.

Thus we now assume $X_{10} = Y_{10} = Z_{10} = \emptyset$, which implies $\card{X_{11}} = r_{\bar x} n$, $\card{Y_{11}} = r_{\bar y} m$, and $\card{Z_{11}} = r_{\bar z} l$.
Since $(x,y,z) \neq (\bar x, \bar y, \bar z)$, at least one of the sets $X_{01}$, $Y_{01}$, $Z_{01}$ is nonempty.
If $\card{X_{01}} \ge 1$, then from \ref{lem card edge sets} we have
\begin{align*}
\delta_{(x,y,z)} \ge \card{X_{01}} \card{Y_{11}} \card{Z_{11}} \ge \card{Y_{11}} \card{Z_{11}} = r_{\bar y} m r_{\bar z} l
\end{align*}
and we are done.
Symmetrically, if $\card{Y_{01}} \ge 1$, we obtain $\delta_{(x,y,z)} \ge r_{\bar x} n r_{\bar z} l$ and we are done.
Symmetrically, if $\card{Z_{01}} \ge 1$, we obtain $\delta_{(x,y,z)} \ge r_{\bar x} n r_{\bar y} m$ and we are done.
\end{prf}

In \cref{lem probs}, we study the probabilities that a vector $(x,y,z) \in \{0,1\}^{n+m+l}$ has value $\obj(x,y,z)$ smaller or larger than $\obj(\bar x, \bar y, \bar z)$ in the fully-random corruption model.

\begin{lemma}
\label{lem probs}
Consider the fully-random corruption model.
Let $(x,y,z) \in \{0,1\}^{n+m+l}$ such that $(x,y,z) \neq (\bar x, \bar y, \bar z)$ and let $\delta := \delta_{(x,y,z)}$.
Then we have:
\begin{enumerate}[label=(\roman*)]
\item
\label{lem probs sum}
$\prob[\obj(x,y,z) > \obj(\bar x, \bar y, \bar z)] = \sum_{\ell = 0}^{\ceil{\delta/2}-1} \binom{\delta}{\ell} \pr^\ell (1-\pr)^{\delta-\ell}.$
\item
\label{lem probs exp good}
If $\pr \le 1/2$, then $\prob[\obj(x,y,z) \le \obj(\bar x, \bar y, \bar z)] \le \exp \pare{- 2 \delta (1/2 - \pr)^2}.$
\item
\label{lem probs exp bad}
If $\pr \ge 1/2$, then $\prob[\obj(x,y,z) \ge \obj(\bar x, \bar y, \bar z)] \le \exp \pare{- 2 \delta (\pr - 1/2)^2}.$
\end{enumerate}
\end{lemma}

\begin{prf}
\ref{lem probs sum}
In the case $\delta = 0$, we have $\obj(x,y,z) = \obj(\bar x, \bar y, \bar z)$ and the sum in the statement is vacuous, thus we are done.
Assume now $\delta > 0$.
Note that $\obj(x,y,z) > \obj(\bar x, \bar y, \bar z)$ if and only if strictly less than half of the $(i,j,k) \in \Delta_{(x,y,z)}$ are corrupted.
We obtain
\begin{align*}
\prob[\obj(x,y,z) > \obj(\bar x, \bar y, \bar z)] = \sum_{\ell = 0}^{\ceil{\delta/2}-1} \binom{\delta}{\ell} \pr^\ell (1-\pr)^{\delta-\ell}.
\end{align*}

\ref{lem probs exp good}
In the case $\delta = 0$, we have $\exp \pare{- 2 \delta (1/2 - \pr)^2} = \exp(0) = 1$, thus we are done.
Assume now $\delta > 0$.
Note that $\obj(x,y,z) \le \obj(\bar x, \bar y, \bar z)$ if and only if at least half of the $(i,j,k) \in \Delta_{(x,y,z)}$ are corrupted.
For every $(i,j,k) \in \Delta_{(x,y,z)}$, let $B_{(i,j,k)}$ be the Bernoulli random variable defined by $B_{(i,j,k)} := 1$, if $(i,j,k)$ is corrupted, and $B_{(i,j,k)} := 0$, if $(i,j,k)$ is not corrupted.
Consider the Binomial random variable
\begin{align*}
S_\delta := \sum_{(i,j,k) \in \Delta_{(x,y,z)}} B_{(i,j,k)}.
\end{align*}
We obtain that $\obj(x,y,z) \le \obj(\bar x, \bar y, \bar z)$ if and only if $S_\delta \ge \delta/2$.
Note that $\avg[S_\delta] = \delta \pr$, thus
\begin{align*}
\prob[\obj(x,y,z) \le \obj(\bar x, \bar y, \bar z)]
= \prob[S_\delta \ge \delta/2]
= \prob[S_\delta - \avg[S_\delta] \ge \delta (1/2 - \pr)].
\end{align*}
By assumption $\pr \le 1/2$ and $\delta >0$, thus in particular $\delta (1/2 - \pr) \ge 0$.
Hence we can use Hoeffding's inequality and obtain
\begin{align*}
\prob[\obj(x,y,z) \le \obj(\bar x, \bar y, \bar z)]
\le \exp \pare{-\frac {2 \delta^2 (1/2 - \pr)^2}{\delta}}
= \exp \pare{- 2 \delta (1/2 - \pr)^2}.
\end{align*}

\ref{lem probs exp bad}
This proof is symmetric to the proof of \ref{lem probs exp good}.
\end{prf}

\subsection{Proof of necessary conditions for recovery}
\label{sec it lb proof}

In this section we use \cref{lem card,lem probs} in~\cref{sec lemmas} to prove \cref{th it lb}.
We start with the following proposition, which provides general necessary conditions under which Problem~\eqref{pr comb} recovers the ground truth.

\begin{proposition}
\label{prop it lb}
Consider the fully-random corruption model.
Assume $\pr > 1/2$ and
\begin{align*}
\lim_{n, m, l \to \infty} \min\{r_{\bar x} n r_{\bar y} m, r_{\bar x} n r_{\bar z} l, r_{\bar y} m r_{\bar z} l\} (\pr - 1/2)^2 = \infty.
\end{align*}
Then with high probability Problem~\eqref{pr comb} does not recover the ground truth.
\end{proposition}

\begin{prf}
Let $(x,y,z) \in \{0,1\}^{n+m+l}$ such that $(x,y,z) \neq (\bar x, \bar y, \bar z)$.
If \eqref{pr comb} recovers the ground truth, then we must have $\obj(x,y,z) \ge \obj(\bar x, \bar y, \bar z)$.
Thus, the probability that \eqref{pr comb} recovers the ground truth is at most $\prob[\obj(x,y,z) \ge \obj(\bar x, \bar y, \bar z)]$.

From \cref{lem probs} \ref{lem probs exp bad} and \cref{lem card} \ref{lem card all}, we obtain
\begin{align*}
\prob[\obj(x,y,z) \ge \obj(\bar x, \bar y, \bar z)] & \le \exp \pare{- 2 \delta_{(x,y,z)} (\pr - 1/2)^2} \\
& \le \exp \pare{- 2 \min\{r_{\bar x} n r_{\bar y} m, r_{\bar x} n r_{\bar z} l, r_{\bar y} m r_{\bar z} l\} (\pr - 1/2)^2}.
\end{align*}
Due to our assumption, as $n,m,l \to \infty$, the argument of the exponential function goes to $-\infty$ and so $\prob[\obj(x,y,z) \ge \obj(\bar x, \bar y, \bar z)] \to 0$ and with high probability \eqref{pr comb} does not recover the ground truth.
\end{prf}

We are now ready to prove \cref{th it lb}.


\begin{prfc}[of \cref{th it lb}]
We first prove the first part of the statement: If $\pr \ge 1/2$, then the probability that Problem~\eqref{pr comb} recovers the ground truth is at most $1/2$.

Let $(x,y,z) \in \{0,1\}^{n+m+l}$ such that $(x,y,z) \neq (\bar x, \bar y, \bar z)$.
If \eqref{pr comb} recovers the ground truth, then we must have $\obj(x,y,z) > \obj(\bar x, \bar y, \bar z)$.
Thus, the probability that \eqref{pr comb} recovers the ground truth is at most $\prob[\obj(x,y,z) > \obj(\bar x, \bar y, \bar z)]$.
From \cref{lem probs} \ref{lem probs sum}, we have
$\prob[\obj(x,y,z) > \obj(\bar x, \bar y, \bar z)] = \sum_{\ell = 0}^{\ceil{\delta/2}-1} \binom{\delta}{\ell} \pr^\ell (1-\pr)^{\delta-\ell},$ where $\delta := \delta_{(x,y,z)}$.
It can then be seen that $\prob[\obj(x,y,z) > \obj(\bar x, \bar y, \bar z)] \le 1/2$ if $\pr \ge 1/2$. This concludes the proof of the first part of the statement.

The second part of the statement, where we assume that $r_{\bar x}, r_{\bar y}, r_{\bar z}, \pr$ are positive constants and $\pr > 1/2$, follows directly from \cref{prop it lb}.
\end{prfc}



\subsection{Proof of sufficient conditions for recovery}
\label{sec it ub proof}

The main goal of this section is to prove \cref{th it ub}.
First, we remark that it follows from \cref{robustLPs} that Problem~\eqref{pr comb} is robust.
We proceed with the following proposition, which provides general sufficient conditions under which Problem~\eqref{pr comb} recovers the ground truth.

\begin{proposition}
\label{prop it ub}
Consider the semi-random corruption model.
Assume $\pr < 1/2$ and
\begin{align*}
\lim_{n, m, l \to \infty} \pare{2 \min\{r_{\bar x} n r_{\bar y} m, r_{\bar x} n r_{\bar z} l, r_{\bar y} m r_{\bar z} l\} \pare{1/2 - \pr}^2 - (n + m + l) \log 2} = \infty.
\end{align*}
Then with high probability \eqref{pr comb} recovers the ground truth.
\end{proposition}


\begin{prf}
Due to \cref{robustLPs}, it suffices to prove the statement of the proposition for the fully-random corruption model.
Denote by $\bar \prob$ the probability that \eqref{pr comb} does not recover the ground truth.
We have
\begin{align*}
\bar \prob
& = \prob [\exists (x,y,z) \in \{0,1\}^{n+m+l} \text{ s.t. } (x,y,z) \neq (\bar x, \bar y, \bar z), \ \obj(x,y,z) \le \obj(\bar x, \bar y, \bar z)] \\
& \le \sum_{\substack{(x,y,z) \in \{0,1\}^{n+m+l} \\ (x,y,z) \neq (\bar x, \bar y, \bar z)}} \prob [\obj(x,y,z) \le \obj(\bar x, \bar y, \bar z)] \\
& \le \sum_{\substack{(x,y,z) \in \{0,1\}^{n+m+l} \\ (x,y,z) \neq (\bar x, \bar y, \bar z)}} \exp \pare{- 2 \delta_{(x,y,z)} \pare{1/2 - \pr}^2} \\
& \le 2^{n+m+l} \exp \pare{- 2 \min\{r_{\bar x} n r_{\bar y} m, r_{\bar x} n r_{\bar z} l, r_{\bar y} m r_{\bar z} l\} \pare{1/2 - \pr}^2}.
\end{align*}
In the first inequality we used the union bound.
In the second inequality we used part~\ref{lem probs exp good} of~\cref{lem probs}.
In the third inequality we used part~\ref{lem card all}  of~\cref{lem card} and the fact that there are $2^{n+m+l}$ vectors $(x,y,z) \in \{0,1\}^{n+m+l}$.

We now show that the last expression goes to zero as $n, m, l \to \infty$.
For notational simplicity we define $\mu(n,m,l) := \min\{r_{\bar x} n r_{\bar y} m, r_{\bar x} n r_{\bar z} l, r_{\bar y} m r_{\bar z} l\}$.
We have
\begin{align*}
& \lim_{n, m, l \to \infty} 2^{n+m+l} \exp \pare{- 2 \mu(n,m,l) \pare{1/2 - \pr}^2} \\
& \quad = \lim_{n, m, l \to \infty} \exp \pare{\log \pare{ 2^{n + m + l} {\exp \pare{- 2 \mu(n,m,l) \pare{1/2 - \pr}^2}} }} \\
& \quad = \lim_{n, m, l \to \infty} \exp \pare{\log \pare{ 2^{n + m + l}} + \log \pare{ {\exp \pare{- 2 \mu(n,m,l) \pare{1/2 - \pr}^2}} }} \\
& \quad = \lim_{n, m, l \to \infty} \exp \pare{(n + m + l) \log 2 - 2 \mu(n,m,l) \pare{1/2 - \pr}^2} \\
& \quad = \exp \pare{ \lim_{n, m, l \to \infty} \pare{(n + m + l) \log 2 - 2 \mu(n,m,l) \pare{1/2 - \pr}^2}}.
\end{align*}
By assumption, the last limit is $- \infty$.
Hence, the original limit is $\exp(- \infty) = 0$.

We have shown that the probability that \eqref{pr comb} does not recover the ground truth goes to zero as $n,m,l \to \infty$.
Therefore, the probability that \eqref{pr comb} recovers the ground truth goes to one as $n,m,l \to \infty$ and so \eqref{pr comb} recovers the ground truth with high probability.
\end{prf}

We are now ready to present the proof of \cref{th it ub}.
The key difference with \cref{prop it ub} is that, in \cref{th it ub}, we assume that $\pr, r_{\bar x}, r_{\bar y}, r_{\bar z}$ are constants, and we obtain a simpler condition in terms of the growth rates of $n,m,l$.

\begin{prfc}[of \cref{th it ub}]
For notational simplicity we define $\mu(n,m,l) := \min\{r_{\bar x} n r_{\bar y} m, r_{\bar x} n r_{\bar z} l, r_{\bar y} m r_{\bar z} l\}$.
Consider the limit in the statement of \cref{prop it ub}.
We have
\begin{align*}
& \lim_{n, m, l \to \infty} \pare{2 \mu(n,m,l) \pare{1/2 - \pr}^2 - (n + m + l) \log 2} \\
& \quad = \lim_{n, m, l \to \infty} 2 \mu(n,m,l) \pare{1/2 - \pr}^2 \pare{1 - \frac{(n + m + l) \log 2}{2 \mu(n,m,l) \pare{1/2 - \pr}^2}} \\
& \quad = \pare{\lim_{n, m, l \to \infty} 2 \mu(n,m,l) \pare{1/2 - \pr}^2} \cdot \pare{\lim_{n, m, l \to \infty} \pare{1 - \frac{(n + m + l) \log 2}{2 \mu(n,m,l) \pare{1/2 - \pr}^2}}}.
\end{align*}
Consider the first limit in the last expression.
Since $r_{\bar x}$, $r_{\bar y}$, $r_{\bar z}$ are constants in $(0,1]$ and $p$ is a constant in $[0,1/2)$, this first limit is $\infty$.
Consider now the second limit.
Using our limit assumption and the fact that $r_{\bar x}$, $r_{\bar y}$, $r_{\bar z}$ are constants in $(0,1]$ and $\pr$ is a constant in $[0,1/2)$, we have
\begin{align*}
\lim_{n, m, l \to \infty} \frac{(n + m + l) \log 2}{2 \mu(n,m,l) \pare{1/2 - \pr}^2} = 0.
\end{align*}
Hence the second limit is $1$.
The original limit is then equal to $\infty$, and the corollary follows from \cref{prop it ub}.
\end{prfc}

\section{Recovery proof for the standard LP}
\label{sec SL recovery}

The main goal of this section is to prove \cref{th LP1}.
To this end, we first obtain deterministic sufficient conditions for recovery. Subsequently, we study the semi-random corruption model.

\subsection{Deterministic recovery guarantee}
In the following, we present deterministic recovery guarantees.
We first focus on the  special case in which the input tensor $\G$ is not corrupted; \ie $\F = \emptyset$.
Subsequently, we consider the problem with corrupted inputs.

\begin{proposition}
\label{unique1 nonoise}
Let $\bar x \in \{0,1\}^n$, $\bar y \in \{0,1\}^m$, $\bar z \in \{0,1\}^l$.
Define $\bar \W = (\bar w_{ijk}) := \bar x \otimes \bar y \otimes \bar z  \in \{0,1\}^{n \times m \times l}$ and let $\G := \bar \W$.
Then $(\bar x, \bar y, \bar z, \bar \W)$ is an optimal solution of Problem~\eqref{LP1}.
Furthermore, $(\bar x, \bar y, \bar z, \bar \W)$ is the unique optimal solution of Problem~\eqref{LP1} if and only if $\bar x, \bar y, \bar z \neq 0$.
\end{proposition}

\begin{prf}
First, we show that $(\bar x, \bar y, \bar z, \bar \W)$ is an optimal solution of Problem~\eqref{LP1}.
Since $\G = \bar \W$, in Problem~\eqref{LP1}, we have $S_0 = \N$ and $S_1 = \P$.
Let $(x, y, z, \W)$ be a feasible solution to Problem~\eqref{LP1}.
The objective value of this solution is nonnegative, since $w_{ijk} \ge 0$ for every $(i,j,k) \in \N$ and $w_{ijk} \le 1$ for every $(i,j,k) \in \P$.
Since $(\bar x, \bar y, \bar z, \bar \W)$ is a feasible solution to Problem~\eqref{LP1} with objective value zero, it is then an optimal solution to Problem~\eqref{LP1}.

In the rest of the proof we show the ``if and only if'' in the statement.
First, we prove the ``only if'': We assume that at least one among $\bar x, \bar y, \bar z$ is the zero vector, and we show that $(\bar x, \bar y, \bar z, \bar \W)$ is not the unique optimal solution of Problem~\eqref{LP1}.
Assume $\bar x=0$.
It is then simple to check that $(\bar x, 1 - \bar y, \bar z, \bar \W)$ is another optimal solution to Problem~\eqref{LP1}.
The cases $\bar y=0$ and $\bar z=0$ are symmetric, so this concludes the proof of the ``only if''.

In the remainder of the proof we show the ``if'' in the statement: We assume $\bar x, \bar y, \bar z \neq 0$ and show that $(\bar x, \bar y, \bar z, \bar \W)$ is the unique optimal solution of Problem~\eqref{LP1}.
Since the objective value of $(\bar x, \bar y, \bar z, \bar \W)$ is zero, it suffices to show that, if $(x, y, z, \W)$ is a feasible solution to Problem~\eqref{LP1} with objective value zero, then we have $(x, y, z, \W) = (\bar x, \bar y, \bar z, \bar \W)$.
Since the objective value of $(x, y, z, \W)$ is zero, we have $w_{ijk} = 0$ for every $(i,j,k) \in \N$ and $w_{ijk} = 1$ for every $(i,j,k) \in \P$. Thus $\W = \bar \W$.

We now show that, for every $i \in [n]$, $\bar x_i =1$ implies $x_i =1$.
Since $\bar y, \bar z \neq 0$, there exist $j \in [m]$, $k \in [l]$ such that $\bar y_j = 1$, $\bar z_k = 1$.
Then we have $(i,j,k) \in \P$ and so $w_{ijk} = 1$.
Constraints \eqref{cp}, \eqref{ub} then imply $x_i=1$.
Next, we show that, for every $i \in [n]$, $\bar x_i =0$ implies $x_i =0$.
Since $\bar y, \bar z \neq 0$, there exist $j \in [m]$, $k \in [l]$ such that $\bar y_j = 1$, $\bar z_k = 1$.
Then we have $(i,j,k) \in \N$ and so $w_{ijk} = 0$.
Constraints \eqref{cn}, \eqref{ub} then imply $x_i = 0$.
We have shown that $x = \bar x$.
Symmetrically, we obtain $y = \bar y$ and $z = \bar z$.
\end{prf}

Henceforth, assume that $\F \neq \emptyset$.
In the following, we first present a sufficient condition under which an optimal solution of Problem~\eqref{LP1} coincides with the ground truth. Next, we investigate the question of uniqueness.
%
For every $i \in [n]$ and $r,s \in \{0,1\}$, define
\begin{align*}
& T^{x,i}_{rs}:=\Big|(j,k) \in [m]\times[l]: \bar y_{j}=r, \bar z_{k}=s, (i,j,k) \in \T \Big|,\\
& F^{x,i}_{rs}:=\Big|(j,k)\in [m]\times[l]: \bar y_{j}=r, \bar z_{k}=s, (i,j,k) \in \F \Big|.
\end{align*}
Parameters $T^{y,j}_{rs}$, $F^{y,j}_{rs}$ for all $j \in [m]$, and
$T^{z,k}_{rs}$, $F^{z,k}_{rs}$ for all $k \in [l]$ are similarly defined.
Finally, define $T^{x,i} = T^{x,i}_{00}+ T^{x,i}_{10}+ T^{x,i}_{01}+T^{x,i}_{11}$
and $F^{x,i} = F^{x,i}_{00}+ F^{x,i}_{10}+ F^{x,i}_{01}+F^{x,i}_{11}$. Parameters
$T^{y,j}$, $F^{y,j}$, $T^{z,k}$, $F^{z,k}$ are similarly defined.

\begin{proposition}
\label{deter1}
Let $\bar x \in \{0,1\}^n$, $\bar y \in \{0,1\}^m$, $\bar z \in \{0,1\}^l$ and define $\bar \W = (\bar w_{ijk}) := \bar x \otimes \bar y \otimes \bar z  \in \{0,1\}^{n \times m \times l}$.
Suppose that $\F \neq \emptyset$.
Then $(\bar x, \bar y, \bar z, \bar \W)$ is an optimal solution of Problem~\eqref{LP1}, if the following conditions are satisfied:
\begin{enumerate}
    \item
    \label{c1-1}
    For each $i \in [n]$ with $\bar x_i = 0$, we have $F^{x,i} > 0$, for each $j \in [m]$ with $\bar y_j = 0$, we have $F^{y,j} > 0$, and for each $k \in [l]$ with $\bar z_k = 0$, we have $F^{z,k} > 0$.
    \item \label{c1-2}
    For each $i \in [n]$ with $\bar x_i =0$, we have
    \begin{equation}
    \label{nc1}
     T^{x,i}_{11} \geq \frac{1}{3}F^{x,i}_{00}+\frac{1}{2}(F^{x,i}_{01}+F^{x,i}_{10})
     +F^{x,i}_{11}.
   \end{equation}
   Symmetrically,~\eqref{nc1} holds when replacing $(x,i,n)$ by
   $(y,j,m)$ and $(z,k,l)$, respectively.
\item \label{c1-3}
For each $i \in [n]$ with $\bar x_i =1$, we have
\begin{align}
\label{nc4}
\begin{split}
\frac{1}{3}T^{x,i}_{11}
& \geq F^{x,i}_{11}+\sum_{\substack{(j,k):(i,j,k) \in \T: \\ \bar y_j= 0, \bar z_k =1}}{ \frac{1}{T^{y,j}_{11}}\Big(\frac{1}{3}F^{y,j}_{00} + \frac{1}{2}(F^{y,j}_{01}+F^{y,j}_{10})+F^{y,j}_{11}\Big)}\\
& \quad +
\sum_{\substack{(j,k):(i,j,k) \in \T: \\ \bar y_j= 1, \bar z_k =0}}
\frac{1}{T^{z,k}_{11}}\Big(\frac{1}{3}F^{z,k}_{00} + \frac{1}{2}(F^{z,k}_{01}+F^{z,k}_{10})+F^{z,k}_{11}\Big).
\end{split}
\end{align}
   Symmetrically,~\eqref{nc4} holds when switching $(x,i,n)$ by
   $(y,j,m)$ and switching $(x,i,n)$ by $(z,k,l)$.
\end{enumerate}
\end{proposition}

\begin{prf}
We start by constructing the dual of Problem~\eqref{LP1}.
Define dual variables $\lambda^x_{ijk}$, $\lambda^y_{ijk}$, $\lambda^z_{ijk}$ for all $(i,j,k) \in S_1$ associated with the first, the second and the third set of constraints in~\eqref{cp}, respectively.  Define $\mu^1_{ijk}, \mu^2_{ijk}$ for all $(i,j,k) \in S_0$ associated with the first and the second set of constraints in~\eqref{cn}, respectively.
Finally, define $u^x_i$ (resp. $l^x_i$) for all $i \in [n]$, $u^y_j$ (resp. $l^y_j$)
for all $j \in [m]$, and $u^z_k$ (resp. $l^z_k$) for all $k \in [l]$, associated with $x_i \leq 1$ (resp. $-x_i \leq 0$), $y_j \leq 1$ (resp. $-y_j \leq 0$), and $z_k\leq 1$ (resp. $-z_k \leq 0$) respectively.
It then follows that the dual of Problem~\eqref{LP1} is given by
\begin{align}\label{dual}
\tag{sD}
\max \quad & |S_1|-2 \sum_{(i,j,k) \in S_0}{\mu^2_{ijk}}-\sum_{i \in [n]}{u^x_i}-\sum_{j \in [m]}{u^y_j}-\sum_{k \in [l]}{u^z_k}\nonumber \\
\st  \quad & \sum_{(j,k):(i,j,k) \in S_0}{\mu^2_{ijk}} + u^x_i-l^x_i = \sum_{(j,k):(i,j,k) \in S_1}{\lambda^x_{ijk}} , \qquad \forall i \in [n]\label{e1}\\
&  \sum_{(i,k):(i,j,k) \in S_0}{\mu^2_{ijk}} + u^y_j - l^y_j=  \sum_{(i,k):(i,j,k) \in S_1}{\lambda^y_{ijk}}, \qquad \forall j \in [m]\label{e2}\\
& \sum_{(i,j):(i,j,k) \in S_0}{\mu^2_{ijk}} + u^z_k - l^z_k=\sum_{(i,j):(i,j,k) \in S_1}{\lambda^z_{ijk}} , \qquad \forall k \in [l]\label{e3}\\
& \lambda^x_{ijk} + \lambda^y_{ijk} + \lambda^z_{ijk} = 1, \qquad \forall (i,j,k) \in S_1 \label{e4}\\
& \mu^1_{ijk} + \mu^2_{ijk} = 1 , \qquad \forall (i,j,k) \in S_0\label{e4m}\\
& \lambda^x_{ijk}\geq 0, \; \lambda^y_{ijk}\geq 0, \; \lambda^z_{ijk}\geq 0, \qquad \forall (i,j,k) \in S_1\label{e5}\\
& \mu^1_{ijk}\geq 0, \mu^2_{ijk} \geq 0 , \qquad \forall (i,j,k) \in S_0 \label{e6}\\
& l^x_i \geq 0, \; \forall i \in [n], \ l^y_j \geq 0, \; \forall j \in [m], \ l^z_k \geq 0, \; \forall k \in [l] \label{usl1}\\
& u^x_i \geq 0, \; \forall i \in [n], \ u^y_j \geq 0, \; \forall j \in [m], \ u^z_k \geq 0, \; \forall k \in [l] \label{e7}.
\end{align}
To prove the optimality of $(\bar x, \bar y, \bar z, \bar \W)$, it suffices to construct a dual feasible point
$(\bar\lambda^x, \bar\lambda^y, \bar\lambda^z, \bar\mu^1,$ $\bar\mu^2, \bar u^x, \bar u^y, \bar u^z)$ of Problem~\eqref{dual} that satisfies complementary slackness. First, to satisfy~\eqref{usl1} we set
$l^x_i =0$ for all $i \in [n]$, $l^y_j= 0$ for all $j \in [m]$ and $l^z_k = 0$ for all $k \in [l]$.
By complementary slackness, we have:
\begin{enumerate}
\item [(I)] For each $(i,j,k) \in \F \cap \N$: (i) if $\bar x_i =1$, we have $\bar \lambda^x_{ijk} = 0$, (ii) if $\bar y_j =1$, we have $\bar \lambda^y_{ijk} = 0$, and (iii) if
$\bar z_k =1$, we have $\bar \lambda^z_{ijk} = 0$.
\item [(II)] For each $(i,j,k) \in \F \cap \P$, we have $\bar \mu^1_{ijk} = 0$; in this case, by~\eqref{e4m}, we get $\bar \mu^2_{ijk} = 1$.
\item  [(III)] For each $(i,j,k) \in \T \cap \N$ with $\bar x_i= \bar y_j = 0$ or  $\bar x_i=\bar z_k = 0$ or $\bar y_j = \bar z_k = 0$, we have $\bar \mu^2_{ijk} = 0$;
in this case by~\eqref{e4m}, we get  $\bar \mu^1_{ijk} = 1$.
\item [(IV)] For each $i \in [n]$, with $\bar x_i = 0$, we have $\bar u^x_i = 0$; for each $j \in [m]$, with $\bar y_j = 0$, we have $\bar u^y_j = 0$;
for each $k \in [l]$, with $\bar z_k = 0$, we have $\bar u^z_k = 0$.
\end{enumerate}
In order to satisfy constraints~\eqref{e4} and~\eqref{e5}, we choose $\bar\lambda^x, \bar\lambda^y, \bar\lambda^z$ as follows:
\begin{align}
\label{lamda}
\begin{split}
& \bar\lambda^x_{ijk} = \bar\lambda^y_{ijk} = \lambda^z_{ijk} = \frac{1}{3},\; \forall (i,j,k) \in \T \cap \P, \\
& \bar\lambda^x_{ijk} = \bar\lambda^y_{ijk} = \lambda^z_{ijk} = \frac{1}{3}, \;\forall (i,j,k) \in \F \cap \N, \; {\rm with} \; \bar x_i = \bar y_j = \bar z_k =0, \\
& \bar\lambda^x_{ijk} = \bar\lambda^y_{ijk} = \frac{1}{2}, \; \forall (i,j,k) \in \F \cap \N, \; {\rm with} \; \bar x_i = \bar y_j = 0, \; \bar z_k =1, \\
& \bar\lambda^x_{ijk} = \bar\lambda^z_{ijk} = \frac{1}{2}, \; \forall (i,j,k) \in \F \cap \N, \; {\rm with} \; \bar x_i = \bar z_k = 0, \; \bar y_j =1, \\
& \bar\lambda^y_{ijk} = \bar\lambda^z_{ijk} = \frac{1}{2}, \; \forall (i,j,k) \in \F \cap \N, \; {\rm with} \; \bar y_j = \bar z_k = 0, \; \bar x_i =1, \\
& \bar\lambda^x_{ijk} = 1, \; \forall (i,j,k) \in \F \cap \N, \; {\rm with}\; \bar x_i =0, \; \bar y_j = \bar z_k = 1, \\
& \bar\lambda^y_{ijk} = 1, \; \forall (i,j,k) \in \F \cap \N, \; {\rm with}\; \bar y_j =0, \; \bar x_i = \bar z_k = 1, \\
& \bar\lambda^z_{ijk} = 1, \; \forall (i,j,k) \in \F \cap \N, \; {\rm with}\; \bar z_k =0,\; \bar x_i = \bar y_j = 1.
\end{split}
\end{align}
Moreover, we let
\begin{align}\label{mu}
&\bar\mu^2_{ijk}  = \alpha_i, \qquad \forall (i,j,k) \in \T \cap \N, \; {\rm with} \; \bar x_i = 0,\; \bar y_j = \bar z_k =1,\nonumber\\
&\bar\mu^2_{ijk}  = \beta_j, \qquad \forall  (i,j,k) \in \T \cap \N, \; {\rm with} \; \bar y_j = 0,\; \bar x_i = \bar z_k =1,\\
&\bar\mu^2_{ijk}  = \gamma_k, \qquad \forall (i,j,k) \in \T \cap \N, \; {\rm with} \; \bar z_k = 0,\; \bar x_i = \bar y_j =1,\nonumber
\end{align}
where parameters $\alpha_i, \beta_j, \gamma_k$ are to be determined later. By constraints~\eqref{e4m} we have
$\bar\mu^1_{ijk}  = 1-\alpha_i$, for all $(i,j,k) \in \T \cap \N$ with
$\bar x_i = 0$, $\bar y_j = \bar z_k =1$,
$\bar\mu^1_{ijk}  = 1-\beta_j$, for all  $(i,j,k) \in \T \cap \N$ with
$\bar y_j = 0$, $\bar x_i = \bar z_k =1$, and
$\bar\mu^1_{ijk}  = 1-\gamma_k$, for all $(i,j,k) \in \T \cap \N$, with
$\bar z_k = 0$, $\bar x_i = \bar y_j =1$.
Hence to satisfy constraints~\eqref{e6} we impose
$$0 \leq \alpha_i \leq 1,\qquad
0 \leq \beta_j \leq 1, \qquad
0 \leq \gamma_k \leq 1.$$
Substituting~\eqref{lamda} and~\eqref{mu} in constraints~\eqref{e1}-\eqref{e3}, the following cases arise:
\begin{itemize}[leftmargin=*]
\item For each $i \in [n]$ with $\bar x_i = 0$, constraints~\eqref{e1} simplify to
\begin{equation*}
\sum_{\substack{(j,k):(i,j,k) \in \T: \\ \bar y_j = \bar z_k =1}}{\alpha_i}  = \sum_{\substack{(j,k):(i,j,k) \in \F: \\ \bar y_j = \bar z_k = 0}}
{\frac{1}{3}}
+ \sum_{\substack{(j,k):(i,j,k) \in \F: \\ (\bar y_j =1, \bar z_k = 0) \lor (\bar y_j =0, \bar z_k = 1)}}{\frac{1}{2}}
+ \sum_{\substack{(j,k):(i,j,k) \in \F: \\ \bar y_j =\bar z_k = 1}}{1}.
\end{equation*}
By Condition~\ref{c1-1}, we have $F^{x,i} >0$, \ie the right-hand side of the above inequality is positive. By inequality~\eqref{nc1} of Condition~\ref{c1-2}, we have $T^{x,i}_{11} > 0$, \ie the left-hand side of the above inequality is positive.  Hence we obtain:
\begin{equation}\label{alpha}
\alpha_i  = \frac{1}{T^{x,i}_{11}}\Big(\frac{1}{3}F^{x,i}_{00}+\frac{1}{2}(F^{x,i}_{01}+F^{x,i}_{10})+F^{x,i}_{11}\Big).
\end{equation}
Clearly, $\alpha_i \geq 0$, hence it suffices to have $\alpha_i \leq 1$, which can be equivalently written as inequality~\eqref{nc1}.
Similarly, for each $j \in [m]$ with $\bar y_j = 0$,
by Condition~\ref{c1-1} we have $F^{y,j} > 0$ and by symmetric counterpart of inequality~\eqref{nc1} of Condition~\ref{c1-2} we have $T^{y,j}_{11} > 0$. Hence to satisfy constraints~\eqref{e2},
we let
\begin{equation}\label{beta}
\beta_j = \frac{1}{T^{y,j}_{11}}\Big(\frac{1}{3}F^{y,j}_{00} + \frac{1}{2}(F^{y,j}_{01}+F^{y,j}_{10})+F^{y,j}_{11}\Big).
\end{equation}
It then follows that the constraint $\beta_j \leq 1$ can be equivalently written as a symmetric counterpart of inequality~\eqref{nc1}.
Finally, for each $k \in [l]$ with $\bar z_k = 0$,
by Condition~\ref{c1-1}, we have $F^{z,k} > 0$ and by symmetric counterpart of inequality~\eqref{nc1} of Condition~\ref{c1-2} we have $T^{z,k}_{11} > 0$. Hence to satisfy constraints~\eqref{e3}, we let
\begin{equation}\label{gamma2}
\gamma_k = \frac{1}{T^{z,k}_{11}}\Big(\frac{1}{3}F^{z,k}_{00} + \frac{1}{2}(F^{z,k}_{01}+F^{z,k}_{10})+F^{z,k}_{11}\Big).
\end{equation}
It then follows that $\gamma_k \leq 1$ can be equivalently
written as a symmetric counterpart inequality~\eqref{nc1}.
\item For each $i \in [n]$ with $\bar x_i = 1$, constraints~\eqref{e1} simplify to
\begin{equation*}
\bar u^x_i =
\sum_{\substack{(j,k):(i,j,k) \in \T \\ \bar y_j = \bar z_k =1}} {\frac{1}{3}}-\sum_{\substack{(j,k):(i,j,k) \in \T: \\ \bar y_j= 0, \bar z_k =1}}{\beta_j}-
\sum_{\substack{(j,k):(i,j,k) \in \T: \\ \bar y_j= 1, \bar z_k =0}}{\gamma_k}
- \sum_{\substack{(j,k): (i,j,k) \in \F \\  \bar y_j = \bar z_k =1}}{1}.
\end{equation*}
Substituting for $\beta_j, \gamma_k$ using~\eqref{beta} and~\eqref{gamma2},
it follows that the constraint $\bar u^x_i \geq 0$ can be equivalently written as inequality~\eqref{nc4} of Condition~\ref{c1-3} in \cref{deter2}.
Similarly, substituting for $\alpha_i, \gamma_k$ using~\eqref{alpha} and~\eqref{gamma2} in equalities~\eqref{e2}, it follows that
for each $j \in [m]$ with $\bar y_j = 1$, the constraint $\bar u^y_j \geq 0$
can be written as a symmetric counterpart of inequality~\eqref{nc4} of Condition~\ref{c1-3}.
Finally, substituting for $\alpha_i, \beta_j$ using~\eqref{alpha} and~\eqref{beta} in equalities~\eqref{e3}, it follows that
for each $k \in [l]$ with $\bar z_k = 1$,  the constraint $\bar u^z_k \geq 0$
can be written as a symmetric counterpart of inequality~\eqref{nc4}  of Condition~\ref{c1-3}.
\end{itemize}
\end{prf}

It is important to note that Condition~\ref{c1-1} in the statement of Proposition~\ref{deter1} is not necessary and is only added to simplify the remaining conditions and the proof. As we will show shortly, for the fully-random corruption model, this condition is not restrictive as it always holds with high probability, provided that $p>0$.
%
%
%
We now provide a sufficient condition under which the ground truth
is the unique optimal solution of Problem~\eqref{LP1}.
To this end, we use Mangasarian's characterization of uniqueness of LP optimal solutions~\cite{Man79}:

\begin{proposition}[Part~(iv) of Theorem~2 in~\cite{Man79}]\label{mang}
Consider an LP whose feasible region is defined by $C x \leq d$.
Let $\bar x$ be an optimal solution of this LP and denote by $\bar u$ the vector of dual optimal solution.
Let $C_i$ denote
the $i$-th row of $C$. Define $K =\{i : C_i \bar x = d_i, \; \bar u_i > 0\}$, $L =\{i : C_i \bar x = d_i, \; \bar u_i = 0\}$. Let $C_K$ and $C_L$
be the
matrices whose rows are $C_i$, $i \in K$ and $C_i$, $i \in L$, respectively.
Then $\bar x$ is the unique optimal solution of the LP, if there exists no nonzero vector $x$ satisfying
\begin{equation}\label{uc}
C_K x = 0, \qquad C_L x \leq 0.
\end{equation}
\end{proposition}

Utilizing the above result, we are now ready to establish our uniqueness condition:
\begin{proposition}\label{unique1}
Suppose that all assumptions of Proposition~\ref{deter1} hold.
Moreover, suppose that inequalities~\eqref{nc1} and~\eqref{nc4} and their symmetric counterparts are strictly satisfied.
Then $(\bar x, \bar y, \bar z, \bar \W)$ is the unique optimal solution of Problem~\eqref{LP1}.
\end{proposition}
\begin{prf}
Let $(\bar x, \bar y, \bar z, \bar \W)$ be an optimal solution of Problem~\eqref{LP1}.
To prove the statement it suffices to show there is no nonzero 
$(x, y, z, \W)$ satisfying condition~\eqref{uc}.
We have:
\begin{itemize}
   \item [(i)]
   Since inequalities~\eqref{nc1} strictly hold, together with part~(III) of complementary slackness in the proof of Proposition~\ref{deter1}, we conclude that
   $\bar \mu^1_{ijk} > 0$ for all $(i,j,k) \in \T \cap \N$, implying $w_{ijk} = 0$ for all $(i,j,k) \in \T \cap \N$.

    \item [(ii)]
    Since inequalities~\eqref{nc4} strictly hold, we have $\bar u^x_i > 0$
    for all $i \in [n]$ with $\bar x_i = 1$. This in turn implies that $x_i = 0$ for all $i \in [n]$ with $\bar x_i = 1$. By symmetry we conclude that $y_j = 0$ for all $j \in [m]$ with $\bar y_j = 1$ and
    $z_k = 0$ for all $k \in [l]$ with $\bar z_k = 1$.

    \item [(iii)]
    By part~(II) of complementary slackness in the proof of Proposition~\ref{deter1}, we have $\bar \mu^2_{ijk} > 0$ for all $(i,j,k) \in \F \cap \P$, implying $w_{ijk} = x_i + y_j + z_k$ for all $(i,j,k) \in \F \cap \P$. By part~(ii) above this implies that $w_{ijk} = 0$ for all $(i,j,k) \in \F \cap \P$.

    \item [(iv)]
    By~\eqref{lamda} we have $\bar \lambda^x_{ijk} > 0$ for all $(i,j,k) \in \T \cap \P$, implying $w_{ijk} = x_i$ in this case. By part~(ii) above, we conclude that $w_{ijk} = 0$ for all $(i,j,k) \in \T \cap \P$.

    \item [(v)]
    By assumption~\ref{c1-1} and inequality~\eqref{nc1} of Proposition~\ref{deter1}, for each $i \in [n]$ with $\bar x_i=0$, we have
    $T^{x,i}_{11} \geq 1$; that is, Problem~\eqref{LP1} contains a constraint of the form $w_{ijk} \geq x_i + y_j + z_k -2$ with $(i,j,k) \in \T$ and  $\bar x_i=0$, $\bar y_j = \bar z_k = 1$. By~\eqref{mu} and assumption~~\ref{c1-1}, we have $\bar \mu^2_{ijk} > 0$.
    Therefore, $w_{ijk} = x_i + y_j + z_k = x_i$, where the second equality follows from part~(ii) above. By part~(i) we have $w_{ijk} = 0$; hence, we conclude that $x_i =0$  for any $i \in [n]$ with $\bar x_i=0$.
    By symmetry, $y_j =0$  for any $j \in [m]$ with $\bar y_j=0$
    and $z_k =0$  for any $k \in [l]$ with $\bar z_k=0$.

    \item [(vi)] By~\eqref{lamda} for any $(i,j,k) \in \F \cap \N$ with $\bar x_i = 0$, we have $\bar \lambda^x_{ijk} > 0$. This in turn implies that we must have $w_{ijk} = x_i$. By part~(v) above we have $x_i = 0$, implying $w_{ijk} = 0$. By symmetry, it follows that $w_{ijk} = 0$ for all $(i,j,k) \in \F \cap \N$.
\end{itemize}
From parts~(i)-(vi) we conclude there is no nonzero $(x,y,z,\W)$ satisfying~\eqref{uc}.
\end{prf}

\subsection{Recovery under the random corruption model}

We now consider the semi-random corruption model and prove~\cref{th LP1}, which provides a sufficient condition in terms of $p, r_{\bar x}, r_{\bar y}, r_{\bar z}$ under which the standard LP recovers the ground truth with high probability.
To this end, we define the following random variables.
For each $i \in [n]$, $j \in [m]$, $k \in [l]$ and $r,s \in \{0,1\}$, denote by $t^{x,i}_{jk \rightarrow rs}$ (resp. $f^{x,i}_{jk \rightarrow rs}$) a random variable whose value equals $1$, if $\bar y_j = r$, $\bar z_k = s$, $(i,j,k) \in \T$
(resp. $(i,j,k) \in \F$), and equals $0$, otherwise.
Random variables $t^{y,j}_{ik \rightarrow rs}$, $f^{y,j}_{ik \rightarrow rs}$, $t^{z,k}_{ij \rightarrow rs}$, $f^{z,k}_{ij \rightarrow rs}$ for all
$i\in [n], j \in [m], k \in [l]$ are similarly defined.
In the remainder of the paper, we denote by $n_{\bar x}$, $n_{\bar y}$, and $n_{\bar z}$ the number of ones in binary vectors $\bar x$, $\bar y$, $\bar z$, i.e., $n_{\bar x} := n r_{\bar x}$, $n_{\bar y} := m r_{\bar y}$, and $n_{\bar z} := l r_{\bar z}$.
For ease of notation, from now on, when we sum over index sets $[n]$, $[m]$, or $[l]$, we omit the index set, with the understanding that indices $i,i'$ are summed over $[n]$, indices $j,j'$ are summed over $[m]$, and indices $k,k'$ are summed over $[l]$.

\begin{prfc}[of \cref{th LP1}]
By~\cref{robustLPs} it suffices to prove the statement under the fully-random corruption model.
First, let us we consider the case where the input tensor is not corrupted; \ie $p = 0$. From Proposition~\ref{unique1 nonoise} it follows that if $r_{\bar x}, r_{\bar y}, r_{\bar z}$ are positive, the standard LP recovers the ground truth.

Henceforth, suppose that $p > 0$.
This assumption implies, in particular that $\F \neq \emptyset$ with high probability.
Denote by $A^0$ the event that Condition~\ref{c1-1} in~\cref{deter1} is satisfied.
Denote by $A^1$, the event that all inequalities of the form~\eqref{nc1} and symmetric counterparts strictly hold.
Moreover, denote by
$A^2$ the event that all inequalities of the form~\eqref{nc4} and symmetric counterparts strictly hold. Denote by $A_{\rm recovery}$ the event that the standard LP recovers the ground truth. Then, by Proposition~\ref{unique1} we have  $\prob[A_{\rm recovery}] \geq \prob [A^0 \cap A^1 \cap A^2]$.
Since $A_{\rm recovery}$ is the intersection of a constant number of
events $A^i$, to establish recovery with high probability, it suffices to prove that each $A^i$, $i \in \{0,1,2\}$
occurs with high probability.

\begin{claim}\label{cl1}
Event $A^0$ occurs with high probability.
\end{claim}

\begin{cpf}
We have $A^0 = A^0_1 \cap A^0_2 \cap A^0_3$, where
the event $A^0_1$ occurs if
\begin{equation}\label{rr1}
\frac{1}{ml} \sum_{j,k}{f_{jk}^{x,i}} > 0, \qquad \forall i \in [n],
\end{equation}
where $f_{jk}^{x,i} = 1$, if $(i,j,k) \in \F$ and $f^{x,i}_{jk} = 0$,
otherwise. The event $A^0_2$ occurs if
$
\frac{1}{nl} \sum_{i,k}{f_{ik}^{y,j}} > 0,$ $\forall j \in [m],
$
where $f_{ik}^{y,j} = 1$, if $(i,j,k) \in \F$ and $f^{y,j}_{ik} = 0$, otherwise. The event $A^0_3$ occurs if
$
\frac{1}{nm} \sum_{i,j}{f_{ij}^{z,k}} > 0,$ $\forall k \in [l],
$
where $f_{ij}^{z,k} = 1$, if $(i,j,k) \in \F$ and $f^{z,k}_{ij} = 0$, otherwise.
We show that event $A^0_1$ occurs with high probability. Using a similar line of arguments, it follows that $A^0_2$ and $A^0_3$ occur with high probability. Denote by $\epsilon$ the expected value of the left-hand side of inequality~\eqref{rr1}. Then:
$$
\epsilon := \avg\Big[\frac{1}{ml} \sum_{j,k}{f_{jk}^{x,i}}\Big] = \frac{1}{ml} \sum_{j,k}{\avg\big[f_{jk}^{x,i}\big]} = p > 0,
$$
where the inequality follows by assumption.
Then:
\begin{align*}
    \prob[A^0_1]
    & = \prob \Big[\bigcap_{i=1}^n \Big\{\frac{1}{ml} \sum_{j,k}{f_{jk}^{x,i}} > 0\Big\}\Big]
     \geq \prob \Big[\bigcap_{i=1}^n \Big\{\Big|\frac{1}{ml} \sum_{j,k}{f_{jk}^{x,i}} -\avg \Big[\frac{1}{ml} \sum_{j,k}{f_{jk}^{x,i}}\Big]\Big| \leq \epsilon\Big\}\Big]\\
     & \geq 1-\sum_{i=1}^n{\prob \Big[\Big|\frac{1}{ml} \sum_{j,k}{f_{jk}^{x,i}} -\avg \Big[\frac{1}{ml} \sum_{j,k}{f_{jk}^{x,i}}\Big]\Big| > \epsilon\Big]}
      \geq 1 - 2 n \exp(-2ml\epsilon^2),
\end{align*}
where the first inequality follows by set inclusion, the second inequality follows by taking the union bound and the last inequality follows from the application of Hoeffding's inequality since the random variables $0 \leq f^{x,i}_{jk} \leq 1$ for all $i\in [n]$, $j\in [m]$, $k \in [l]$ are independent. The proof then follows since by assumption $\epsilon = p$ is a positive constant independent of $n,m,l$ and since the limit assumptions in the theorem imply that, as $n, m, l \to \infty$, we have $n \exp (-ml) \to 0$.
\end{cpf}

\begin{claim}\label{cl2}
Event $A^1$ occur with high probability.
\end{claim}
\begin{cpf}
Denote by $A^1_1$ the event that inequalities~\eqref{nc1} are strictly satisfied. By symmetry, to show that $A^1$ occurs with high probability, it suffices to show that $A^1_1$ occurs with high probability.
For each $i \in [n]$, define
\begin{equation*}
Y^{x,i}_{jk} := t^{x,i}_{jk \rightarrow 11}-\frac{1}{3}f^{x,i}_{jk \rightarrow 00}-\frac{1}{2} f^{x,i}_{jk \rightarrow 01}-\frac{1}{2} f^{x,i}_{jk \rightarrow 10}- f^{x,i}_{jk \rightarrow 11}.
\end{equation*}
It then follows that event $A^1$ occurs if
\begin{equation}\label{frp}
\frac{1}{ml} \Big(\sum_{j,k} {Y^{x,i}_{jk}}\Big) > 0, \qquad \forall i \in [n]\; {\rm with} \; \bar x_i = 0.
\end{equation}
Denote by $\epsilon$ the expected value of the left-side of~\eqref{frp}. For any $i \in [n]$ with $\bar x_i = 0$, we have:
\begin{align}\label{ep}
\epsilon & = \frac{1}{ml}\sum_{j,k}\Big(\avg[t^{x,i}_{jk \rightarrow 11}]-\frac{1}{3}\avg[f^{x,i}_{jk \rightarrow 00}]-\frac{1}{2}\avg[f^{x,i}_{jk \rightarrow 01}]-\frac{1}{2} \avg[f^{x,i}_{jk \rightarrow 10}]-\avg[f^{x,i}_{jk \rightarrow 11}]\Big)\nonumber\\
& = (1-p) r_{\bar y} r_{\bar z} - \frac{1}{3} p (1-r_{\bar y})(1-r_{\bar z})
-\frac{1}{2}p (r_{\bar y}(1-r_{\bar z})+r_{\bar z} (1-r_{\bar y}))-p r_{\bar y} r_{\bar z}.
\end{align}
Since $p > 0$, it follows that inequality $\epsilon > 0$ can be equivalently written as
\begin{equation}\label{cond1}
p < \frac{6 r_{\bar y}r_{\bar z}}{(8 r_{\bar y} r_{\bar z}+r_{\bar y}+r_{\bar z}+2)}.
\end{equation}
It is simple to check that the above condition is implied by inequality~\eqref{cond2}.
%
Notice that by symmetry, the inequalities obtained by replacing $(x,i,n)$ by $(y,j,m)$ in~\eqref{nc1} strictly hold in expectation if
$p < \frac{6 r_{\bar x} r_{\bar z}}{(8 r_{\bar x} r_{\bar z}+r_{\bar x}+r_{\bar z}+2)}$, and the inequalities obtained by replacing $(x,i,n)$ by $(z,k,l)$ in~\eqref{nc1} strictly hold in expectation if
$p < \frac{6 r_{\bar x} r_{\bar y}}{(8 r_{\bar x} r_{\bar y}+r_{\bar x}+r_{\bar y}+2)}$. Since by assumption $r_{\bar x} \geq r_{\bar y} \geq r_{\bar z}$, it can be checked that the latter two are implied by inequality~\eqref{cond1}.

Define $I_0 = \{i \in [n]: \bar x_i = 0\}$.
We now show that event $A_1^1$ occurs with high probability:
\begin{align*}
\prob[A_1^1]
& =\prob\Bigg[\bigcap_{i\in I_0} {\Bigg\{\frac{1}{ml}{\sum_{j,k}{Y^{x,i}_{jk}}} > 0\Bigg\}}\Bigg]
=\prob\Bigg[\bigcap_{i \in I_0}{\Bigg\{\frac{1}{ml}{\sum_{j,k}{Y^{x,i}_{jk}}}-\avg\Big[\frac{1}{ml}\sum_{j,k}{Y^{x,i}_{jk}}\Big]  > -\epsilon}\Bigg\}\Bigg]\\
& \geq \prob\Bigg[\bigcap_{i \in I_0} {\Bigg\{\Bigg|\frac{1}{ml}{\sum_{j,k}{Y^{x,i}_{jk}}}-\avg\Big[\frac{1}{ml}\sum_{j,k}{Y^{x,i}_{jk}}\Big]\Bigg|  < \epsilon \Bigg\} }\Bigg]\\
& \geq 1-\sum_{i\in I_0} {\prob\Bigg[\Bigg|{\frac{1}{ml}{\sum_{j,k}{Y^{x,i}_{jk}}}-\avg\Big[\frac{1}{ml}\sum_{j,k}{Y^{x,i}_{jk}}\Big]\Bigg| \geq \epsilon\Bigg]}}
\geq  1 - 2 n\exp\Big(-\frac{ml \epsilon^2}{2}\Big),
\end{align*}
where the first inequality follows from set inclusion, the second inequality follows from taking the union bound, and the last inequality follows from the application of Hoeffding's inequality using the fact that $Y^{x,i}_{jk}$, for all $i,j,k$ are independent random variables and $-1 \leq Y^{x,i}_{jk} \leq 1$. Since by assumption $p, r_{\bar y}, r_{\bar z}$ are positive constants, from~\eqref{ep} it follows that $\epsilon$ is a constant. Moreover, as we detailed above, by~\eqref{cond2} we have $\epsilon > 0$. The proof then follows since the limit assumptions in the theorem imply that, as $n, m, l \to \infty$, we have $n \exp (-ml) \to 0$.
\end{cpf}

\begin{claim}\label{cl3}
Events $A^2$ occurs with high probability.
\end{claim}

\begin{cpf}
Denote by $A^2_1$ the event that inequalities~\eqref{nc4} are strictly satisfied. By symmetry, to show that $A^2$ occurs with high probability, it suffices to show that $A^2_1$ occurs with high probability.
First notice that by Condition~2 of Proposition~\ref{deter1}, inequalities~\eqref{nc4} can be equivalently written as:
\begin{align*}
&\frac{1}{ml}(\frac{1}{3}T^{x,i}_{11}-F^{x,i}_{11})
- \frac{1}{m}\sum_{j \in [m]: \bar y_j= 0}\Bigg(\frac{1}{l} \sum_{\substack{k \in [l]: \bar z_k =1, \\ (i,j,k) \in \T}}{\min\left\{\frac{1}{T^{y,j}_{11}}\Big(\frac{1}{3}F^{y,j}_{00} + \frac{1}{2}(F^{y,j}_{01}+F^{y,j}_{10})+F^{y,j}_{11}\Big), \; 1\right\}}\Bigg)\nonumber\\
& \quad -\frac{1}{l}
\sum_{k \in [l]: \bar z_k =0} \Bigg(\frac{1}{m}\sum_{\substack{j \in [m]: \bar y_j =1 \\ (i,j,k) \in \T}} {\min\left\{
\frac{1}{T^{z,k}_{11}}\Big(\frac{1}{3}F^{z,k}_{00} + \frac{1}{2}(F^{z,k}_{01}+F^{z,k}_{10})+F^{z,k}_{11}\Big),1\right\}}\Bigg) \geq 0.
\end{align*}
We next define some random variables associated with the above inequality.
For each $i \in [n]$ with $\bar x_i =1$, and for each $(j,k) \in [m] \times [l]$, define
\begin{equation}\label{nuzu}
\nu^{y,j}_{ik} :=
\frac{\sum_{i',k'}\Big(\frac{1}{3} f^{y,j}_{i'k'\rightarrow 00}+\frac{1}{2} f^{y,j}_{i'k'\rightarrow 01}+\frac{1}{2} f^{y,j}_{i'k'\rightarrow 10}+f^{y,j}_{i'k'\rightarrow 11}\Big)}{\sum_{i',k'}{t^{y,j}_{i'k'\rightarrow 11}}},
\end{equation}
if $\bar y_j = 0$, $\bar z_k =1$, $(i,j,k) \in \T$, and define $\nu^{y,j}_{ik}: =0$, otherwise.
Since to define~\eqref{nuzu} we assume $(i,j,k) \in \T$, it follows that $t^{y,j}_{ik \rightarrow 11} =1$. Hence the denominator of~\eqref{nuzu} can be equivalently written as $1+\sum_{(i',k')\in [n] \times [l] \setminus (i, k)}{t^{y,j}_{i'k'\rightarrow 11}}$.
Subsequently, define
\begin{align}\label{defnu}
&\nu^{y,j}_{x,i} := \frac{1}{l}\sum_{k}{\nu^{y,j}_{ik}}, \qquad
\bar \nu^{y,j}_{x,i} := \frac{1}{l}\sum_{k}{\min\{\nu^{y,j}_{ik},1\}}.
\end{align}
Clearly
$0 \leq \bar \nu^{y,j}_{x,i} \leq r_{\bar z} \leq 1$ and $0 \leq \bar \nu^{y,j}_{z,k} \leq r_{\bar x} \leq 1$. The random variables $\nu^{z,k}_{ij}$, and $\bar \nu^{z,k}_{x,i}$ are similarly defined.
%
It then follows that the event $A^2_1$ occurs if, for every $i \in [n]$ with $\bar x_i = 1$,
\begin{equation}\label{rr4}
\frac{1}{ml}\sum_{j,k}{\Big(\frac{1}{3} t^{x,i}_{jk \rightarrow 11}-f^{x,i}_{jk \rightarrow 11}\Big)}
-\frac{1}{m}\sum_{j}{\bar \nu^{y,j}_{x,i}}
-\frac{1}{l}\sum_{k}{\bar \nu^{z,k}_{x,i}} >0.
\end{equation}
To prove the statement, it suffices to show that inequalities~\eqref{rr4}
hold with high probability. Denote by $\epsilon$ the expected value of the left-hand side of~\eqref{rr4}.
From the definition of $\nu^{y,j}_{x,i}, \nu^{z,k}_{x,i}$ given by~\eqref{defnu},
it follows that
\begin{align}
\epsilon & \geq  \avg\Big[\frac{1}{ml}\sum_{j,k}{\Big(\frac{1}{3} t^{x,i}_{jk \rightarrow 11}-f^{x,i}_{jk \rightarrow 11}\Big)}
-\frac{1}{m}\sum_{j}{\nu^{y,j}_{x,i}}
-\frac{1}{l}\sum_{k}{\nu^{z,k}_{x,i}}\Big] \nonumber\\
& = \frac{1}{ml} \sum_{j,k} \Big(\frac{1}{3}\avg[t^{x,i}_{jk \rightarrow 11}]- \avg[f^{x,i}_{jk \rightarrow 11}]-\avg[\nu^{y,j}_{ik}]- \avg[\nu^{z,k}_{ij}]\Big) := \bar \epsilon \label{eval}.
\end{align}
Hence $\epsilon > 0$, if $\bar \epsilon > 0$. In the following, we obtain a lower bound on $\bar \epsilon$.
To this end, we first obtain an upper bound on $\avg[\nu^{y,j}_{ik}]$.
Using a similar line of arguments, an upper bound on $\avg[\nu^{z,k}_{ij}]$ can be calculated.
Recall that by definition for each $i \in [n]$ with $\bar x_i = 1$, $\nu^{y,j}_{ik} =0$ unless
$\bar y_j = 0$, $\bar z_k =1$ and $(i,j,k) \in \T$
It then follows that
\begin{equation}\label{ff1}
    \avg[\nu^{y,j}_{ik}] = \avg\Big[\nu^{y,j}_{ik}\Big|(i,j,k) \in \T\Big] (1-p).
\end{equation}
Using $\sum_{i',k'}{(t^{y,j}_{i'k'\rightarrow 11} + f^{y,j}_{i'k'\rightarrow 11})} = n_{\bar x} n_{\bar z}$, we obtain
\begin{align}
\label{ff3}
   &\avg[\nu^{y,j}_{ik}|(i,j,k) \in \T]
   =\avg\Bigg[\frac{\sum_{i',k'}\Big(\frac{1}{3} f^{y,j}_{i'k'\rightarrow 00}+\frac{1}{2} f^{y,j}_{i'k'\rightarrow 01}+\frac{1}{2} f^{y,j}_{i'k'\rightarrow 10}\Big)+n_{\bar x} n_{\bar z}}{\Big(1+\sum_{(i',k')\in [n] \times [l] \setminus (i, k)}{t^{y,j}_{i'k'\rightarrow 11}}\Big)}-1\Bigg]\nonumber\\
   &=\avg\Bigg[\frac{1}{1+\sum_{(i',k')\in [n] \times [l] \setminus (i, k)}{t^{y,j}_{i'k'\rightarrow 11}}}\Bigg]
   \avg\Bigg[\sum_{i',k'}\Big(\frac{1}{3} f^{y,j}_{i'k'\rightarrow 00}+\frac{1}{2} f^{y,j}_{i'k'\rightarrow 01}+\frac{1}{2} f^{y,j}_{i'k'\rightarrow 10}\Big)+n_{\bar x} n_{\bar z}\Bigg]-1\nonumber\\
   &=\frac{1}{(1-p)n_{\bar x} n_{\bar z}} (1-p^{n_{\bar x} n_{\bar z}})\Big(\frac{1}{3}p (n-n_{\bar x})(l-n_{\bar z})+\frac{1}{2}p ((n-n_{\bar x})n_{\bar z}+n_{\bar x}(l-n_{\bar z}))+n_{\bar x} n_{\bar z}\Big)-1\nonumber\\
   &=\frac{1-p^{n_{\bar x} n_{\bar z}}}{1-p} \Big(\frac{p}{3} (\frac{1}{r_{\bar x} r_{\bar z}}+\frac{1}{2 r_{\bar x}}+\frac{1}{2 r_{\bar z}}-2)+1\Big)-1 \nonumber \\
   & \leq \frac{1}{1-p} \Big(\frac{p}{3} (\frac{1}{r_{\bar x} r_{\bar z}}+\frac{1}{2 r_{\bar x}}+\frac{1}{2 r_{\bar z}}-2)+1\Big)-1,
\end{align}
where in the third line we used the fact that for a binomial random variable $X$ with parameters $(N, \bar p)$, we have $\avg[\frac{1}{1+X}] = \frac{1}{(N+1)\bar p} (1-(1-\bar p)^{N+1})$. The last inequality follows since $0 \leq 1-p^{n_{\bar x} n_{\bar z}} \leq 1$ and
$\frac{1}{r_{\bar x} r_{\bar z}}+\frac{1}{2 r_{\bar x}}+\frac{1}{2 r_{\bar z}}-2 \geq 0$.
Substituting~\eqref{ff3} in~\eqref{ff1} yields:
\begin{equation}\label{ff4}
    \avg[\nu^{y,j}_{ik}] \leq \frac{p}{3} \Big(\frac{1}{r_{\bar x} r_{\bar z}}+\frac{1}{2 r_{\bar x}}+\frac{1}{2 r_{\bar z}}+1\Big).
\end{equation}
Using a similar line of arguments we obtain
\begin{equation}\label{ff5}
    \avg[\nu^{z,k}_{ij}] \leq \frac{p}{3} \Big(\frac{1}{r_{\bar x} r_{\bar y}}+\frac{1}{2 r_{\bar x}}+\frac{1}{2 r_{\bar y}}+1\Big).
\end{equation}
Substituting~\eqref{ff4} and~\eqref{ff5} in~\eqref{eval} yields:
\begin{align*}
\bar \epsilon & \geq
\frac{1}{3}(1-p) r_{\bar y} r_{\bar z}-p r_{\bar y} r_{\bar z}
-\frac{p}{3} \Big(\frac{1}{r_{\bar x} r_{\bar z}}+\frac{1}{2 r_{\bar x}}+\frac{1}{2 r_{\bar z}}+1\Big) (1-r_{\bar y})r_{\bar z} \\
& \quad -\frac{p}{3} \Big(\frac{1}{r_{\bar x} r_{\bar y}}+\frac{1}{2 r_{\bar x}}+\frac{1}{2 r_{\bar y}}+1\Big) (1-r_{\bar z})r_{\bar y} \\
&=\frac{r_{\bar y} r_{\bar z}}{6} \Big(2-\frac{p}{r_{\bar x}r_{\bar y} r_{\bar z}}( 4 r_{\bar x} r_{\bar y} r_{\bar z} +r_{\bar x} r_{\bar y} + r_{\bar x} r_{\bar z} - 2 r_{\bar y} r_{\bar z} + 2 r_{\bar x} -r_{\bar y} -r_{\bar z} +4)\Big):= \tilde \epsilon.
\end{align*}
Hence, if $\tilde \epsilon > 0$, we have $\epsilon > 0$.
Since $r_{\bar x}, r_{\bar y}, r_{\bar z}$ are positive, it can be checked that inequality $\tilde \epsilon > 0$ is equivalent to inequality~\eqref{cond2}.
Using a similar line of argument it can be shown that inequalities obtained by switching $(x,i,n)$ with $(y,j,m)$ in~\eqref{nc4} hold in expectation, if
$p < 2 r_{\bar x} r_{\bar y} r_{\bar z}/(4r_{\bar x} r_{\bar y} r_{\bar z} +r_{\bar x} r_{\bar y} + r_{\bar y} r_{\bar z}- 2 r_{\bar x} r_{\bar z}+ 2 r_{\bar y} -r_{\bar x} - r_{\bar z}+4)$, and inequalities obtained by switching $(x,i,n)$ with $(z,k,l)$ in~\eqref{nc4} hold in expectation, if
$p < 2r_{\bar x} r_{\bar y} r_{\bar z}/(4r_{\bar x} r_{\bar y} r_{\bar z} +r_{\bar x} r_{\bar z} + r_{\bar y} r_{\bar z}- 2 r_{\bar x} r_{\bar y}+ 2 r_{\bar z} -r_{\bar x} - r_{\bar y}+4)$. Since by assumption $r_{\bar x} \geq r_{\bar y} \geq r_{\bar z}$, these two inequalities are implied by inequality~\eqref{cond2}.

Define $I_1 = \{i \in [n]: \bar x_i = 1\}$.
We now show that event $A^2_1$ occurs with high probability:
\begin{align*}
&\prob[A^2_1]
=\prob\Bigg[\bigcap_{i \in I_1} {\Bigg\{\frac{1}{ml}\sum_{j,k}{\Big(\frac{1}{3} t^{x,i}_{jk \rightarrow 11}-f^{x,i}_{jk \rightarrow 11}\Big)}
-\frac{1}{m}\sum_{j}{\bar \nu^{y,j}_{x,i}}
-\frac{1}{l}\sum_{k}{\bar \nu^{z,k}_{x,i}} > 0}\Bigg\}\Bigg]   \\
& \geq \prob\Bigg[\bigcap_{i \in I_1} \Bigg\{\frac{1}{ml}\sum_{j,k}{\Big(\frac{1}{3} t^{x,i}_{jk \rightarrow 11}-f^{x,i}_{jk \rightarrow 11}\Big)}
-\avg\Big[\frac{1}{ml}\sum_{j,k}{\Big(\frac{1}{3} t^{x,i}_{jk \rightarrow 11}-f^{x,i}_{jk \rightarrow 11}\Big)}
\Big]\\
& \quad -\frac{1}{m}\sum_{j}{\bar \nu^{y,j}_{x,i}}
+\avg\Big[\frac{1}{m}\sum_{j}{\bar \nu^{y,j}_{x,i}}\Big]
-\frac{1}{l}\sum_{k}{\bar \nu^{z,k}_{x,i}}
+\avg\Big[\frac{1}{l}\sum_{k}{\bar \nu^{z,k}_{x,i}}\Big]> -\tilde\epsilon
\Bigg\}\Bigg]\\
& \geq \prob\Bigg[\bigcap_{i \in I_1} \Bigg\{\frac{1}{ml}\sum_{j,k}{\Big(\frac{1}{3} t^{x,i}_{jk \rightarrow 11}-f^{x,i}_{jk \rightarrow 11}\Big)} -\avg\Big[\frac{1}{ml}\sum_{j,k}{\Big(\frac{1}{3} t^{x,i}_{jk \rightarrow 11}-f^{x,i}_{jk \rightarrow 11}\Big)}
\Big]> -\frac{\tilde\epsilon}{3} \Bigg\} \\
& \quad \cap
\bigcap_{i \in I_1} \Bigg\{-\frac{1}{m}\sum_{j}{\bar \nu^{y,j}_{x,i}}
+\avg\Big[\frac{1}{m}\sum_{j}{\bar \nu^{y,j}_{x,i}}\Big] >  -\frac{\tilde\epsilon}{3}  \Bigg\} \cap
\bigcap_{i \in I_1} \Bigg\{-\frac{1}{l}\sum_{k}{\bar \nu^{z,k}_{x,i}}
+\avg\Big[\frac{1}{l}\sum_{k}{\bar \nu^{z,k}_{x,i}}\Big]>
-\frac{\tilde\epsilon}{3}
\Bigg\}\Bigg]\\
& \geq \prob\Bigg[\bigcap_{i \in I_1} \Bigg\{\Bigg|\frac{1}{ml}\sum_{j,k}{\Big(\frac{1}{3} t^{x,i}_{jk \rightarrow 11}-f^{x,i}_{jk \rightarrow 11}\Big)} -\avg\Big[\frac{1}{ml}\sum_{j,k}{\Big(\frac{1}{3} t^{x,i}_{jk \rightarrow 11}-f^{x,i}_{jk \rightarrow 11}\Big)}
\Big]\Bigg| < \frac{\tilde\epsilon}{3} \Bigg\} \\
& \quad  \cap
\bigcap_{i \in I_1} \Bigg\{\Bigg|\frac{1}{m}\sum_{j}{\bar \nu^{y,j}_{x,i}}
-\avg\Big[\frac{1}{m}\sum_{j}{\bar \nu^{y,j}_{x,i}}\Big]\Bigg| <  \frac{\tilde\epsilon}{3}  \Bigg\} \cap
\bigcap_{i \in I_1} \Bigg\{\Bigg|\frac{1}{l}\sum_{k}{\bar \nu^{z,k}_{x,i}}
-\avg\Big[\frac{1}{l}\sum_{k}{\bar \nu^{z,k}_{x,i}}\Big]\Bigg| <\frac{\tilde\epsilon}{3}
\Bigg\}\Bigg]\\
& \geq 1- \sum_{i \in I_1} \prob\Bigg[\Bigg|\frac{1}{ml}\sum_{j,k}{\Big(\frac{1}{3} t^{x,i}_{jk \rightarrow 11}-f^{x,i}_{jk \rightarrow 11}\Big)} -\avg\Big[\frac{1}{ml}\sum_{j,k}{\Big(\frac{1}{3} t^{x,i}_{jk \rightarrow 11}-f^{x,i}_{jk \rightarrow 11}\Big)}
\Big]\Bigg| \geq \frac{\tilde\epsilon}{3} \Bigg] \\
& \quad - \sum_{i \in I_1} \prob\Bigg[
\Bigg|\frac{1}{m}\sum_{j}{\bar \nu^{y,j}_{x,i}}
-\avg\Big[\frac{1}{m}\sum_{j}{\bar \nu^{y,j}_{x,i}}\Big]\Bigg| \geq  \frac{\tilde\epsilon}{3}  \Bigg] - \sum_{i \in I_1} \prob\Bigg[
\Bigg|\frac{1}{l}\sum_{k}{\bar \nu^{z,k}_{x,i}}
-\avg\Big[\frac{1}{l}\sum_{k}{\bar \nu^{z,k}_{x,i}}\Big]\Bigg| \geq \frac{\tilde\epsilon}{3}
\Bigg]\\
& \geq 1 - 2n\exp\Big(-\frac{m l \tilde\epsilon^2}{8}\Big) - 2n\exp\Big(-\frac{2m\tilde\epsilon^2}{9}\Big) -2n\exp\Big(-\frac{2l\tilde\epsilon^2}{9}\Big),
\end{align*}
where the first inequality follows since $\epsilon \geq \tilde \epsilon > 0$, the second and third inequalities follow from set inclusion.
The fourth inequality follows from taking the union bound.
The last inequality follows from the application of the Hoeffding inequality by noting that (i) random variables
$s^i_{jk} := \frac{1}{3} t^{x,i}_{jk \rightarrow 11}-f^{x,i}_{jk \rightarrow 11}$ for all $i \in [n],j \in [m],k \in [l]$ are independent and $-1 \leq s^i_{jk} \leq \frac{1}{3}$, (ii) random variables $\bar \nu^{y,j}_{x,i}$
for all $j \in [m]$ are independent with
$0 \leq \bar \nu^{y,j}_{x,i} \leq 1$ and, (iii) random variables $\bar \nu^{z,k}_{x,i}$
for all $k \in [l]$ are independent with
$0 \leq \bar \nu^{z,k}_{x,i} \leq 1$.
The proof then follows from the fact that $\tilde \epsilon$ is a positive constant and because the limit assumptions in the theorem imply that, as $n, m, l \to \infty$, we have $n \exp (-m)$, $n \exp (-l)$, $n \exp (-ml)$ go to zero.
\end{cpf}
\end{prfc}



%

\section{Recovery proof for the flower LP}
\label{sec FR recovery}

The main goal of this section is to prove \cref{th LP2}.
To this end, we first obtain a deterministic sufficient condition for recovery.
Subsequently, we study the semi-random corruption model.
Since the feasible region of the flower LP is a subset of the feasible region of the standard LP, if all assumptions of Proposition~\ref{unique1} are satisfied, then then
the ground truth is the unique optimal solution of flower LP. 
Similarly, if all assumptions of Theorem~\ref{th LP1} are satisfied, then the flower LP recovers the ground truth with high probability.

\subsection{Deterministic recovery guarantee}
To obtain a deterministic condition for recovery, we first present a sufficient condition under which an optimal solution of Problem~\eqref{LP3} coincides with the ground truth. Next, we investigate the question of uniqueness.
For notational simplicity, for any $r\in \{0,1\}$, we define
\begin{align*}
&T^{x,y,i,j}_r = \Big|k \in [l]: \bar z_k = r, \; (i,j,k) \in \T\Big|.
\end{align*}
Parameters $T^{x,z,i,k}_r$, $T^{y,z,j,k}_r$ are similarly defined. Moreover for each $r,s,t \in \{0,1\}$ we define
\begin{align*}
    T_{rst}^{x,i} = \Big|(j,k,i') \in [m] \times [l] \times [n]:  \bar y_j = r, \bar z_k =s, \bar x_{i'} =t,   \; (i,j,k) \in \T, (i',j,k) \in \T\Big|.
\end{align*}
Parameters $T_{rst}^{y,j}$ and $T_{rst}^{z,k}$ are similarly defined. Finally, define
$$
\bar T^{x,y}_1 = \frac{1}{n m} \sum_{i,j}{T^{x,y,i,j}_1}.
$$
Parameters $\bar T^{x,z}_1$ and $\bar T^{y,z}_1$ are similarly defined. Since the feasible region of the flower LP is a subset of the feasible region of the standard LP, 
by Proposition~\ref{unique1 nonoise},
if $\F = \emptyset$, Problem~\eqref{LP3} recovers the ground truth provided that $\bar x, \bar y, \bar z \neq 0$. Therefore, in the following, we consider the case with $\F \neq \emptyset$.

\begin{proposition}\label{deter2}
Let $\bar x \in \{0,1\}^n$, $\bar y \in \{0,1\}^m$, $\bar z \in \{0,1\}^l$ and define $\bar \W = (\bar w_{ijk}) := \bar x \otimes \bar y \otimes \bar z  \in \{0,1\}^{n \times m \times l}$.
Suppose that $\F \neq \emptyset$.
Then $(\bar x, \bar y, \bar z, \bar \W)$ is an optimal solution of Problem~\eqref{LP3}, if in addition to Condition~\ref{c1-1} of Proposition~\ref{deter1},
the following conditions are satisfied:
\begin{enumerate}
    \item \label{c2-1}
    For each $i \in [n],j \in [m]$ with $\bar x_i = \bar y_j= 1$, we
    have  $T^{x,y,i,j}_1 \geq 1$,
    for each $j \in [m],k \in [l]$ with $\bar y_j = \bar z_k = 1$, we have $T^{y,z,j,k}_1 \geq 1$ and, for each
    $i \in [n],k \in [l]$ with $\bar x_i = \bar z_k = 1$, we have
    $T^{x,z,i,k}_1 \geq 1$.


   \item \label{c2-3}
Let $\alpha < 1$ be a constant arbitrarily close to 1. Then for each $i \in [n]$ with $\bar x_i = 0$, we have
      $T^{x,i}_{111} \geq \alpha T^{x,i}_{11} \bar T^{y,z}_1$,
      for each $j \in [m]$ with $\bar y_j = 0$, we have
      $T^{y,j}_{111} \geq \alpha T^{y,j}_{11} \bar T^{x,z}_1$, and
      for each $k \in [l]$ with $\bar z_k = 0$, we have
      $T^{z,k}_{111} \geq \alpha T^{z,k}_{11} \bar T^{x,y}_1$.

    \item  \label{c2-2}
    For each $(i,j,k) \in \T$ with $\bar x_i = 0$, $\bar y_j = \bar z_k = 1$, we have
       \begin{equation} \label{aa1}
   \min\Big\{T^{x,i}_{11}, \frac{T^{x,i}_{111}}{T^{y,z,j,k}_1}\Big\} \geq
     \frac{1}{3}F^{x,i}_{00}+\frac{1}{2}(F^{x,i}_{01}+F^{x,i}_{10})+F^{x,i}_{11}.
    \end{equation}
    Moreover,~\eqref{aa1} holds when $(x,i,n)$ is replaced
    with $(y,j,m)$ and $(z,k,l)$, respectively.

\item \label{c2-4}
For each $(i,j,k) \in \T \cap \P$, define
\begin{align*}
    \bar \gamma_{ijk} & :=\frac{1}{3}\Bigg(2 - \frac{1}{3}\frac{n_{\bar x}}{T^{y,z,j,k}_1}
    -\frac{1}{3}\frac{n_{\bar y}}{T^{x,z,i,k}_1}-\frac{1}{3}\frac{n_{\bar z}}{T^{x,y,i,j}_1}\nonumber\\
    & \quad -\frac{1}{\alpha \bar T^{y,z}_1}\sum_{\substack{i': \bar x_{i'} = 0\\ (i',j,k) \in \T}}{\frac{1}{T_{11}^{x,i'}}\Big(\frac{1}{3}F^{x,i'}_{00}+\frac{1}{2}(F^{x,i'}_{01}+F^{x,i'}_{10})+F^{x,i'}_{11}\Big)} \nonumber\\
    & \quad -\frac{1}{\alpha \bar T^{x,z}_1}\sum_{\substack{j': \bar y_{j'} =0\\ (i,j',k) \in \T}} { \frac{1}{T_{11}^{y,j'}}\Big(\frac{1}{3}F^{y,j'}_{00}+\frac{1}{2}(F^{y,j'}_{01}+F^{y,j'}_{10})+F^{y,j'}_{11}\Big)}\nonumber\\
    & \quad -\frac{1}{\alpha \bar T^{x,y}_1}\sum_{\substack{k': \bar z_{k'} =0 \\ (i,j,k') \in \T}} {\frac{1}{T_{11}^{z,k'}}\Big(\frac{1}{3}F^{z,k'}_{00}+\frac{1}{2}(F^{z,k'}_{01}+F^{z,k'}_{10})+F^{z,k'}_{11}\Big)}\Bigg).
\end{align*}
Then for each $(i,j,k) \in \T \cap \P$, we have
\begin{equation}\label{newass}
3\bar \gamma_{ijk} \geq \max\Bigg\{\frac{F^{x,i}_{11}}{T^{x,i}_{11}}, \;
\frac{F^{y,j}_{11}}{T^{y,j}_{11}}, \; \frac{F^{z,k}_{11}}{T^{z,k}_{11}}\Bigg\}.
\end{equation}
\end{enumerate}

\end{proposition}
\begin{prf}
We start by constructing the dual of Problem~\eqref{LP3}.
Define dual variables $\lambda^x_{ijk}, \lambda^y_{ijk}, \lambda^z_{ijk}$ for all $(i,j,k) \in S_1$ associated with the first, the second and the third set of constraints in~\eqref{f1}, respectively.  Define $\mu^1_{ijk}, \mu^2_{ijk}$ for all $(i,j,k) \in S_0$ associated with first and second set of constraints in~\eqref{f2}, respectively.
Define $f^x_{ijki'}$ for all $(i',j,k) \in S_1$ and for all $(i,j,k) \in S_0$
associated with constraints~\eqref{f3}, $f^y_{ijkj'}$ for all $(i,j',k) \in S_1$ and for all $(i,j,k) \in S_0$ associated with constraints~\eqref{f4}, and $f^z_{ijkk'}$ for all $(i,j,k') \in S_1$ and for all $(i,j,k) \in S_0$
associated with constraints~\eqref{f5}.
Finally, define $u^x_i$ (resp. $l^x_i$) for all $i \in [n]$, $u^y_j$
(resp. $l^y_j$) for all $j \in [m]$, and $u^z_k$ (resp. $l^z_k$) for all $k \in [l]$, associated with $x_i \leq 1$ (resp $-x_i \leq 0$), $y_j \leq 1$ (resp $-y_j \leq 0$), $z_k \leq 1$ (resp $-z_k \leq 0$), respectively.
For notational simplicity, let
$F^x = \{(i,j,k,i'): (i',j,k) \in S_1, (i,j,k) \in S_0\}$,
$F^y = \{(i,j,k,j'): (i,j',k) \in S_1, (i,j,k) \in S_0\}$, and
$F^z = \{(i,j,k,k'): (i,j,k') \in S_1, (i,j,k) \in S_0\}$.
The dual of Problem~\eqref{LP3} is given by:
\begin{align}\label{Fdual}
\tag{fD}
\max \  & |S_1|-2 \sum_{\substack{(i,j,k)\\ \in S_0}}{\mu^2_{ijk}}
-\sum_{\substack{(i,j,k,i') \\\in F^x}}{f^x_{ijki'}}
-\sum_{\substack{(i,j,k,j') \\ \in F^y}}{f^y_{ijkj'}}
-\sum_{\substack{(i,j,k,k') \\ \in F^z}}{f^x_{ijkk'}}
-\sum_{i}{u^x_i}-\sum_{j}{u^y_j}-\sum_{k}{u^z_k}\nonumber \\
\st \  &  \sum_{(j,k):(i,j,k) \in S_0}{\mu^2_{ijk}}+ \sum_{\substack{(j,k,i'): \\ (i,j,k,i') \in F^x}}{f^x_{ijki'}} + u^x_i - l^x_i= \sum_{(j,k):(i,j,k) \in S_1}{\lambda^x_{ijk}} , \qquad \forall i \in [n]\label{g1}\\
&  \sum_{(i,k):(i,j,k) \in S_0}{\mu^2_{ijk}}+\sum_{\substack{(i,k,j'): \\ (i,j, k,j') \in F^y}}{f^y_{ijkj'}} + u^y_j - l^y_j=  \sum_{(i,k):(i,j,k) \in S_1}{\lambda^y_{ijk}}, \qquad \forall j \in [m]\label{g2}\\
& \sum_{(i,j):(i,j,k) \in S_0}{\mu^2_{ijk}}+\sum_{\substack{(i,j,k'): \\ (i,j,k,k') \in F^z}}{f^z_{ijkk'}} + u^z_k -l^z_k=\sum_{(i,j):(i,j,k) \in S_1}{\lambda^z_{ijk}} , \qquad \forall k \in [l]\label{g3}\\
& \lambda^x_{ijk} + \lambda^y_{ijk} + \lambda^z_{ijk}
+\sum_{\substack{i':(i',j,k,i) \\ \in F^x}}{f^x_{i'jki}}
+\sum_{\substack{j':(i,j',k,j) \\ \in F^y}}{f^y_{ij'kj}}
+\sum_{\substack{k':(i,j,k',k)\\  \in F^z}}{f^z_{ijk'k}} = 1, \; \forall (i,j,k) \in S_1 \label{g4}\\
& \mu^1_{ijk}+ \mu^2_{ijk} +\sum_{\substack{i':(i,j,k,i') \\ \in F^x}}{f^x_{ijki'}}
+\sum_{\substack{j':(i,j,k,j') \\ \in F^y}}{f^y_{ijkj'}}
+\sum_{\substack{k':(i,j,k,k') \\ \in F^z}}{f^z_{ijkk'}}
= 1 , \; \forall (i,j,k) \in S_0\label{g4m}\\
& \lambda^x_{ijk}\geq 0, \; \lambda^y_{ijk}\geq 0, \; \lambda^z_{ijk}\geq 0, \qquad \forall (i,j,k) \in S_1\label{g5}\\
& \mu^1_{ijk}\geq 0, \;  \mu^2_{ijk}\geq 0 \qquad \forall (i,j,k) \in S_0 \label{g6}\\
& f^x_{ijki'} \geq 0, \; \forall (i,j,k,i') \in F^x, \
  f^y_{ijkj'} \geq 0, \; \forall (i,j,k,j') \in F^y, \
  f^z_{ijkk'} \geq 0, \; \forall (i,j,k,k') \in F^z \label{g7}\\
 & l^x_i \geq 0, \; \forall i \in [n], \; l^y_j \geq 0, \; \forall j \in [m],\; l^z_k \geq 0, \; \forall k \in [l] \label{usl2}\\
& u^x_i \geq 0, \; \forall i \in [n], \; u^y_j \geq 0, \; \forall j \in [m],\; u^z_k \geq 0, \; \forall k \in [l] \label{g8}.
\end{align}
To prove the optimality of $(\bar x, \bar y, \bar z, \bar \W)$, it suffices to construct a dual feasible point
$(\bar\lambda^x, \bar\lambda^y, \bar\lambda^z, \bar\mu^1,$ $\bar \mu^2, \bar f^x,\bar f^y, \bar f^z, \bar l^x, \bar l^y, \bar l^z, \bar u^x, \bar u^y, \bar u^z)$ for Problem~\eqref{Fdual} that satisfies complementary slackness.
First,  we set
$\bar l^x_i = 0$ for all $i \in [n]$, $\bar l^y_j = 0$ for all $j \in [m]$, $\bar l^z_k = 0$ for all $k \in [l]$ and $\bar \mu^2_{ijk} = 0$ for all $(i,j,k) \in \S_0$.
By complementary slackness, we have:
\begin{enumerate}
\item [(I)] For each $(i,j,k) \in \F \cap \N$ with (i) $\bar x_i =1$, we have $\bar \lambda^x_{ijk} = 0$, (ii) $\bar y_j =1$, we have $\bar \lambda^y_{ijk} = 0$, and (iii) $\bar z_k =1$, we have $\bar \lambda^z_{ijk} = 0$.
\item [(II)] For each $(i,j,k) \in \F \cap \P$, we have $\bar \mu^1_{ijk} = 0$.
\item  [(III)] For each $(i',j,k) \in \F \cap \N$
and $(i,j,k) \in \F \cap \P$, we have $\bar f^x_{ijki'} = 0$;
for each $(i,j',k) \in \F \cap \N$
and $(i,j,k) \in \F \cap \P$, we have $\bar f^y_{ijkj'} = 0$;
for each $(i,j,k') \in \F \cap \N$
and $(i,j,k) \in \F \cap \P$, we have $\bar f^z_{ijkk'} = 0$.
\item [(IV)] For each $i \in [n]$, with $\bar x_i = 0$, we have $\bar u^x_i = 0$; for each $j \in [m]$, with $\bar y_j = 0$, we have $\bar u^y_j = 0$;
for each $k \in [l]$, with $\bar z_k = 0$, we have $\bar u^z_k = 0$.
\end{enumerate}
\paragraph{Simplifications.} To construct the dual certificate, we make the following simplifications:
\begin{itemize}
\item for each $(i,j,k) \in \T \cap \P$, we let
\begin{equation}\label{n1}
\bar \lambda^x_{ijk} =  \bar \lambda^y_{ijk} = \bar \lambda^z_{ijk} = \gamma_{ijk}.
\end{equation}
We establish the non-negativity of $\lambda^x_{ijk},  \bar \lambda^y_{ijk}, \bar \lambda^z_{ijk}$ later when we determine $\gamma_{ijk}$.
\item for each $(i,j,k) \in \F \cap \N$, we set:
\begin{align}
\label{el}
\begin{split}
&\bar \lambda^x_{ijk} =  \bar \lambda^y_{ijk} = \bar \lambda^z_{ijk} = \frac{1}{3}, \qquad {\rm if}\; \bar x_i = \bar y_j = \bar z_k = 0 \\
&\bar \lambda^y_{ijk}  = \bar \lambda^z_{ijk} = \frac{1}{2}, \qquad \text{if } \bar x_i = 1, \; \bar y_j =\bar z_k = 0 \\
&\bar \lambda^x_{ijk}  = \bar \lambda^z_{ijk} = \frac{1}{2}, \qquad \text{if } \bar x_i = 0, \; \bar y_j =1, \; \bar z_k = 0\\
&\bar \lambda^x_{ijk}  = \bar \lambda^y_{ijk} = \frac{1}{2}, \qquad \text{if } \bar x_i = \bar y_j =0, \; \bar z_k = 1.
\end{split}
\end{align}
\item for each $(i,j,k) \in \T \cap \N$, we set:
\begin{align}\label{smplfy}
&\bar f^x_{ijki'} = 0, \qquad \forall (i',j,k) \in \F \cap \N\nonumber\\
&\bar f^y_{ijkj'} = 0, \qquad \forall (i,j',k) \in \F \cap \N\\
&\bar f^z_{ijkk'} = 0, \qquad \forall (i,j,k') \in \F \cap \N\nonumber
\end{align}
and
\begin{align}\label{smplfy2}
&\bar f^x_{ijki'} = \alpha_i^x, \qquad \forall (i',j,k) \in \T \cap \P \nonumber\\
&\bar f^y_{ijkj'} = \alpha_j^y, \qquad \forall (i,j',k) \in \T \cap \P\\
&\bar f^z_{ijkk'} = \alpha_k^z, \qquad \forall (i,j,k') \in \T \cap \P,\nonumber
\end{align}
where $\alpha_i^x,\alpha_j^y, \alpha_k^z$ are to be determined later.
\item for each $(i,j,k) \in \F \cap \P$, we set:
\begin{align}\label{hv}
&\bar f^x_{ijki'} = \frac{1}{3T^{y,z,j,k}_1}, \qquad \forall (i',j,k) \in \T \cap \P\nonumber\\
&\bar f^y_{ijkj'} = \frac{1}{3T^{x,z,i,k}_1}, \qquad \forall (i,j',k) \in \T \cap \P\\
&\bar f^z_{ijkk'} = \frac{1}{3T^{x,y,i,j}_1}, \qquad \forall (i,j,k') \in \T \cap \P \nonumber,
\end{align}
where by Condition~\ref{c2-1}, we have $T^{y,z,j,k}_1 \geq 1$,
$T^{x,z,i,k}_1 \geq 1$, and $T^{x,y,i,j}_1 \geq 1$.
\end{itemize}
%
%
Using these simplifications, in the following we establish dual feasibility. For clarity of presentation, we consider different type of constraints of Problem~\eqref{Fdual}, separately.

\paragraph{Constraints~\eqref{g1}--\eqref{g3}:} By complementary slackness, equations~\eqref{el} and~\eqref{smplfy2},
it follows that for each $i \in [n]$ with $\bar x_i = 0$,
constraints~\eqref{g1} simplify to

\begin{equation*}
\sum_{\substack{(j,k,i'): (i,j,k) \in \T,\\ (i',j,k) \in \T: \; \bar x_{i'}=\bar y_j = \bar z_k = 1}}{ \alpha_i^x} = \sum_{\substack{(j,k):(i,j,k) \in \F\\ \bar y_j =\bar z_k = 0}}{\frac{1}{3}}+
\sum_{\substack{(j,k):(i,j,k) \in \F:\\ (\bar y_j =1, \bar z_k = 0) \lor
(\bar y_j =0, \bar z_k = 1)}} {\frac{1}{2}}+
\sum_{\substack{(j,k):(i,j,k) \in \F: \\ \bar y_j = \bar z_k  = 1}} {1}.
\end{equation*}
By Condition~\ref{c1-1} of Proposition~\ref{deter1}, $F^{x,i} > 0$; \ie the right-hand of the above equality is positive. By inequality~\eqref{aa1} of Condition~\ref{c2-2}, $T_{111}^{x,i} > 0$; \ie the left-hand of the above equality is positive.
Hence:
\begin{equation}\label{alfa2x}
\alpha^x_i  = \frac{1}{T_{111}^{x,i}}\Big(\frac{1}{3}F^{x,i}_{00}+\frac{1}{2}(F^{x,i}_{01}+F^{x,i}_{10})+F^{x,i}_{11}\Big).
\end{equation}
%
%
Clearly $\alpha^x_i \geq 0$, satisfying constraints~\eqref{g5}. Substituting~\eqref{n1} and~\eqref{hv}
in constraints~\eqref{g1} and using complementary slackness, for each $i \in [n]$ with $\bar x_i = 1$ we obtain
\begin{align}\label{g11p}
\bar u^x_i & = \sum_{\substack{(j,k):(i,j,k) \in \T:\\ \bar y_j = \bar z_k=1}}{\gamma_{ijk}}-
\frac{1}{3}\sum_{\substack{(j,k,i'): (i',j,k) \in \T, \\ (i,j,k) \in \F:\\ \bar x_{i'}=\bar y_j = \bar z_k = 1}} {\frac{1}{T^{y,z,j,k}_1}}\nonumber\\
& =\sum_{\substack{(j,k):(i,j,k) \in \T,\\ \bar y_j = \bar z_k=1}}{\gamma_{ijk}}-
\frac{1}{3}\sum_{\substack{(j,k): (i,j,k) \in \F,\\ \bar y_j = \bar z_k = 1}} \Bigg(\sum_{\substack{i': (i',j,k) \in \T,\\ \bar x_{i'}=1}}{\frac{1}{T^{y,z,j,k}_1}}\Bigg) \nonumber \\
& = \sum_{\substack{(j,k):(i,j,k) \in \T,\\ \bar y_j = \bar z_k=1}}{\gamma_{ijk}}-
\frac{1}{3} F^{x,i}_{11}.
\end{align}
We establish non-negativity of $\bar u^x_i$ after we determine $\gamma_{ijk}$.
%
%
Similarly, substituting~\eqref{el} and~\eqref{smplfy2} in constraints~\eqref{g2} and using complementary slackness,
for each $j \in [m]$ with $\bar y_j = 0$ we get:
\begin{equation*}
\label{g20p}
\sum_{\substack{(i,k,j'): (i,j,k) \in \T,\\ (i,j',k) \in \T: \bar x_i = \bar y_{j'} = \bar z_k =1}}{\alpha_j^y} = \sum_{\substack{(i,k):(i,j,k) \in \F \\ \bar x_i =\bar z_k = 0}}{\frac{1}{3}}+
\sum_{\substack{(i,k):(i,j,k) \in \F:\\ (\bar x_i =1, \bar z_k = 0) \lor
(\bar x_i =0, \bar z_k = 1)}} {\frac{1}{2}}+
\sum_{\substack{(i,k):(i,j,k) \in \F: \\ \bar x_i = \bar z_k  = 1}} {1}.
\end{equation*}
By Condition~\ref{c1-1} of Proposition~\ref{deter1}, $F^{y,j} > 0$; \ie the right-hand of the above equality is positive.  By symmetric counterpart of inequality~\eqref{aa1} of Condition~\ref{c2-2}, $T_{111}^{y,j} > 0$;
\ie the left-hand of the above equality is positive. Hence:
\begin{equation}\label{alfa2y}
\alpha^y_i  = \frac{1}{T_{111}^{y,j}}\Big(\frac{1}{3}F^{y,j}_{00}+\frac{1}{2}(F^{y,j}_{01}+F^{y,j}_{10})+F^{y,j}_{11}\Big).
\end{equation}
Clearly $\alpha^y_j \geq 0$, satisfying constraints~\eqref{g5}.
Substituting~\eqref{n1} and~\eqref{hv} in constraints~\eqref{g2},
for each $j \in [m]$ with $\bar y_j = 1$, we get
\begin{equation*}
\label{g21p}
 \bar u^y_j = \sum_{\substack{(i,k):(i,j,k) \in \T:\\ \bar x_i= \bar z_k =1}}{\gamma_{ijk}}
 -\frac{1}{3}F^{y,j}_{11}.
\end{equation*}
We establish non-negativity of $\bar u^y_j$ after we determine $\gamma_{ijk}$.
Finally, substituting~\eqref{el} and~\eqref{smplfy2}  in constraints~\eqref{g3} for each $k \in [l]$ with $\bar z_k = 0$, we obtain
\begin{equation*}
\sum_{\substack{(i,j,k'): (i,j,k) \in \T,\\ (i,j,k') \in \T: \bar x_i = \bar y_j = \bar z_{k'}=1}}{\alpha_k^z} = \sum_{\substack{(i,j):(i,j,k) \in \F\\ \bar x_i =\bar y_j = 0}}{\frac{1}{3}}+
\sum_{\substack{(i,j):(i,j,k) \in \F:\\ (\bar x_i =1, \bar y_j = 0) \lor
(\bar x_i =0, \bar y_j = 1)}} {\frac{1}{2}}+
\sum_{\substack{(i,j):(i,j,k) \in \F: \\ \bar x_i = \bar y_j  = 1}} {1}. \end{equation*}
By Condition~\ref{c1-1} of Proposition~\ref{deter1}, $F^{z,k} > 0$; \ie the right-hand of the above equality is positive. By a symmetric counterpart of inequality~\eqref{aa1} of Condition~\ref{c2-2}, $T_{111}^{z,k} > 0$\ie the right-hand of the above equality is positive. Hence
\begin{equation}\label{alfa2z}
\alpha^z_k  = \frac{1}{T_{111}^{z,k}}\Big(\frac{1}{3}F^{z,k}_{00}+\frac{1}{2}(F^{z,k}_{01}+F^{z,k}_{10})+F^{z,k}_{11}\Big).
\end{equation}
Substituting~\eqref{n1} and~\eqref{hv} in constraints~\eqref{g3}, for each $k \in [l]$ with $\bar z_k = 1$, we obtain
\begin{equation*}
\bar u^z_k = \sum_{\substack{(i,j):(i,j,k) \in \T: \\ \bar x_i = \bar y_j=1}}{\gamma_{ijk}}-
\frac{1}{3} F^{z,k}_{11}.
\end{equation*}
We establish non-negativity of $\bar u^z_k$ after we determine $\gamma_{ijk}$.
%
%
\paragraph{Constraints~\eqref{g4}:} The following cases arise:
\begin{itemize}
    \item If $(i,j,k) \in \T \cap \P$, substituting~\eqref{n1},~\eqref{smplfy2} and~\eqref{hv} in constraints~\eqref{g4} we obtain:
\begin{align*}
    & 3\gamma_{ijk} +\sum_{i': (i',j,k) \in \T \cap \N}{\alpha_{i'}^x}
     +\sum_{i': (i',j,k) \in \F \cap \P}{\frac{1}{3T^{y,z,j,k}_1}}
     +\sum_{j': (i,j',k) \in \T \cap \N} {\alpha_{j'}^y} \\
     & \quad  +\sum_{j': (i,j',k) \in \F \cap \P} {\frac{1}{3T^{x,z,i,k}_1}}
     +\sum_{k': (i,j,k') \in \T \cap \N} {\alpha_{k'}^z}
     +\sum_{k': (i,j,k') \in \F \cap \P} {\frac{1}{3T^{x,y,i,j}_1}}=1.
    \end{align*}
    Substituting for $\alpha^x_i, \alpha^y_j, \alpha^z_k$ using~\eqref{alfa2x},~\eqref{alfa2y}, and~\eqref{alfa2z}, respectively, we obtain:
    \begin{align*}
    \gamma_{ijk} & =\frac{1}{3}\Bigg(2 - \frac{1}{3}\frac{n_{\bar x}}{T^{y,z,j,k}_1}
    -\frac{1}{3}\frac{n_{\bar y}}{T^{x,z,i,k}_1}-\frac{1}{3}\frac{n_{\bar z}}{T^{x,y,i,j}_1}\nonumber\\
    & \quad -\sum_{\substack{i': \bar x_{i'} = 0\\ (i',j,k) \in \T}}{\frac{1}{T_{111}^{x,i'}}\Big(\frac{1}{3}F^{x,i'}_{00}+\frac{1}{2}(F^{x,i'}_{01}+F^{x,i'}_{10})+F^{x,i'}_{11}\Big)} \nonumber\\
    & \quad -\sum_{\substack{j': \bar y_{j'} =0\\ (i,j',k) \in \T}} { \frac{1}{T_{111}^{y,j'}}\Big(\frac{1}{3}F^{y,j'}_{00}+\frac{1}{2}(F^{y,j'}_{01}+F^{y,j'}_{10})+F^{y,j'}_{11}\Big)}\nonumber\\
    & \quad -\sum_{\substack{k': \bar z_{k'} =0 \\ (i,j,k') \in \T}} {\frac{1}{T_{111}^{z,k'}}\Big(\frac{1}{3}F^{z,k'}_{00}+\frac{1}{2}(F^{z,k'}_{01}+F^{z,k'}_{10})+F^{z,k'}_{11}\Big)}\Bigg).
\end{align*}
    From Condition~\ref{c2-3} it follows that $\gamma_{ijk} \geq \bar \gamma_{ijk}$
    for all $(i,j,k) \in \T \cap \P$. Hence the non-negativity of $\gamma_{ijk}$, \ie non-negativity of $\bar \lambda^x_{ijk}, \bar \lambda^y_{ijk},\bar\lambda^z_{ijk}$ for each $(i,j,k) \in \T \cap \P$, follows from Condition~\ref{c2-4}.
    \item If $(i,j,k) \in \F \cap \N$ such that $\bar x_i = \bar y_j = \bar z_k =0$, substituting~\eqref{el} and~\eqref{smplfy} in constraints~\eqref{g4} we obtain: $\frac{1}{3} + \frac{1}{3} + \frac{1}{3} = 1$; \ie these constraints are satisfied.
    \item If $(i,j,k) \in \F \cap \N$ such that $\bar x_i = 1$, $\bar y_j = \bar z_k =0$, substituting~\eqref{el} and~\eqref{smplfy} in constraints~\eqref{g4}, and using
    condition~(I) of complementary slackness (\ie $\bar \lambda^x_{ijk} = 0$) we obtain:
    $\frac{1}{2} + \frac{1}{2} = 1$; \ie these constraints are satisfied. By symmetry, constraints~\eqref{g4} are satisfied if $(i,j,k) \in \F \cap \N$ and
    (i) $\bar x_i = 0$, $\bar y_j = 1$,
    $\bar z_k =0$, or (ii) $\bar x_i = \bar y_j =0$, $\bar z_k =1$.
    \item If $(i,j,k) \in \F \cap \N$ such that $\bar x_i = \bar y_j =1$, $\bar z_k =0$, by condition~(I) of complementary slackness we have $\bar \lambda^x_{ijk} =\bar \lambda^y_{ijk} = 0$. Hence, by~\eqref{smplfy}, constraints~\eqref{g4} simplify to $\bar \lambda^z_{ijk} =1$. By symmetry,
    if $(i,j,k) \in \F \cap \N$ such that $\bar x_i = 1$, $\bar y_j =0$, $\bar z_k =1$, constraints~\eqref{g4} simplify to
    $\bar \lambda^y_{ijk} =1$, and
    if $(i,j,k) \in \F \cap \N$ such that $\bar x_i = 0$, $\bar y_j =1, \bar z_k =1$, constraints~\eqref{g4} simplify to
    $\bar \lambda^x_{ijk} =1$.

\end{itemize}
\paragraph{Constraints~\eqref{g4m}:} Two cases arise:
\begin{itemize}
    \item If $(i,j,k) \in \F \cap \P$, then by condition~(II) of complementary slackness we have $\bar \mu^1_{ijk} = 0$; hence by~\eqref{hv} we obtain:
    \begin{align*}
    \sum_{\substack{i': \bar x_{i'}=1, \\(i',j,k) \in \T}} {\frac{1}{3T^{y,z,j,k}_1}}+
    \sum_{\substack{j': \bar y_{j'}=1, \\(i,j',k) \in \T }} {\frac{1}{3T^{x,z,i,k}_1}}+
    \sum_{\substack{k':\bar z_{k'}=1,\\ (i,j,k') \in \T}}
    {\frac{1}{3T^{x,y,i,j}_1}} =
    \frac{1}{3} + \frac{1}{3} + \frac{1}{3}=
    1.
    \end{align*}
 That is, in this case, constraints~~\eqref{g4m} are satisfied.
 \item If $(i,j,k) \in \T \cap \N$, by projecting out variables $\mu^1_{ijk}$ and using~\eqref{smplfy2} we obtain:
    $$
\sum_{i':(i',j,k) \in \T \cap \P}{\alpha^x_i}+
    \sum_{j':(i,j',k) \in \T \cap \P}{\alpha^y_j}+
    \sum_{k':(i,j,k') \in \T \cap \P}{\alpha^z_k} \leq 1.
    $$
    This inequality in turn corresponds to the following cases (in all remaining cases, it simplifies to $0 \leq 1$):
    \begin{itemize}
    \item If $\bar x_i = 0$, $\bar y_j = \bar z_k =1$, we get
    \begin{equation*}
    \sum_{\substack{i': \bar x_{i'}=1, \\ (i',j,k) \in \T}}{\alpha^x_{i}} \leq 1.
    \end{equation*}
    By~\eqref{alfa2x}, the above inequality is implied by inequality~\eqref{aa1}.

    \item If $\bar x_i = 1$, $\bar y_j = 0$, $\bar z_k =1$, we get
    $$
    \sum_{\substack{j': \bar y_{j'}=1, \\(i,j',k) \in \T}}{\alpha^y_j} \leq 1.
    $$
    By~\eqref{alfa2y}, the above inequality
    is implied by a symmetric counterpart of inequality~\eqref{aa1}.
\item If $\bar x_i = \bar y_j =1$, $\bar z_k = 0$, we get
        \begin{equation*}
    \sum_{\substack{k': \bar z_{k'}=1, \\ (i,j,k') \in \T}}{\alpha^z_k} \leq 1.
    \end{equation*}
    By~\eqref{alfa2z}, the above inequality is implied by a symmetric counterpart of inequality~\eqref{aa1}.
     \end{itemize}
\end{itemize}
\paragraph{Final step.} It remains to establish non-negativity of $\bar u^x_i,\bar u^y_j,\bar u^z_k$. Recall that for each $i \in [n]$, $\bar u^x_i$ is given by~\eqref{g11p} if $\bar x_i = 1$ and equals zero, otherwise. Since $\gamma_{ijk} \geq \bar \gamma_{ijk}$, it follows that if
\begin{equation}\label{mmn}
3\sum_{\substack{(j,k): (i,j,k) \in \T \cap \P}}{\bar \gamma_{ijk}} \geq F^{x,i}_{11},
\end{equation}
then we have  $\bar u^x_i \geq 0$. By Condition~\ref{c2-4}, for each $(i,j,k) \in \T \cap \P$, we have
$3\bar \gamma_{ijk} \geq F^{x,i}_{11}/T^{x,i}_{11}$. Summing both sides of this inequality over all $(j,k)$ for which $(i,j,k) \in \T \cap \P$, we obtain inequality~\eqref{mmn}. The non-negativity of $\bar u^y_j,\bar u^z_k$ follows from a similar line of arguments.
\end{prf}

Next, utilizing Proposition~\ref{mang}, we present our uniqueness condition:

\begin{proposition}\label{unique2}
Suppose that all assumptions of Proposition~\ref{deter2} hold; moreover, suppose that inequalities~\eqref{aa1} (and symmetric counterparts) and inequalities~\eqref{newass} are strictly satisfied.
 Then $(\bar x, \bar y, \bar z, \bar \W)$ is the unique optimal solution of Problem~\eqref{LP3}.
\end{proposition}
\begin{prf}
Let $(\bar x, \bar y, \bar z, \bar \W)$ be an optimal solution of Problem~\eqref{LP1}.
To prove the statement it suffices to show there is no nonzero
$(x, y, z, \W)$ satisfying condition~\eqref{uc}.
\begin{itemize}
    \item [(i)] Since inequality~\eqref{newass} is strictly satisfied, we have $\bar u^x_i > 0$ for all $i \in [n]$ with $\bar x_i = 1$. This in turn implies that $x_i = 0$ for all $i \in [n]$ with $\bar x_i = 1$.
    By symmetry, we obtain $y_j = 0$ for all $j \in [m]$ with $\bar y_j = 1$ and $z_k = 0$ for all $k \in [l]$ with $\bar z_k = 1$.
    \item [(ii)] Since inequality~\eqref{newass} is strictly satisfied, we have
    $\gamma_{ijk} = \lambda^x_{ijk} > 0$ for all $(i,j,k) \in \T \cap \P$. Hence, we have $w_{ijk} = x_i = 0$ for all $(i,j,k) \in \T \cap \P$, where the second equality follows from part~(i) above since $\bar x_i = 1$.
    \item [(iii)] Since inequalities~\eqref{aa1} are strictly satisfied, we have $\bar \mu^1_{ijk} > 0$ for all $(i,j,k) \in \T \cap \N$, implying $w_{ijk} = 0$ for all $(i,j,k) \in \T \cap \N$.
    \item [(iv)] By~\eqref{hv} we have $\bar f^x_{ijki'} > 0$ for all $(i,j,k) \in \F \cap \P$ and $(i',j,k) \in \T \cap \P$; implying
    $w_{ijk}-w_{i'jk} = x_i$. By part~(i) $x_i =0$, while by part~(ii),
    $w_{i'jk} = 0$, implying $w_{ijk} = 0$. By assumption~(1) of Proposition~\ref{deter2},
    for any $j,k$ with $\bar y_j = \bar z_k =1$, there exists $(i',j,k) \in \T \cap \P$. Therefore, we conclude that $w_{ijk} = 0$ for all $(i,j,k) \in \F \cap \P$.

    \item [(v)] By assumption $\bar f^x_{ijki'} = \alpha^x_i > 0$ for all $(i,j,k) \in \T \cap \N$ and $(i',j,k) \in \T \cap \P$; this in turn implies that $w_{ijk}-w_{i'jk} = x_i$. By part~(ii) we have $w_{i'jk} = 0$ and by part~(iii) we have $w_{ijk} = 0$. This implies that $x_i = 0$. For any $i \in [n]$ with $\bar x_i = 0$, by condition~1 of Proposition~\ref{deter1}, $F^{x,i} > 0$ and hence by inequality~\eqref{aa1} of Proposition~\ref{deter2}, we have $T^{x,i}_{111} >  0$. That is for each $i \in [n]$ with $\bar x_i = 0$, there exists $j,k,i'$ with $\bar y_j = \bar z_k = \bar x_{i'} =1$ such that $(i,j,k) \in \T \cap \N$ and $(i',j,k) \in \T \cap \P$. Hence $x_i = 0$ for all $i \in [n]$ with $\bar x_i = 0$. By symmetry we conclude that $y_i = 0$ for all $j \in [m]$ with $\bar y_j = 0$
    and $z_k = 0$ for all $k \in [l]$ with $\bar z_k = 0$.

    \item [(vi)] By~\eqref{el} for any $(i,j,k) \in \F \cap \N$ with
    $\bar x_i = 0$, we have $w_{ijk} = x_i$.
    By part~(v) we have $x_i = 0$, implying $w_{ijk}=0$. By symmetry we conclude that $w_{ijk}=0$ for all $(i,j,k) \in \F \cap \N$.

\end{itemize}
From parts~(i)-(vi) we conclude there is no nonzero $(x,y,z,\W)$ satisfying~\eqref{uc}.
\end{prf}

\subsection{Recovery under the semi-random corruption model}

We now consider the semi-random corruption model and prove~\cref{th LP2}, which provides a sufficient condition in terms of $p, r_{\bar x}, r_{\bar y}, r_{\bar z}$ under which the flower LP recovers the ground truth with high probability. To this end, we make use of the following two lemmas:

\begin{lemma}\label{om3}
Let $X_1, X_2$ be nonnegative independent random variables with $\avg[X_1] > 0$ and $\avg[X_2] > 0$. Suppose that for every $t > 0$ we have
\begin{align}\label{wtfp}
\prob[|X_1 -\avg[X_1]| \geq t] \leq f_1(t), \quad \prob[|X_2 -\avg[X_2]| \geq t] \leq f_2(t),
\end{align}
for some functions $f_1(\cdot), f_2(\cdot)$.
Then for any
$t > 0$, we have
\begin{align*}
& \prob\Big[\Big|X_1 X_2-\avg[X_1 X_2]\Big|\geq t\Big] \\
& \qquad \leq f_1\left( \frac{t}{\sqrt{t \avg[X_2]/\avg[X_1]} + 2 \avg[X_2]} \right)
+ f_2\left( \frac{t}{\sqrt{t \avg[X_1] / \avg[X_2]} + 2 \avg[X_1]} \right).
\end{align*}
\end{lemma}
\begin{prf}
Let $\beta := 2+\sqrt{t/(\avg[X_1] \avg[X_2])}$.
We claim that if $X_1 X_2 - \avg[X_1] \avg[X_2] \geq t$,
then at least one the following inequalities hold:
$X_1\geq \avg[X_1] + \frac{t}{\beta \avg[X_2]}$,
$X_2\geq \avg[X_2] + \frac{t}{\beta \avg[X_1]}$.
To see this, assume that neither of these inequalities hold, then we have
\begin{align*}
X_1 X_2
& < \Big(\avg[X_1] + \frac{t}{\beta \avg[X_2]}\Big) \Big(\avg[X_2] + \frac{t}{\beta \avg[X_1]}\Big) \\
& = \avg[X_1]\avg[X_2]+\frac{2t}{\beta}+ \frac{t^2}{\beta^2 \avg[X_1]\avg[X_2]} \leq \avg[X_1] \avg[X_2] + t,
\end{align*}
where the last inequality follows since
$\frac{2t}{\beta} + \frac{t^2}{\beta^2 \avg[X_1]\avg[X_2]} - t \le 0$ for every $t \ge 0$.
We have obtained a contradiction, hence we have shown our claim.
It then follows that:
\begin{align*}
\prob\Big[X_1 X_2-\avg[X_1 X_2]\geq t\Big]  & \leq \prob\Big[\Big\{X_1 -\avg[X_1]\geq  \frac{t}{\beta \avg[X_2]}\Big\} \cup
\Big\{X_2 -\avg[X_2]\geq \frac{t}{\beta \avg[X_1]}\Big\}\Big]\\
& \leq \prob\Big[X_1 -\avg[X_1]\geq \frac{t}{\beta \avg[X_2]}\Big]
+ \prob\Big[X_2-\avg[X_2] \geq \frac{t}{\beta \avg[X_1]}\Big],
\end{align*}
where the first inequality follows from set inclusion, and the second inequality follows from taking the union bound.
Similarly it can be shown that
$$
\prob\Big[X_1 X_2-\avg[X_1 X_2]\leq -t\Big] \leq \prob\Big[X_1 -\avg[X_1]\leq - \frac{t}{\beta \avg[X_2]}\Big] +  \prob\Big[X_2 -\avg[X_2]\leq - \frac{t}{\beta \avg[X_1]}\Big].
$$
Therefore we have
\begin{align*}
\prob\Big[\Big|X_1 X_2-\avg[X_1 X_2]\Big|\geq t\Big] & \leq
\prob\Big[\Big|X_1-\avg[X_1]\Big|\geq \frac{t}{\beta \avg[X_2]}\Big]
+ \prob\Big[\Big|X_2-\avg[X_2]\Big|\geq \frac{t}{\beta \avg[X_1]}\Big]\\
& \leq f_1\Big(\frac{t}{\beta \avg[X_2]}\Big)+f_2\Big(\frac{t}{\beta \avg[X_1]}\Big),
\end{align*}
where the second inequality follows from inequalities~\eqref{wtfp}.
\end{prf}

\begin{lemma}\label{om}
Suppose that $X$ is a positive random variable and that we have
\begin{equation}\label{wtf1}
\prob[|X -\avg[X]| \geq t] \leq f(t),
\end{equation}
for some function $f(\cdot)$ and $t \geq 0$. Then, for every $s \ge 0$, we have
\begin{equation}\label{wtf4}
\prob\Big[\Big|\frac{1}{X}-\frac{1}{\avg[X]}\Big|\geq s\Big] \leq f(g(s)),
\end{equation}
where $g(s)=\avg[X]-\frac{1}{1/\avg[ X]+s}$.
\end{lemma}
\begin{prf}
Let $s \geq 0$. We have
\begin{align}\label{wtf2}
    \prob\Big[\frac{1}{X}-\frac{1}{\avg[X]}\geq s\Big] = \prob\Big[X\leq \frac{1}{1/\avg[X]+s}\Big]=\prob\Big[X\leq \avg[X]-g(s)\Big],
\end{align}
where the first equality follows since $X > 0$ and the second equality follows from the definition of $g(s)$. By symmetry it can be shown that
\begin{align}\label{wtf3}
    \prob\Big[\frac{1}{X}-\frac{1}{\avg[X]}\leq -s\Big] = \prob\Big[ X\geq \avg[X]+g(s)\Big].
\end{align}
Combining~\eqref{wtf1},~\eqref{wtf2}, and~\eqref{wtf3}, we obtain~\eqref{wtf4}:
\begin{align*}
    \prob\Big[\Big|\frac{1}{X}-\frac{1}{\avg[X]}\Big|\geq s\Big]
    = \prob[|X -\avg[X]| \geq g(s)]
    \le f(g(s)).
\end{align*}
\end{prf}
To proceed with the proof of Theorem~\ref{th LP2}, we introduce some random variables:
for each $i \in [n]$
and $j \in [m]$ and $r \in \{0,1\}$, denote by $t^{x,y,i,j}_{k \rightarrow r}$ (resp. $t^{x,z,i,k}_{j \rightarrow r}$ and $t^{y,z,j,k}_{i \rightarrow r}$) a random variable whose value equals 1, if  $\bar z_k =r$ (resp. $\bar y_j = r$ and $\bar x_i =r$), and $(i,j,k) \in \T$, and equals 0, otherwise.

\begin{prfc}[of \cref{th LP2}]
By~\cref{robustLPs} it suffices to prove the statement under the fully-random corruption model.
First let $p =0$. Then by Theorem~\ref{th it ub} the standard LP recovers the ground truth if $r_{\bar x}, r_{\bar y}, r_{\bar z}$ are positive. Since the feasible region of the flower LP is a subset of that of the standard LP, we conclude that for $p =0$, the flower LP recovers the ground truth tensor, if $r_{\bar x}, r_{\bar y}, r_{\bar z}$ are positive.

Henceforth, we assume that $p > 0$.
By proof of Claim~\ref{cl1}, Condition~\ref{c1-1} of Proposition~\ref{deter1} holds with high probability.
Denote by $A^1$ the event that Condition~\ref{c2-1} of Proposition~\ref{deter2} is satisfied,
denote by $A^2$ the event that Condition~\ref{c2-3} of Proposition~\ref{deter2} is satisfied,
denote by $A^3$ the event that inequalities~\eqref{aa1} and all symmetric counterparts
are strictly satisfied, and denote by $A^4$ the event that inequalities~\eqref{newass} are strictly satisfied.
To establish recovery with high probability, it suffices to show that each $A^i$ occurs
with high probability.

\begin{claim}\label{cl4}
Event $A^1$ occurs with high probability.
\end{claim}
\begin{cpf}
We have $A^1 =
\bigcup_{i=1}^3 {A^1_i}$, where the event $A^1_1$ occurs if
$\frac{1}{l} \sum_{k} {t^{x,y,i,j}_{k \rightarrow 1}} > 0$ for all $i \in [n], j \in [m]$,the event $A^1_2$ occurs if
$\frac{1}{m} \sum_{j} {t^{x,z,i,k}_{j \rightarrow 1}} > 0$ for all
$i \in [n], k \in [l]$, the event $A^1_3$ occurs if
$
\frac{1}{n} \sum_{i} {t^{y,z,j,k}_{i \rightarrow 1}} > 0$, for all
$j \in [m], k \in [l]$.
In the following, we prove that event $A^1_1$ occurs with high probability. By symmetry, it follows that $A^1_2$ and $A^1_3$ occur with high probability as well.
First notice that for each $i \in [n]$, $j \in [m]$ we have
$$
\epsilon := \avg\Big[\frac{1}{l} \sum_{k} {t^{x,y,i,j}_{k \rightarrow 1}}\Big] = \frac{1}{l} \sum_{k} {\avg[t^{x,y,i,j}_{k \rightarrow 1}]} = r_{\bar z} (1-p) >0,
$$
where the inequality follows since by assumption $r_{\bar z} > 0$
and $p < 1$. We then have
\begin{align*}
    \prob[A^1_1]
    & = \prob \Big[\bigcap_{i=1}^n \bigcap_{j=1}^m\Big\{\frac{1}{l} \sum_{k} {t^{x,y,i,j}_{k \rightarrow 1}} > 0\Big\}\Big]
    = \prob \Big[\bigcap_{i=1}^n \bigcap_{j=1}^m\Big\{\frac{1}{l} \sum_{k} {t^{x,y,i,j}_{k \rightarrow 1}} -\avg \Big[\frac{1}{l} \sum_{k} {t^{x,y,i,j}_{k \rightarrow 1}}\Big] > -\epsilon\Big\}\Big]\\
    & \geq \prob \Big[\bigcap_{i=1}^n \bigcap_{j=1}^m\Big\{\Big|\frac{1}{l} \sum_{k} {t^{x,y,i,j}_{k \rightarrow 1}} -\avg \Big[\frac{1}{l} \sum_{k} {t^{x,y,i,j}_{k \rightarrow 1}}\Big]\Big| \leq \epsilon\Big\}\Big]\\
     & \geq 1-\sum_{i,j} {\prob \Big[\Big|\frac{1}{l} \sum_{k} {t^{x,y,i,j}_{k \rightarrow 1}} -\avg \Big[\frac{1}{l} \sum_{k} {t^{x,y,i,j}_{k \rightarrow 1}}\Big]\Big| > \epsilon\Big]}
    \geq 1 - 2 n m \exp(-2l\epsilon^2),
\end{align*}
where the first inequality follows from set inclusion, the second inequality follows from taking the union bound, and the last inequality follows from the application of Hoeffding's inequality by noting that for all $i \in [n]$, $j \in [m]$ and $k \in [l]$, random variables $t^{x,y,i,j}_{k \rightarrow 1}$ are independent and $0 \leq t^{x,y,i,j}_{k \rightarrow 1} \leq 1$.
The proof then follows since $\epsilon$ is a constant and since the limit assumptions in the theorem imply that, as $n, m, l \to \infty$, we have $nm \exp (-l) \to 0$.
\end{cpf}

\begin{claim}\label{a3}
Event $A^2$ occurs with high probability for any $0 < \alpha < 1$.
\end{claim}
\begin{cpf}
To prove the statement, by symmetry, it suffices to show that for all $i \in [n]$ with $\bar x_i = 0$, with high probability, we have
\begin{align}
\label{tr1}
\frac{1}{nml}\Big(\sum_{j,k} {\Big(t^{x,i}_{jk\rightarrow 11}\sum_{i'}{t^{y,z,j,k}_{i' \rightarrow 1}}\Big)}
-\alpha \Big(\sum_{j,k}{t^{x,i}_{jk\rightarrow 11}}\Big) \Big(\frac{1}{ml}\sum_{i',j,k}{t^{y,z,j,k}_{i' \rightarrow 1}}\Big)\Big) \geq 0.
\end{align}
Denoting by $\epsilon$ the expectation of the left hand side of inequality~\eqref{tr1}, we obtain
\begin{align*}
\epsilon
&= \frac{1}{nml}\Big(\sum_{j,k} {\avg[t^{x,i}_{jk\rightarrow 11}]\sum_{i'}{\avg[t^{y,z,j,k}_{i' \rightarrow 1}}]}
-\alpha \Big(\sum_{j,k}{\avg[t^{x,i}_{jk\rightarrow 11}}]\Big) \Big(\frac{1}{ml}\sum_{i',j,k}{\avg[t^{y,z,j,k}_{i' \rightarrow 1}}]\Big)\Big) \\
&= r_{\bar y} r_{\bar z} (1-p) r_{\bar x} (1-p) - \alpha r_{\bar y} r_{\bar z} (1-p) r_{\bar x} (1-p)
= r_{\bar x} r_{\bar y} r_{\bar z} (1-p)^2 (1-\alpha),
\end{align*}
where the first equality follows from the independence of random variables since $\bar x_i =0$ while $\bar x_{i'} = 1$. Hence $\epsilon$ is a positive constant since $r_{\bar x}, r_{\bar y}, r_{\bar z}, p, \alpha$ are all positive constants and $\alpha < 1$, $p < 1$.
Let $I_0:=\{i \in [n]: \bar x_i=0\}$.
For each $i \in I_0$, define
$$
Y^{x,i}_{jk} =\frac{1}{n}t^{x,i}_{jk\rightarrow 11}\sum_{i'}{t^{y,z,j,k}_{i' \rightarrow 1}}, \quad \forall (j,k) \in [m] \times [l].
$$
Then inequality~\eqref{tr1} can be written as:
$$
\frac{1}{ml}\sum_{j,k} {Y^{x,i}_{jk}} - \alpha\Big(\frac{1}{ml}\sum_{j,k}{t^{x,i}_{jk\rightarrow 11}} \Big) \Big(\frac{1}{nml}\sum_{i',j,k}{t^{y,z,j,k}_{i' \rightarrow 1}} \Big) \geq 0.
$$
It then follows that
\begin{align*}
& \prob\Bigg[\bigcap_{i \in I_0}\Bigg\{\frac{1}{ml}\sum_{j,k} {Y^{x,i}_{jk}} - \alpha\Big(\frac{1}{ml}\sum_{j,k}{t^{x,i}_{jk\rightarrow 11}} \Big) \Big(\frac{1}{nml}\sum_{i',j,k}{t^{y,z,j,k}_{i' \rightarrow 1}} \Big) \geq 0 \Bigg\} \Bigg] \\
&= \prob\Bigg[\bigcap_{i \in I_0}\Bigg\{\frac{1}{ml}\sum_{j,k} {Y^{x,i}_{jk}}-\avg\Bigg[\frac{1}{ml}\sum_{j,k} {Y^{x,i}_{jk}}\Bigg]
-\alpha\Big(\frac{1}{ml}\sum_{j,k}{t^{x,i}_{jk\rightarrow 11}} \Big)\\
& \quad \Big(\frac{1}{nml}\sum_{i',j,k}{t^{y,z,j,k}_{i' \rightarrow 1}} \Big) + \alpha \avg\Bigg[\Big(\frac{1}{ml}\sum_{j,k}{t^{x,i}_{jk\rightarrow 11}} \Big) \Big(\frac{1}{nml}\sum_{i',j,k}{t^{y,z,j,k}_{i' \rightarrow 1}} \Big)\Bigg]\geq -\epsilon \Bigg\} \Bigg]\\
& \geq \prob\Bigg[\bigcap_{i \in I_0}\Bigg\{\Bigg|\frac{1}{ml}\sum_{j,k} {Y^{x,i}_{jk}}-\avg\Bigg[\frac{1}{ml}\sum_{j,k} {Y^{x,i}_{jk}}\Bigg]\Bigg| \leq \frac{\epsilon}{2}\Bigg\} \bigcap \bigcap_{i \in I_0} \\
& \quad \Bigg\{\Bigg|\Big(\frac{1}{ml}\sum_{j,k}{t^{x,i}_{jk\rightarrow 11}} \Big) \Big(\frac{1}{nml}\sum_{i',j,k}{t^{y,z,j,k}_{i' \rightarrow 1}} \Big) - \avg\Bigg[\Big(\frac{1}{ml}\sum_{j,k}{t^{x,i}_{jk\rightarrow 11}} \Big) \Big(\frac{1}{nml}\sum_{i',j,k}{t^{y,z,j,k}_{i' \rightarrow 1}} \Big)\Bigg]\Bigg| \leq \frac{\epsilon}{2 \alpha} \Bigg\} \Bigg]\\
&\geq 1- \sum_{i \in I_0}\prob\Bigg[\Bigg|\frac{1}{ml}\sum_{j,k} {Y^{x,i}_{jk}}-\avg\Bigg[\frac{1}{ml}\sum_{j,k} {Y^{x,i}_{jk}}\Bigg]\Bigg| > \frac{\epsilon}{2}\Bigg]\\
& \quad -\sum_{i\in I_0}{\prob\Bigg[\Bigg|\Big(\frac{1}{ml}\sum_{j,k}{t^{x,i}_{jk\rightarrow 11}} \Big) \Big(\frac{1}{nml}\sum_{i',j,k}{t^{y,z,j,k}_{i' \rightarrow 1}} \Big) - \avg\Bigg[\Big(\frac{1}{ml}\sum_{j,k}{t^{x,i}_{jk\rightarrow 11}} \Big) \Big(\frac{1}{nml}\sum_{i',j,k}{t^{y,z,j,k}_{i' \rightarrow 1}} \Big)\Bigg]\Bigg| > \frac{\epsilon}{2 \alpha} \Bigg]},
\end{align*}
where the first inequality follows from set inclusion, and the second inequality follows from taking the union bound.
Now consider the expression on the right-hand side of the last inequality; let us denote this expression by $\zeta$.  Consider the first summation in $\zeta$.
From the application of Hoeffding's inequality and using the fact that the random variables $0 \leq Y^{x,i}_{jk} \leq 1$ for all $(j,k) \in [m] \times [l]$ are independent, it follows that
\begin{equation}\label{r1}
 \sum_{i \in I_0}\prob\Bigg[\Bigg|\frac{1}{ml}\sum_{j,k} {Y^{x,i}_{jk}}-\avg\Bigg[\frac{1}{ml}\sum_{j,k} {Y^{x,i}_{jk}}\Bigg]\Bigg| > \frac{\epsilon}{2}\Bigg] \leq    2n \exp\Big(-\frac{ml \epsilon^2}{2}\Big).
\end{equation}
To bound the second summation in $\zeta$,  we use \cref{om3} by defining:
$$
Z^{x,i} := \frac{1}{ml}\sum_{j,k}{t^{x,i}_{jk\rightarrow 11}}, \; \forall i \in [n] \text{ with } \bar x_i = 0,
\qquad
W := \frac{1}{nml}\sum_{i',j,k}{t^{y,z,j,k}_{i' \rightarrow 1}}.
$$
Notice that $Z^{x,i}$ and $W$ are nonnegative independent random variables with $\avg[Z]:=\avg[Z^{x,i}] = r_{\bar y} r_{\bar z} (1-p) > 0$ and $\avg[W] = r_{\bar x} (1-p) > 0$. Moreover, by Hoeffding's inequality we have:
\begin{align*}
\prob\Big[|Z^{x,i}-\avg[Z^{x,i}]| > \frac{\epsilon}{2 \alpha}\Big] \leq 2 \exp\Big(-\frac{ml\epsilon^2}{2 \alpha^2}\Big),
\qquad \prob\Big[|W-\avg[W]| > \frac{\epsilon}{2 \alpha}\Big] \leq  \exp\Big(-\frac{nml\epsilon^2}{2 \alpha^2}\Big).
\end{align*}
Hence utilizing \cref{om3} yields
\begin{align}
\label{r2}
\begin{split}
\sum_{i \in I_0}{\prob\Big[|Z^{x,i} W-\avg[Z^{x,i}W]| > \frac{\epsilon}{2 \alpha}\Big]}
& \leq 2n \exp\Big(-\frac{ml\epsilon^2}{2 \alpha^2 (\sqrt{2\alpha\epsilon \avg[W]/\avg[Z]}+2\avg[W])^2}\Big) \\
& \quad +2n \exp\Big(-\frac{nml\epsilon^2}{2 \alpha^2 (\sqrt{2\alpha\epsilon \avg[Z]/\avg[W]}+2\avg[Z])^2}\Big).
\end{split}
\end{align}
Combining~\eqref{r1} and~\eqref{r2},
the proof then follows since $\epsilon, \alpha, \avg[W], \avg[Z]$ are positive constants and since the limit assumptions in the theorem imply that, as $n, m, l \to \infty$, we have $n \exp (-ml)$, $n \exp (-nml)$ go to zero.
\end{cpf}

\begin{claim}\label{a2}
Event $A^3$ occurs with high probability.
\end{claim}

\begin{cpf}
Denote by $A^3_1$ the event that inequalities~\eqref{aa1} are strictly satisfied.
By symmetry, to show that $A^3$ occurs with high probability, it suffices to show that $A^3_1$ occurs with high probability. Under the random corruption model,
$A^3_1$ occurs if, (i) inequalities~\eqref{frp} are satisfied and (ii) for each $(i,j,k) \in \T$ with $\bar x_i = 0$, $\bar y_j = \bar
z_k = 1$, we have
\begin{align}\label{weve1}
\frac{1}{ml}\sum_{j',k'} {Y^{x,i}_{j'k'}}-
\frac{1}{nml}\sum_{j',k'} {\Big(\frac{1}{3} f^{x,i}_{j'k' \rightarrow 00}+\frac{1}{2} f^{x,i}_{j'k' \rightarrow 01}+\frac{1}{2} f^{x,i}_{j'k' \rightarrow 10} +f^{x,i}_{j'k' \rightarrow 11}}\Big) \sum_{i'}{t^{y,z,j,k}_{i' \rightarrow 1}} > 0,
\end{align}
where we define
$$
Y^{x,i}_{j'k'} := \frac{1}{n} t^{x,i}_{j',k'\rightarrow 11}\sum_{i'}{t^{y,z,j',k'}_{i' \rightarrow 1}}, \qquad \forall j' \in [m], k' \in [l].
$$
By proof of~\cref{cl2}, inequalities~\eqref{frp} are satisfied with high probability, if inequality~\eqref{cond1} holds. It is simple to check that inequality~\eqref{cond1} is implied by inequality~\eqref{condaux simple}. We now show that if~\eqref{cond1} holds, inequalities~\eqref{weve1} are satisfied with high probability as well.
Denote by $\epsilon$ the expected value of the left-hand side of inequality~\eqref{weve1}. We have
\begin{align*}
    \epsilon
    &= \frac{1}{nml}\sum_{j',k'} {\avg\Big[t^{x,i}_{j',k'\rightarrow 11}\sum_{i'}{t^{y,z,j',k'}_{i' \rightarrow 1}}\Big]}\\
    &\quad-\frac{1}{nml}\avg\Bigg[\sum_{j',k'} {\Big(\frac{1}{3} f^{x,i}_{j'k' \rightarrow 00}+\frac{1}{2} f^{x,i}_{j'k' \rightarrow 01}+\frac{1}{2} f^{x,i}_{j'k' \rightarrow 10} +f^{x,i}_{j'k' \rightarrow 11}}\Big) \sum_{i'}{t^{y,z,j,k}_{i' \rightarrow 1}}\Bigg]\\
&= r_{\bar x} r_{\bar y} r_{\bar z} (1-p)^2 -\frac{1}{3}r_{\bar x} (1-r_{\bar y}) (1-r_{\bar z})p (1-p) -\frac{1}{2}r_{\bar x} (1-r_{\bar y})r_{\bar z} p (1-p)\\
&\quad -\frac{1}{2}r_{\bar x} r_{\bar y} (1-r_{\bar z}) p(1-p) - r_{\bar x} r_{\bar y} r_{\bar z} p (1-p),
\end{align*}
where the second equality follows from the independence of random variables as we have $\bar x_{i'} =1$ while
$\bar x_{i} = 0$. Since by assumption $r_{\bar x} > 0$ and $p < 1$, the inequality $\epsilon > 0$ can be equivalently written as:
$$
r_{\bar y} r_{\bar z} (1-2 p) -\frac{1}{3}(1-r_{\bar y}) (1-r_{\bar z})p -\frac{1}{2}(1-r_{\bar y})r_{\bar z} p -\frac{1}{2}r_{\bar y} (1-r_{\bar z}) p
- r_{\bar y} r_{\bar z} p >0,
$$
which is in turn equivalent to inequality~\eqref{cond1}.
 Define
 $$
 \bar f^{x,i}_{jk} =\frac{1}{3} f^{x,i}_{jk \rightarrow 00}+\frac{1}{2} f^{x,i}_{jk \rightarrow 01}+\frac{1}{2} f^{x,i}_{jk \rightarrow 10} +f^{x,i}_{jk \rightarrow 11}, \qquad \forall (i,j,k) \in [n] \times [m] \times [l],
 $$
and $\M := \{(i,j,k): (i,j,k) \in \T, \; \bar x_i = 0, \; \bar y_j = \bar z_k = 1\}$. Then we have:
\begin{align*}
&\prob[A^3_1]= \prob\Bigg[\bigcap_{\substack{(i,j,k) \in \M}} \Bigg\{\frac{1}{ml}\sum_{j',k'} {Y^{x,i}_{j'k'}}-
\Big(\frac{1}{ml}\sum_{j',k'} \bar f^{x,i}_{j'k'}\Big) \Big(\frac{1}{n}\sum_{i'}{t^{y,z,j,k}_{i' \rightarrow 1}}\Big) > 0\Bigg\} \Bigg]\\
& = \prob\Bigg[\bigcap_{\substack{(i,j,k) \in \M}} \Bigg\{\frac{1}{ml}\sum_{j',k'} {Y^{x,i}_{j'k'}}-\avg\Big[\frac{1}{ml}\sum_{j',k'} {Y^{x,i}_{j'k'}}\Big]-
\Big(\frac{1}{ml}\sum_{j',k'} \bar f^{x,i}_{j'k'}\Big) \Big(\frac{1}{n}\sum_{i'}{t^{y,z,j,k}_{i' \rightarrow 1}}\Big) \\
& \quad + \avg\Big[\Big(\frac{1}{ml}\sum_{j',k'} \bar f^{x,i}_{j'k'}\Big) \Big(\frac{1}{n}\sum_{i'}{t^{y,z,j,k}_{i' \rightarrow 1}}\Big)\Big]> -\epsilon\Bigg\} \Bigg]\\
& \geq \prob\Bigg[\bigcap_{\substack{(i,j,k) \in \M}} \Bigg\{\Big|\frac{1}{ml}\sum_{j',k'} {Y^{x,i}_{j'k'}}-\avg\Big[\frac{1}{ml}\sum_{j',k'} {Y^{x,i}_{j'k'}}\Big]\Big| < \frac{\epsilon}{2}\Bigg\}\Bigg] \cdot \\
& \quad \prob\Bigg[\bigcap_{\substack{(i,j,k) \in \M}} \Bigg\{\Big|\Big(\frac{1}{ml}\sum_{j',k'} \bar f^{x,i}_{j'k'}\Big) \Big(\frac{1}{n}\sum_{i'}{t^{y,z,j,k}_{i' \rightarrow 1}}\Big)
- \avg\Big[\Big(\frac{1}{ml}\sum_{j',k'} \bar f^{x,i}_{j'k'}\Big) \Big(\frac{1}{n}\sum_{i'}{t^{y,z,j,k}_{i' \rightarrow 1}}\Big)\Big]\Big| < \frac{\epsilon}{2}\Bigg\}\Bigg]\\
& \geq 1- \sum_{(i,j,k) \in \M }\prob\Bigg[\Bigg|\frac{1}{ml}\sum_{j',k'} {Y^{x,i}_{j'k'}}-\avg\Big[\frac{1}{ml}\sum_{j',k'} {Y^{x,i}_{j'k'}}\Big]\Bigg| > \frac{\epsilon}{2} \Bigg]\\
& \quad - \sum_{(i,j,k) \in \M }\prob\Bigg[\Bigg|\Big(\frac{1}{ml}\sum_{j',k'} \bar f^{x,i}_{j'k'}\Big) \Big(\frac{1}{n}\sum_{i'}{t^{y,z,j,k}_{i' \rightarrow 1}}\Big)
- \avg\Big[\Big(\frac{1}{ml}\sum_{j',k'} \bar f^{x,i}_{j'k'}\Big) \Big(\frac{1}{n}\sum_{i'}{t^{y,z,j,k}_{i' \rightarrow 1}}\Big)\Big]\Bigg| > \frac{\epsilon}{2} \Bigg] := \zeta,
\end{align*}
where the first inequality follows from set inclusion, and the second inequality follows from taking the union bound. To bound $\zeta$, first note that $Y^{x,i}_{j',k'}$, $j' \in [m]$, $k' \in [l]$ are independent random variables and $0 \leq Y^{x,i}_{j',k'} \leq 1$. Hence, utilizing
Hoeffding's inequality, we obtain:
\begin{equation*}
\prob\Bigg[\Bigg|\frac{1}{ml}\sum_{j',k'} {Y^{x,i}_{j'k'}}-\avg\Big[\frac{1}{ml}\sum_{j',k'} {Y^{x,i}_{j'k'}}\Big]\Bigg| > \frac{\epsilon}{2} \Bigg] \leq
2\exp\Big(-\frac{ml \epsilon^2}{2}\Big).
\end{equation*}
To bound the terms in the second summation, we make use of \cref{om3} by defining
$$
Z_1 := \frac{1}{ml}\sum_{j',k'} \bar f^{x,i}_{j'k'},
\qquad
Z_2 := \frac{1}{n}\sum_{i'}{t^{y,z,j,k}_{i' \rightarrow 1}}.
$$
First note that $Z_1, Z_2$ are nonnegative and independent random variables with
\begin{equation*}
\avg[Z_1] = \frac{p}{3} (r_{\bar y} r_{\bar z} +\frac{r_{\bar y}+r_{\bar z}}{2} +1 )> 0, \qquad \avg[Z_2] = r_{\bar x} (1-p) > 0.
\end{equation*}
Moreover, the random variables $\bar f^{x,i}_{j'k'}$ for all $(j',k') \in [m] \times [l]$ are independent and  $0 \leq \bar f^{x,i}_{j'k'} \leq 1$. Similarly, the random variables $t^{y,z,j,k}_{i' \rightarrow 1}$ for all $i' \in [n]$ are independent and $0 \leq t^{y,z,j,k}_{i' \rightarrow 1} \leq 1$. Hence applying Hoeffding's inequality we obtain:
$$
    \prob\Big[|Z_1-\avg[Z_1]| > \frac{\epsilon}{2}\Big] \leq 2 \exp\Big(-\frac{ml\epsilon^2}{2}\Big), \qquad
    \prob\Big[|Z_2-\avg[Z_2]| > \frac{\epsilon}{2}\Big] \leq 2 \exp\Big(-\frac{n\epsilon^2}{2}\Big).
$$
Utilizing~\cref{om3} we obtain
$$
\prob\Big[|Z_1 Z_2-\avg[Z_1 Z_2]| > \frac{\epsilon}{2}\Big] \leq 2
\exp\Bigg(-\frac{ml\epsilon^2}{\Big(\sqrt{\frac{2\epsilon\avg[Z_2]}{\avg[Z_1]}}+4\avg[Z_2]\Big)^2}\Bigg)+2 \exp\Bigg(-\frac{n\epsilon^2}{\Big(\sqrt{\frac{2\epsilon\avg[Z_1]}{\avg[Z_2]}}+4\avg[Z_1]\Big)^2}\Bigg).
$$
Therefore, we have
$$
\zeta \geq 1 - 2nml\exp\Bigg(-\frac{ml\epsilon^2}{\Big(\sqrt{\frac{2\epsilon\avg[Z_2]}{\avg[Z_1]}}+4\avg[Z_2]\Big)^2}\Bigg)-2nml \exp\Bigg(-\frac{n\epsilon^2}{\Big(\sqrt{\frac{2\epsilon\avg[Z_1]}{\avg[Z_2]}}+4\avg[Z_1]\Big)^2}\Bigg).
$$
The proof then follows since $\epsilon, \avg[Z_1], \avg[Z_2]$ are positive constants and since the limit assumptions in the theorem imply that, as $n, m, l \to \infty$, we have $nml \exp (-n)$, $nml \exp (-ml)$ go to zero.
%
\end{cpf}

\begin{claim}\label{whp6}
Event $A^4$ occurs with high probability.
\end{claim}
\begin{cpf}
By Condition~2 in \cref{deter2}, $\bar \gamma_{ijk}$ can be equivalently written as
\begin{align*}
\bar \gamma_{ijk} & =\frac{1}{3}\Bigg(2 - \frac{1}{3}\frac{n_{\bar x}}{T^{y,z,j,k}_1}
    -\frac{1}{3}\frac{n_{\bar y}}{T^{x,z,i,k}_1}-\frac{1}{3}\frac{n_{\bar z}}{T^{x,y,i,j}_1}\nonumber\\
    & \quad -\frac{1}{\alpha \bar T^{y,z}_1}\sum_{\substack{i': \bar x_{i'} = 0\\ (i',j,k) \in \T}}{\min\Bigg\{\frac{1}{T_{11}^{x,i'}}\Big(\frac{1}{3}F^{x,i'}_{00}+\frac{1}{2}(F^{x,i'}_{01}+F^{x,i'}_{10})+F^{x,i'}_{11}\Big),\;1\Bigg\}} \nonumber\\
    & \quad -\frac{1}{\alpha \bar T^{x,z}_1}\sum_{\substack{j': \bar y_{j'} =0\\ (i,j',k) \in \T}} {\min\Bigg\{ \frac{1}{T_{11}^{y,j'}}\Big(\frac{1}{3}F^{y,j'}_{00}+\frac{1}{2}(F^{y,j'}_{01}+F^{y,j'}_{10})+F^{y,j'}_{11}\Big),\;1\Bigg\}}\nonumber\\
    & \quad -\frac{1}{\alpha \bar T^{x,y}_1}\sum_{\substack{k': \bar z_{k'} =0 \\ (i,j,k') \in \T}}{\min\Bigg\{  \frac{1}{T_{11}^{z,k'}}\Big(\frac{1}{3}F^{z,k'}_{00}+\frac{1}{2}(F^{z,k'}_{01}+F^{z,k'}_{10})+F^{z,k'}_{11}\Big),\;1\Bigg\}}\Bigg).
\end{align*}
In the following, we show the validity of the inequality:
\begin{equation} \label{nequ}
    3 \bar \gamma_{ijk} - \frac{n_{\bar y} n_{\bar z}}{T^{x,i}_{11}} > -1, \qquad \forall (i,j,k) \in \T \cap \P.
\end{equation}
By symmetry, this in turn implies event $A^4$ occurs with high probability.
Under the random corruption model, inequalities~\eqref{nequ} are satisfied, if for each $(i,j,k) \in \T \cap \P$:
%
\begin{align}
\label{hgamma}
\begin{split}
   &3-\frac{1}{3}\frac{n_{\bar x}}{\sum_{i'}{t^{y,z,j,k}_{i' \rightarrow 1}}}
   -\frac{1}{3}\frac{n_{\bar y}}{\sum_{j'}{t^{x,z,i,k}_{j' \rightarrow 1}}}
   -\frac{1}{3}\frac{n_{\bar z}}{\sum_{k'}{t^{x,y,i,j}_{k' \rightarrow 1}}}-\frac{n_{\bar y} n_{\bar z}}{\sum_{j',k'}{t^{x,i}_{j'k'\rightarrow 1}}} \\
    & \quad -\frac{1}{\alpha}\frac{nml}{\sum_{i',j', k'}{t^{y,z,j', k'}_{i' \rightarrow 1}}} \Big(\frac{\sum_{i'}{\bar\nu^{x,i'}_{jk}}}{n} \Big)
    -\frac{1}{\alpha}\frac{nml}{\sum_{i',j', k'}{t^{x,z,i', k'}_{j' \rightarrow 1}}}\Big(\frac{\sum_{j'}{\bar\nu^{y,j'}_{ik}}}{m}\Big) \\
    & \quad -\frac{1}{\alpha}\frac{nml}{\sum_{i',j', k'}{t^{x,y,i', j'}_{k' \rightarrow 1}}}\Big(\frac{\sum_{k'}{\bar\nu^{z,k'}_{ij}}}{l}\Big) > 0,
    \end{split}
\end{align}
where $\bar \nu^{x,i'}_{jk} = \min\{\nu^{x,i'}_{jk}, 1\}$, $\bar \nu^{y,j'}_{ik} = \min\{\nu^{y,j'}_{ik}, 1\}$, $\bar \nu^{z,k'}_{ij} = \min\{\nu^{z,k'}_{ij}, 1\}$, and $\nu^{x,i'}_{jk}, \nu^{y,j'}_{ik},\nu^{z,k'}_{ij}$ are defined by~\eqref{defnu}.

First, we observe that if we have \eqref{condaux simple},
then we also have
\begin{equation}
\label{condaux}
p < \frac{\alpha r_{\bar x} r_{\bar y} r_{\bar z}}{1+(3\alpha-1)r_{\bar x} r_{\bar y} r_{\bar z}},
\end{equation}
for some $\alpha < 1$ arbitrarily close to 1.
This is because the function $\alpha \mapsto \frac{\alpha}{1+(3\alpha-1)r_{\bar x} r_{\bar y} r_{\bar z}}$ is continuous in $[0,1]$.

Denote by $\avg_g$ the expected value of the left-hand side of inequality~\eqref{hgamma}. To prove the statement, we first show that
for each $(i,j,k) \in \T \cap \P$, we have
\begin{equation}\label{gala}
\avg_g \geq \epsilon:=
3 - \frac{2}{1-p}-\frac{1}{\alpha}\frac{p}{(1-p)}\Big(\frac{1}{r_{\bar x} r_{\bar y} r_{\bar z}}-1\Big),
\end{equation}
which implies $\avg_g > 0$ if condition~\eqref{condaux} is satisfied.
We prove~\eqref{gala} via a number of steps:
\begin{step}
\label{cl5}
For each $(i,j,k) \in \T \cap \P$, we have
$$
\avg\Bigg[{\frac{n_{\bar x}}{\sum_{i'}{t^{y,z,j,k}_{i' \rightarrow 1}}}}\Bigg] \leq \frac{1}{1-p}.
$$
\end{step}
\begin{spf}
For each $(i,j,k) \in \T \cap \P$ we have
\begin{align*}
    & \avg\Bigg[{\frac{n_{\bar x}}{\sum_{i'}{t^{y,z,j,k}_{i' \rightarrow 1}}}}\Bigg]= n_{\bar x} \avg\Bigg[\frac{1}{1+\sum_{i' \in [n] \setminus \{i\}}{t^{y,z,j,k}_{i' \rightarrow 1}}}\Bigg] = \frac{n_{\bar x}}{n_{\bar x} (1-p)} (1-p^{n_{\bar x}}) \leq \frac{1}{1-p},
\end{align*}
where the first equality follows since by assumption $(i,j,k) \in \T \cap \P$, \ie $t^{y,z,j,k}_{i \rightarrow 1} = 1$ , the second equality follows since
for a binomial random variable $X$ with parameters $(n,p)$ we have
$\avg[\frac{1}{1+X}] = \frac{1}{(n+1)p}(1-(1-p)^{n+1})$, and the inequality follows since $0 < p \leq 1$.
\end{spf}

\begin{step}\label{cl6}
For each $(i,j,k) \in \T \cap \P$ we have
\begin{align*}
\epsilon' & := \avg \Bigg[\frac{nml}{\sum_{i',j', k'}{t^{y,z,j', k'}_{i' \rightarrow 1}}} \Big(\frac{\sum_{i'}{\bar\nu^{x,i'}_{jk}}}{n} \Big)\Bigg]
\leq \frac{1}{3}\frac{p}{(1-p)} (\frac{1}{r_{\bar x}}-1) \Big(\frac{1}{r_{\bar y} r_{\bar z}}+\frac{1}{2 r_{\bar y}}+\frac{1}{2 r_{\bar z}}+1\Big).
\end{align*}
\end{step}
\begin{spf}
For each $(i,j,k) \in \T \cap \P$, we have:
\begin{align*}
    \epsilon' & \leq \avg \Bigg[\frac{nml}{\sum_{i',j', k'}{t^{y,z,j', k'}_{i' \rightarrow 1}}} \Big(\frac{\sum_{i'}{\nu^{x,i'}_{jk}}}{n} \Big)\Bigg]
    = \avg\Bigg[\frac{nml}{\sum_{i',j', k'}{t^{y,z,j', k'}_{i' \rightarrow 1}}}\Bigg] \avg\Bigg[\frac{\sum_{i'}{\nu^{x,i'}_{jk}}}{n}\Bigg]\\
    & = \avg\Bigg[\frac{nml}{1+\sum_{(i',j', k') \in [n] \times [m] \times [l] \setminus \{(i,j,k)\}}{t^{y,z,j', k'}_{i' \rightarrow 1}}}\Bigg] \frac{\sum_{i'}{\avg[\nu^{x,i'}_{jk}]}}{n}\\
    & \leq \Big(\frac{1}{r_{\bar x}}\frac{1-p^{n_{\bar x}}}{1-p}\Big) (1-r_{\bar x})\frac{p}{3} \Big(\frac{1}{r_{\bar y} r_{\bar z}}+\frac{1}{2 r_{\bar y}}+\frac{1}{2 r_{\bar z}}+1\Big)\\
     & \leq \frac{1}{(1-p)} (\frac{1}{r_{\bar x}}-1) \frac{p}{3} \Big(\frac{1}{r_{\bar y} r_{\bar z}}+\frac{1}{2 r_{\bar y}}+\frac{1}{2 r_{\bar z}}+1\Big),
\end{align*}
where the first inequality follows from the definition of $\bar\nu^{x,i'}_{jk}$,
the first equality follows from the independence of random variables, the second equality follows since by assumption $(i,j,k) \in \T \cap \P$, the second inequality follows since for a binomial random variable $X$ with parameters $(n,p)$ we have
$\avg[\frac{1}{1+X}] = \frac{1}{(n+1)p}(1-(1-p)^{n+1})$,  and by proof of \cref{cl3} we have
$$\avg[\nu^{x,i'}_{jk}] \leq \frac{p}{3} \Big(\frac{1}{r_{\bar y} r_{\bar z}}+\frac{1}{2 r_{\bar y}}+\frac{1}{2 r_{\bar z}}+1\Big),$$
and the last inequality follows since $0 \leq p \leq 1$ and $n_{\bar x}> 0$.
\end{spf}
Therefore, by Steps~\ref{cl5} and~\ref{cl6}:
\begin{align*}
 \avg_g & \geq \epsilon:= 2 - \frac{1}{1-p} -\frac{1}{3\alpha}\frac{ p}{(1-p)}\Bigg(
  (\frac{1}{r_{\bar x}}-1) \Big(\frac{1}{r_{\bar y} r_{\bar z}}+\frac{1}{2 r_{\bar y}}+\frac{1}{2 r_{\bar z}}+1\Big)+\\
 & \quad (\frac{1}{r_{\bar y}}-1) \Big(\frac{1}{r_{\bar x} r_{\bar z}}+\frac{1}{2 r_{\bar x}}+\frac{1}{2 r_{\bar z}}+1\Big)+(\frac{1}{r_{\bar z}}-1) \Big(\frac{1}{r_{\bar x} r_{\bar y}}+\frac{1}{2 r_{\bar x}}+\frac{1}{2 r_{\bar y}}+1\Big)\Bigg)-\frac{1}{1-p}\\
 & = 3 - \frac{2}{1-p}-\frac{1}{\alpha}\frac{p}{(1-p)}\Big(\frac{1}{r_{\bar x} r_{\bar y} r_{\bar z}}-1\Big).
\end{align*}
It then follows that if condition~\eqref{condaux} holds,
we have $\avg_g \geq \epsilon > 0$.

We now show that inequalities~\eqref{hgamma} are satisfied with high probability; to this end,
utilizing \cref{om}, we first show that the first four terms in inequalities~\eqref{hgamma} concentrate around their expectations:
\begin{step}\label{om2}
Let $\epsilon' > 0$. Then
$$
\sum_{(i,j,k) \in \T \cap \P}\prob\Bigg[ \Bigg|{\frac{n_{\bar z}}{\sum_{k'}{t^{x,y,i,j}_{k' \rightarrow 1}}}}-\avg\Bigg[{\frac{n_{\bar z}}{\sum_{k'}{t^{x,y,i,j}_{k' \rightarrow 1}}}}\Bigg]\Bigg| > \epsilon'\Bigg] \rightarrow 0, \qquad {\rm as} \quad n,m,l \rightarrow \infty.
$$
\end{step}
\begin{spf}
By a simple application of Hoeffding's inequality, for any $t \geq 0$, we have
\begin{equation}\label{eom2}
\prob\Bigg[ \Bigg|{\frac{\sum_{k'}{t^{x,y,i,j}_{k' \rightarrow 1}}}{n_{\bar z}}}-\avg\Bigg[{\frac{\sum_{k'}{t^{x,y,i,j}_{k' \rightarrow 1}}}{n_{\bar z}}}\Bigg]\Bigg| > t\Bigg] \leq 2 \exp(-2n_{\bar z} t^2).
\end{equation}
For each $(i,j,k) \in \T \cap \P$, define:
$$
\Delta := \Bigg|\frac{n_{\bar z}}{\avg\Big[\sum_{k'}{t^{x,y,i,j}_{k' \rightarrow 1}}\Big]} -\avg\Bigg[{\frac{n_{\bar z}}{\sum_{k'}{t^{x,y,i,j}_{k' \rightarrow 1}}}}\Bigg]\Bigg|= \Bigg|\frac{1}{1- p}-\frac{1-p^{n_{\bar z}}}{1-p}\Bigg|=
\frac{p^{n_{\bar z}}}{1-p}.
$$
Note that since $0 < p < 1$, it follows that $\lim_{n_{\bar z}\rightarrow \infty} \Delta = 0$. Then we have
\begin{align*}
&\sum_{(i,j,k) \in \T \cap \P}\prob\Bigg[ \Bigg|{\frac{n_{\bar z}}{\sum_{k'}{t^{x,y,i,j}_{k' \rightarrow 1}}}}-\avg\Bigg[{\frac{n_{\bar z}}{\sum_{k'}{t^{x,y,i,j}_{k' \rightarrow 1}}}}\Bigg]\Bigg| > \epsilon'\Bigg] \\
& = \sum_{(i,j,k) \in \T \cap \P}\prob\Bigg[ \Bigg|{\frac{n_{\bar z}}{\sum_{k'}{t^{x,y,i,j}_{k' \rightarrow 1}}}}-{\frac{n_{\bar z}}{\avg\Big[\sum_{k'}{t^{x,y,i,j}_{k' \rightarrow 1}}\Big]}}+{\frac{n_{\bar z}}{\avg\Big[\sum_{k'}{t^{x,y,i,j}_{k' \rightarrow 1}}\Big]}}-\avg\Bigg[{\frac{n_{\bar z}}{\sum_{k'}{t^{x,y,i,j}_{k' \rightarrow 1}}}}\Bigg]\Bigg| > \epsilon'\Bigg]\\
& \leq \sum_{(i,j,k) \in \T \cap \P}\prob\Bigg[ \Bigg|{\frac{n_{\bar z}}{\sum_{k'}{t^{x,y,i,j}_{k' \rightarrow 1}}}}-{\frac{n_{\bar z}}{\avg\Big[\sum_{k'}{t^{x,y,i,j}_{k' \rightarrow 1}}\Big]}}\Bigg|+\Delta > \epsilon'\Bigg]
\leq 2nml \exp\Big(-2l r_{\bar z} (g(\epsilon'-\Delta))^2\Big),
\end{align*}
where the first inequality follows from the application of triangle inequality, and the second inequality follows from inequality~\eqref{eom2} and \cref{om} by noting that
since $(i,j,k) \in \T \cap \P$, we have
$\sum_{k'}{t^{x,y,i,j}_{k' \rightarrow 1}} = 1+ \sum_{k' \in [l] \setminus \{k\}}{t^{x,y,i,j}_{k' \rightarrow 1}} > 0$.
Since $\Delta \rightarrow 0$ as $l \rightarrow \infty$ and since $(g(\epsilon'-\Delta))^2$ is a positive constant, the proof follows since the limit assumptions in the theorem imply that, as $n, m, l \to \infty$, we have $nml\exp (-l)$ go to zero.
\end{spf}

Next, utilizing \cref{om3,om}, we show that the last three terms in inequalities~\eqref{hgamma} concentrate around their expectation:
\begin{step}\label{om4}
Let $ \epsilon'$ be a positive constant.
Define
$$
\F_{ijk} = \prob\Bigg[\Bigg|\frac{nml}{\sum_{i',j',k'}{t^{y,z,j', k'}_{i' \rightarrow 1}}} \Big(\frac{\sum_{i'}{\bar\nu^{x,i'}_{jk}}}{n} \Big)-\avg\Bigg[\frac{nml}{\sum_{i',j',k'}{t^{y,z,j', k'}_{i' \rightarrow 1}}} \Big(\frac{\sum_{i'}{\bar\nu^{x,i'}_{jk}}}{n} \Big)\Bigg]\Bigg| > \epsilon' \Bigg].
$$
Then
$$
\lim_{n,m.l \rightarrow \infty}\sum_{(i,j,k) \in \T \cap \P} {\F_{ijk}}=0.
$$
\end{step}
\begin{spf}
Define the random variables
$$
X = \frac{nml}{\sum_{i',j',k'}{t^{y,z,j', k'}_{i' \rightarrow 1}}},
\qquad Y_{jk} = \frac{\sum_{i'}{\bar\nu^{x,i'}_{jk}}}{n}.
$$
Observe that $X$ and $Y_{jk}$ are independent random variables since by definition of $t^{y,z,j', k'}_{i' \rightarrow 1}$, we have $\bar x_{i'} =1$ while by definition of  $\bar\nu^{x,i'}_{jk}$ we have $\bar x_{i'} = 0$. Moreover, $X > 0$ since by assumption $(i,j,k) \in \T \cap \P$ implying $t^{y,z,j,k}_{i \rightarrow 1} = 1$; clearly $Y_{jk} \geq 0$ and $\avg[Y_{jk}] \geq 0$. By \cref{om} and the proof of \cref{om2}, we have
\begin{equation*}
\prob\Big[ \Big|X - \avg[X]\Big| > \epsilon'\Big] \leq 2 \exp\Big(-2nml (g(\epsilon'-\Delta))^2\Big),
\end{equation*}
where
$$
\Delta := \Bigg|\avg\Bigg[\frac{nml}{\sum_{i',j',k'}{t^{y,z,j', k'}_{i' \rightarrow 1}}}\Bigg]-\frac{nml}{\avg\Big[\sum_{i',j',k'}{t^{y,z,j', k'}_{i' \rightarrow 1}}\Big]}\Bigg| =
\frac{1}{r_{\bar x}} \frac{p^{mln_{\bar x}}}{1-p}.
$$
Note that since $r_{\bar x}$ is a positive constant, we have $\lim_{n,m,l\rightarrow \infty} \Delta = 0$. Define $\bar Y = \avg[Y_{jk}]$. Since $\bar\nu^{x,i'}_{jk}$ for all $i' \in [n]$
are independent random variables and $0 \leq \bar\nu^{x,i'}_{jk} \leq 1$,
by a simple application of Hoeffding's inequality we have
\begin{equation*}
\prob\Big[ \Big|Y_{jk} - \bar Y \Big| > \epsilon'\Big] \leq 2 \exp(-2n\epsilon'^2),
\end{equation*}
Then
\begin{align*}
& \sum_{(i,j,k) \in \T \cap \P} {\F_{ijk}}= \sum_{(i,j, k) \in \T \cap \P}{\prob[|X Y_{jk}-\avg[X]\avg[Y_{jk}]| > \epsilon']}\\
& \leq 2nml \exp\Big(-2nml \Big(g\Big(\frac{\epsilon'}{\sqrt{\epsilon' \bar Y/\avg[X]}+2\bar Y}-\Delta\Big)\Big)^2\Big)+2nml \exp\Big(-\frac{2n\epsilon'^2}{(\sqrt{\epsilon'\avg[X]/\bar Y}+2 \avg[X])^2}\Big),
\end{align*}
where the inequality follows from \cref{om3}. The proof then follows since
$\Delta \rightarrow 0$ as $n,m,l \rightarrow \infty$,
$\epsilon', \avg[X], \bar Y, (g(\cdot))^2$ are positive constants and since the limit assumptions in the theorem imply that, as $n, m, l \to \infty$, we have $nml \exp (-nml)$, $nml \exp (-n)$ go to zero.
\end{spf}

We are now ready to prove that inequalities~\eqref{hgamma} hold with high probability. For each $(i,j,k) \in \T \cap \P$, denote by $\kappa_{ijk}$
the left hand side of inequality~\eqref{hgamma}.
Let $\tilde \kappa_{ijk} = \kappa_{ijk}-3$. Then we have
\begin{align*}
& \prob\Big[\bigcap_{(i,j,k) \in \T \cap \P}\{\kappa_{ijk} > 0\}\Big] \geq
\prob\Big[\bigcap_{(i,j,k) \in \T \cap \P}\{\tilde \kappa_{ijk}-\avg[\tilde \kappa_{ijk}] > -\epsilon\}\Big]\\
& = \prob\Bigg[\bigcap_{(i,j,k) \in \T \cap \P}\Bigg\{\frac{1}{3}\Bigg(\avg\Big[{\frac{n_{\bar x}}{\sum_{i'}{t^{y,z,j,k}_{i' \rightarrow 1}}}}\Big]-{\frac{n_{\bar x}}{\sum_{i'}{t^{y,z,j,k}_{i' \rightarrow 1}}}}\Bigg) +\frac{1}{3}\Bigg(\avg\Big[{\frac{n_{\bar y}}{\sum_{j'}{t^{x,z,i,k}_{j' \rightarrow 1}}}}\Big]-{\frac{n_{\bar y}}{\sum_{j'}{t^{x,z,i,k}_{j' \rightarrow 1}}}}\Bigg)\\
&\quad +\frac{1}{3}\Bigg(\avg\Big[{\frac{n_{\bar z}}{\sum_{k'}{t^{x,y,i,j}_{k' \rightarrow 1}}}}\Big]-{\frac{n_{\bar z}}{\sum_{k'}{t^{x,y,i,j}_{k' \rightarrow 1}}}}\Bigg)+\Bigg(\avg\Big[\frac{n_{\bar y} n_{\bar z}}{\sum_{j',k'}{t^{x,i}_{j'k'\rightarrow 1}}}\Big]-\frac{n_{\bar y} n_{\bar z}}{\sum_{j',k'}{t^{x,i}_{j'k'\rightarrow 1}}}\Bigg)\\
&\quad +\frac{1}{\alpha}\Bigg(\avg\Bigg[\frac{nml}{\sum_{i',j',k'}{t^{y,z,j', k'}_{i' \rightarrow 1}}} \Big(\frac{\sum_{i'}{\bar\nu^{x,i'}_{jk}}}{n} \Big)\Bigg]-\frac{nml}{\sum_{i',j',k'}{t^{y,z,j', k'}_{i' \rightarrow 1}}} \Big(\frac{\sum_{i'}{\bar\nu^{x,i'}_{jk}}}{n} \Big)\Bigg)\\
&\quad +\frac{1}{\alpha}\Bigg(\avg\Bigg[\frac{nml}{\sum_{i',j',k'}{t^{x,z,i', k'}_{j' \rightarrow 1}}}\Big(\frac{\sum_{j'}{\bar\nu^{y,j'}_{ik}}}{m}\Big)\Bigg]-\frac{nml}{\sum_{i',j',k'}{t^{x,z,i', k'}_{j' \rightarrow 1}}}\Big(\frac{\sum_{j'}{\bar\nu^{y,j'}_{ik}}}{m}\Big)\Bigg)\\
&\quad +\frac{1}{\alpha}\Bigg(\avg\Bigg[\frac{nml}{\sum_{i',j',k'}{t^{x,y,i', j'}_{k' \rightarrow 1}}}\Big(\frac{\sum_{k'}{\bar\nu^{z,k'}_{ij}}}{l}\Big)\Bigg]-\frac{nml}{\sum_{i',j',k'}{t^{x,y,i', j'}_{k' \rightarrow 1}}}\Big(\frac{\sum_{k'}{\bar\nu^{z,k'}_{ij}}}{l}\Big)\Bigg) > -\epsilon\Bigg\}\Bigg]\\
& \geq  1- \sum_{(i,j,k) \in \T \cap \P} {\prob\Bigg[\Bigg|{\frac{n_{\bar x}}{\sum_{i'}{t^{y,z,j,k}_{i' \rightarrow 1}}}}-\avg\Bigg[{\frac{n_{\bar x}}{\sum_{i'}{t^{y,z,j,k}_{i' \rightarrow 1}}}}\Bigg]\Bigg| > \frac{\epsilon}{3}\Bigg]}\\
&\quad - \sum_{(i,j,k) \in \T \cap \P}{\prob\Bigg[\Bigg|{\frac{n_{\bar y}}{\sum_{j'}{t^{x,z,i,k}_{j' \rightarrow 1}}}}-\avg\Bigg[{\frac{n_{\bar y}}{\sum_{j'}{t^{x,z,i,k}_{j' \rightarrow 1}}}}\Bigg]\Bigg|>\frac{\epsilon}{3} \Bigg]}\\
&\quad - \sum_{(i,j,k) \in \T \cap \P}{\prob\Bigg[ \Bigg|{\frac{n_{\bar z}}{\sum_{k'}{t^{x,y,i,j}_{k' \rightarrow 1}}}}-\avg\Bigg[{\frac{n_{\bar z}}{\sum_{k'}{t^{x,y,i,j}_{k' \rightarrow 1}}}}\Bigg]\Bigg| > \frac{\epsilon}{3}\Bigg]}\\
&\quad - \sum_{(i,j,k) \in \T \cap \P}{\prob\Bigg[ \Bigg|\frac{n_{\bar y} n_{\bar z}}{\sum_{j',k'}{t^{x,i}_{j'k'\rightarrow 1}}}-\avg\Big[\frac{n_{\bar y} n_{\bar z}}{\sum_{j',k'}{t^{x,i}_{j'k'\rightarrow 1}}}\Big]\Bigg| > \frac{\epsilon}{3}\Bigg]}\\
&\quad - \sum_{(i,j,k) \in \T \cap \P}{\prob\Bigg[\Bigg|\frac{nml}{\sum_{i',j',k'}{t^{y,z,j', k'}_{i' \rightarrow 1}}} \Big(\frac{\sum_{i'}{\bar\nu^{x,i'}_{jk}}}{n} \Big)-\avg\Bigg[\frac{nml}{\sum_{i',j',k'}{t^{y,z,j', k'}_{i' \rightarrow 1}}} \Big(\frac{\sum_{i'}{\bar\nu^{x,i'}_{jk}}}{n} \Big)\Bigg]\Bigg| > \frac{\alpha \epsilon}{9} \Bigg]}\\
&\quad - \sum_{(i,j,k) \in \T \cap \P}{\prob\Bigg[\Bigg|\frac{nml}{\sum_{i',j',k'}{t^{x,z,i', k'}_{j' \rightarrow 1}}}\Big(\frac{\sum_{j'}{\bar\nu^{y,j'}_{ik}}}{m}\Big)-\avg\Bigg[\frac{nml}{\sum_{i',j',k'}{t^{x,z,i', k'}_{j' \rightarrow 1}}}\Big(\frac{\sum_{j'}{\bar\nu^{y,j'}_{ik}}}{m}\Big)\Bigg]\Bigg| > \frac{\alpha \epsilon}{9} \Bigg]}\\
&\quad - \sum_{(i,j,k) \in \T \cap \P}{\prob\Bigg[\Bigg|\frac{nml}{\sum_{i',j',k'}{t^{x,y,i', j'}_{k' \rightarrow 1}}}\Big(\frac{\sum_{k'}{\bar\nu^{z,k'}_{ij}}}{l}\Big)-\avg\Bigg[\frac{nml}{\sum_{i',j',k'}{t^{x,y,i', j'}_{k' \rightarrow 1}}}\Big(\frac{\sum_{k'}{\bar\nu^{z,k'}_{ij}}}{l}\Big)\Bigg]\Bigg| > \frac{\alpha \epsilon}{9} \Bigg]}\\
&\geq \ 1- 2nml \exp\Big(-2n r_{\bar x} (g(\frac{\epsilon}{3}-\Delta_a))^2\Big)
\quad -2nml \exp\Big(-2m r_{\bar y} (g(\frac{\epsilon}{3}-\Delta_b))^2\Big)\\
&\quad -2nml \exp\Big(-2l r_{\bar z} (g(\frac{\epsilon}{3}-\Delta_c))^2\Big)
-2nml \exp\Big(-2ml r_{\bar y} r_{\bar z} (g(\frac{\epsilon}{3}-\Delta_d))^2\Big)\\
&\quad -2nml \exp\Big(-2nml \Big(g\Big(\frac{\alpha\epsilon/9-\Delta'_a}{\chi_1}\Big)\Big)^2\Big) -2nml\exp\Big(-\frac{2n\epsilon^2}{\chi_2}\Big)\\
&\quad -2nml \exp\Big(-2nml \Big(g\Big(\frac{\alpha\epsilon/9-\Delta'_b}{{\chi_3}}\Big)\Big)^2\Big) -2nml\exp\Big(-\frac{2m\epsilon^2}{\chi_4}\Big)\\
&\quad -2nml \exp\Big(-2nml \Big(g\Big(\frac{\alpha\epsilon/9-\Delta'_c}{{\chi_5}}\Big)\Big)^2\Big)-2nml \exp\Big(-\frac{2l\epsilon^2}{\chi_6}\Big),
\end{align*}
where we define $\Delta_a :=p^{n_{\bar x}}/(1-p)$,
$\Delta_b :=p^{n_{\bar y}}/(1-p)$, $\Delta_c :=p^{n_{\bar z}}/(1-p)$, $\Delta_d :=p^{n_{\bar y} n_{\bar z}}/(1-p)$, $\Delta'_a:= 1/r_{\bar x} p^{mln_{\bar x}}/(1-p)$,
$\Delta'_b:= 1/r_{\bar y} p^{mln_{\bar y}}/(1-p)$, $\Delta'_c:= 1/r_{\bar z} p^{mln_{\bar z}}/(1-p)$ and $\chi_i$, $i \in \{1,\ldots,6\}$ are positive constants as defined in the proof of Step~\ref{om4}.
The first inequality follows since $\avg[\kappa_{ijk}] \geq \epsilon$, the second inequality follows from set inclusion and taking the union bound, and the fourth inequality follows from Steps~\ref{om2} and~\ref{om4}.
The proof then follows since as $n,m,l \rightarrow \infty$, we have $\Delta_a,\Delta_b,\Delta_c, \Delta_d,\Delta'_a, \Delta'_b, \Delta'_c \rightarrow 0$, $(g(\cdot))^2$ is a positive constant,
and the limit assumptions in the theorem imply that, as $n, m, l \to \infty$, we have $nml \exp (-n)$, $nml \exp (-m)$, $nml \exp (-l)$, $nml \exp (-ml)$, $nml \exp (-nml) \rightarrow 0$.
\end{cpf}
\end{prfc}

\section{Facets of the multilinear polytope of rank-one BTF}
\label{sec facets}


\begin{prfc}[of \cref{prop facets}]
Facetness of inequalities $w_{ijk} \geq 0$ for all $(i,j,k) \in [n]\times [m] \times [l]$ follows from Proposition~2 of~\cite{dPKha16}. To prove the facetness of the remaining inequalities, we employ the following standard strategy: denote by $g_1$ an inequality defining the feasible region of the flower LP. Consider a nontrivial valid inequality $g_2$ for $\MP_{G^{\rm BT}}$ that is satisfied tightly at all points in $\S_{G^{\rm BT}}$ that are binding for $g_1$. We then show that $g_1$ and $g_2$
coincide up to a positive scaling which by full-dimensionality of
the multilinear polytope (see Proposition~1 in~\cite{dPKha16}) implies $g_1$
defines a facet of $\MP_{G^{\rm BT}}$. In the following, we consider various points in $\S_{G^{\rm BT}}$ that are binding for $g_1$. For brevity, we refer to any such point as a binding feasible point (BFP).
It then suffices to consider the following inequalities:

\paragraph{$w_{111} \leq x_1$:} Let
\begin{equation}\label{vineq}
\sum_{i}{a_i x_i} + \sum_{j}{b_j y_j} + \sum_{k}{c_k z_k} + \sum_{i,j,k}{d_{ijk} w_{ijk}} \leq \alpha,
\end{equation}
be a nontrivial valid inequality for $\MP_{G^{\rm BT}}$ that is satisfied tightly at all points in $\S_{G^{\rm BT}}$ that are binding for $w_{111} \leq x_1$.
First, consider a BFP with $x=y=z = 0$.
Substituting this point in~\eqref{vineq} we obtain
\begin{equation}\label{bo0}
\alpha = 0.
\end{equation}
Next consider a BFP where all but one element in $(x,y,z)$ are zero, assuming that one component is different from $x_1$. Substituting such a point in~\eqref{vineq} and using~\eqref{bo0} we obtain
\begin{equation}\label{bo1}
a_i = b_j = c_k = 0, \qquad \forall i \in [n] \setminus \{1\},  j \in [m], k \in [l].
\end{equation}
Now, consider a BFP with $x_1 = y_1 = z_1 = 1$, and $x_i = y_j = z_k = 0$ for $i \in [n] \setminus \{1\}$, $j \in [m] \setminus \{1\}$, $k \in [l] \setminus \{1\}$. Substituting in~\eqref{vineq} and using~\eqref{bo0} and~\eqref{bo1}, we obtain
\begin{equation}
\label{adt1}
a_1 + d_{111} = 0.
\end{equation}
Consider a BFP with $x_{\tilde i} = y_{\tilde j} = z_{\tilde k} = 1$ for some $\tilde i \in [i] \setminus \{1\}$, $\tilde j \in [m]$, $\tilde k \in [l]$.  Substituting in~\eqref{vineq} and using~\eqref{bo0} and~\eqref{bo1} we obtain
\begin{equation}
\label{adt2}
d_{ijk} = 0, \qquad \forall i \in [n] \setminus \{1\},  j \in [m], k \in [l].
\end{equation}
Next consider a BFP with $x_1=y_1=z_1=  y_{\tilde j} =1$
for some $\tilde j \in [m] \setminus \{1\}$. Substituting in~\eqref{vineq}
and using~\eqref{adt1} gives $d_{1\tilde j 1} = 0$. Using a similar line of argument we conclude that
\begin{equation}\label{bo5}
d_{1j1} = d_{11k} =0, \qquad \forall j \in [m] \setminus \{1\}, k \in [l] \setminus \{1\}.
\end{equation}
Finally consider a BFP with $x_1=y_1=z_1=  y_{\tilde j} =z_{\tilde k}=1$
for some $\tilde j \in [m] \setminus \{1\}$, $\tilde k \in [l] \setminus \{1\}$. Substituting in~\eqref{vineq} and using~\eqref{adt1} and~\eqref{bo5} gives $d_{1\tilde j \tilde k} = 0$. Hence, we have
\begin{equation}\label{bo6}
d_{1jk} = 0, \qquad \forall j \in [m] \setminus \{1\}, k \in [l] \setminus \{1\}.
\end{equation}
From~\eqref{bo0}-\eqref{bo6} it follows that inequality~\eqref{vineq} is of the form $\beta w_{111} \leq \beta x_1$ for some $\beta > 0$, implying
$w_{111} \leq x_1$ defines a facet of $\MP_{G^{BT}}$.

\paragraph{$w_{111} \geq x_1 + y_1 + z_1 -2$:} Let~\eqref{vineq} be a nontrivial valid inequality for $\MP_{G^{\rm BT}}$ that is satisfied tightly at all points in $\S_{G^{\rm BT}}$ that are binding for $w_{111} \geq x_1 + y_1 + z_1 -2$. First, consider a BFP with $x_1 = y_1 = 1$ (resp. $x_1 = z_1 = 1$ and $y_1 = z_1 = 1$) and $x_i = y_j = z_k = 0$, otherwise. Substituting in~\eqref{vineq} yields:
\begin{equation}\label{bo7}
    a_1 + b_1 = a_1 + c_1 = b_1 + c_1 = \alpha.
\end{equation}
Next, consider a BFP with $x_1 = y_1 = z_1 = 1$ and $x_i = y_j = z_k = 0$, otherwise. Substituting in~\eqref{vineq} and using~\eqref{bo7} yields:
\begin{equation}\label{bo8}
    a_1 = b_1 = c_1 = -d_{111} = \frac{\alpha}{2}.
\end{equation}
Consider a BFP with $x_1 = y_1 = x_{\tilde i} = 1$
for some $\tilde i \in [n] \setminus \{1\}$ and $x_i = y_j = z_k = 0$, otherwise. Substituting in~\eqref{vineq} and
using~\eqref{bo8} gives $a_1 + b_1 + a_{\tilde i} = \alpha + a_{\tilde i} = \alpha$, implying $a_{\tilde i} = 0$. Using a similar line of arguments we get
\begin{equation}\label{bo9}
    a_i = b_j = c_k = 0, \qquad \forall i \in [n] \setminus \{1\}, \; j \in [m] \setminus \{1\}, \; k \in [l] \setminus \{1\}.
\end{equation}
Next consider a BFP with $x_1=y_1=z_1=  y_{\tilde j} =1$
for some $\tilde j \in [m] \setminus \{1\}$. Substituting in~\eqref{vineq} and using~\eqref{bo8} and~\eqref{bo9} gives $a_1+b_1+c_1+d_{111}+d_{1\tilde j 1} = \alpha + d_{1\tilde j 1} = \alpha$, implying that $d_{1\tilde j 1} = 0$. Using a similar line of argument yields~\eqref{bo5}.

Finally consider a BFP with $x_1=y_1=z_1=  y_{\tilde j} =z_{\tilde k}=1$
for some $\tilde j \in [m] \setminus \{1\}$ and $\tilde k \in [l] \setminus \{1\}$. Substituting in~\eqref{vineq} and using~\eqref{bo5},~\eqref{bo8}, and~\eqref{bo9} gives $a_1+b_1+c_1+d_{111}+d_{1\tilde j \tilde k} = \alpha+d_{1\tilde j \tilde k} = \alpha$, implying~\eqref{bo6}. Therefore, by~\eqref{bo5}-\eqref{bo6} and~\eqref{bo7}--\eqref{bo9}, we conclude that~\eqref{vineq} can be written as $\alpha (x_1 + y_1 + z_1 -w_{111}) \leq 2\alpha $ for some $\alpha > 0$ and this completes the proof.

\paragraph{$w_{211}-w_{111} \leq 1 - x_1$:} Let~\eqref{vineq} be a nontrivial valid inequality for $\MP_{G^{\rm BT}}$ that is satisfied tightly at all points in $\S_{G^{\rm BT}}$ that are binding for $w_{211}-w_{111} \leq 1 - x_1$.
Consider a BFP with $x_1 = 1$, and $x_i=y_j=z_k=0$, otherwise.  Substituting this point in~\eqref{vineq} gives
\begin{equation} \label{bo10}
a_1 = \alpha.
\end{equation}
Next, consider a BFP with $x_1 = x_{\tilde i} = 1$, for some $\tilde i \in [n] \setminus \{1\}$, and $x_i=y_j=z_k=0$, otherwise.  Substituting this point in~\eqref{vineq} and using~\eqref{bo10} gives $a_1 + a_{\tilde i} = \alpha + a_{\tilde i} = \alpha$, implying $a_{\tilde i} = 0$.
Using a similar line of arguments we obtain:
\begin{equation} \label{bo11}
a_i = b_j = c_k =0, \qquad \forall i \in [n]\setminus \{1\}, j \in [m], k \in [l].
\end{equation}
Consider a BFP with $x_2 = y_1 = z_1 = 1$, and $x_i=y_j=z_k=0$, otherwise.  Substituting in~\eqref{vineq} and using~\eqref{bo11} yields
\begin{equation} \label{bo12}
d_{211} = \alpha.
\end{equation}
Consider a BFP with $x_1 = x_2 = y_1 = z_1 = 1$, and $x_i=y_j=z_k=0$, otherwise.  Substituting this point in~\eqref{vineq} and using~\eqref{bo10}-\eqref{bo12}  yields
\begin{equation} \label{adt3}
d_{111} = -\alpha.
\end{equation}
Consider a BFP with $x_1 = y_{\tilde j} = z_{\tilde k} = 1$, for some $(\tilde j, \tilde k) \in [m]\times [l] \setminus \{(1,1)\}$
and $x_i=y_j=z_k=0$, otherwise.
Substituting this point in~\eqref{vineq} and using~\eqref{bo10} and~\eqref{bo11} gives
\begin{equation} \label{bo15}
d_{1jk} = 0, \qquad \forall (j, k) \in [m]\times [l] \setminus \{(1,1)\}.
\end{equation}
Consider a BFP with $x_1 = x_{\tilde i}= y_{\tilde j} = z_{\tilde k} = 1$, for some $\tilde i \in [n] \setminus \{1\}$ and $(\tilde j, \tilde k) \in [m]\times [l] \setminus \{(1,1)\}$
and $x_i=y_j=z_k=0$, otherwise.
Substituting in~\eqref{vineq} and using~\eqref{bo10},~\eqref{bo11}, and~\eqref{bo15} gives
\begin{equation}
\label{adt4}
d_{ijk} = 0, \qquad \forall i \in [n] \setminus \{1\}, (j, k) \in [m]\times [l] \setminus \{(1,1)\}.
\end{equation}
Finally, consider a BFP with $x_2 = x_{\tilde i}= y_1= z_1 = 1$, for some $\tilde i \in [n] \setminus \{1,2\}$ and $x_i=y_j=z_k=0$, otherwise.
Substituting in~\eqref{vineq} and using~\eqref{bo11} and~\eqref{bo12} yields
\begin{equation} \label{bo17}
d_{i11} = 0, \qquad \forall i \in [n] \setminus \{1,2\}.
\end{equation}
Therefore, from~\eqref{bo10}-\eqref{bo17} it follows that inequality~\eqref{vineq} can be equivalently written as $\alpha (x_1 + w_{211} -w_{111}) \leq \alpha $ for some $\alpha > 0$ and this completes the proof.
\end{prfc}


\section*{Acknowledgments}
A.~Del~Pia is partially funded by AFOSR grant FA9550-23-1-0433. A.~Khajavirad is partially funded by AFOSR grant FA9550-23-1-0123. Any opinions, findings, and conclusions or
recommendations expressed in this material are those of the authors and do not necessarily reflect the views of the Air Force Office of Scientific Research.

\ifthenelse{\boolean{MOR}}
{
     \bibliographystyle{informs2014}
     \bibliography{biblio}
}
{
\bibliographystyle{plainurl}

}

\end{document}